\documentclass[stsy]{informs3_noinforms} 

\OneAndAHalfSpacedXII 

\usepackage{hyperref}

\usepackage{natbib}
 \bibpunct[, ]{(}{)}{,}{a}{}{,}%
 \def\bibfont{\small}%
 \usepackage{bbm}
\usepackage{array} 
\newcolumntype{L}[1]{>{\raggedright\arraybackslash}p{#1}}
\newcolumntype{C}[1]{>{\centering\arraybackslash}p{#1}}

\usepackage{multirow}
\usepackage{booktabs}
\usepackage[table]{xcolor}
\usepackage{makecell}
\usepackage{nicematrix}
\definecolor{lightgray}{RGB}{230,230,230}
\definecolor{celestialblue}{rgb}{0.29, 0.59, 0.82}
\usepackage{booktabs}
\usepackage{colortbl}
\usepackage{hhline}
\usepackage{arydshln}

\definecolor{columbiablue}{rgb}{0.61, 0.87, 1.0}

\TheoremsNumberedThrough     
\ECRepeatTheorems

\EquationsNumberedThrough    

\MANUSCRIPTNO{}

  \newtheorem{assump}{Assumption}
 \newtheorem{observation}{Observation}
\begin{document}



\RUNTITLE{Polynomial time algorithm for optimal stopping}


\TITLE{Polynomial time algorithm for optimal stopping with fixed accuracy}

\ARTICLEAUTHORS{%
\AUTHOR{Yilun Chen }
\AFF{The School of Data Science, CUHK-Shenzhen,  \EMAIL{chenyilun@cuhk.edu.cn}}
\AUTHOR{David A. Goldberg}
\AFF{The School of Operations Research and Information Engineering, Cornell University, \EMAIL{dag369@cornell.edu}}
} 
\ABSTRACT{The problem of high-dimensional path-dependent optimal stopping (OS) is important to multiple academic communities and applications.  Modern OS tasks often have a large number of decision epochs, and complicated non-Markovian dynamics, making them especially challenging.  Standard approaches, often relying on ADP, duality, deep learning and other heuristics, have shows very strong empirical performance, yet have limited rigorous guarantees (which may scale exponentially in the problem parameters and/or require previous knowledge of basis functions or additional continuity assumptions).  Although past work has placed these problems in the framework of computational complexity and polynomial-time approximability, those analyses were limited to simple one-dimensional problems.
\\\indent For long-horizon complex OS problems, is a polynomial time solution even theoretically
possible?  We prove that given access to an efficient simulator of the underlying information process, and fixed accuracy $\epsilon$, there exists an algorithm that returns an $\epsilon$-optimal solution (both stopping policies and approximate optimal values) with computational complexity scaling polynomially in the time horizon $T$ and underlying dimension $D$ (in our computational model).   Like the first polynomial-time (approximation) algorithms for several other well-studied problems, our theoretical guarantees are polynomial yet impractical (with a runtime scaling exponentially in $\epsilon^{-2}$).  Our approach is based on a novel expansion for the optimal value which may be of independent interest.
}
\maketitle

\section{Introduction}
\subsection{Overview of problem}
We consider the problem of discrete-time optimal stopping (OS): a decision maker (DM) observes a stochastically evolving information process. At each time period, an irrevocable decision of whether or not to stop must be made by the DM. A stop decision results in a reward/cost, which depends on the underlying stochastic process. The DM's goal is to find the strategy that maximizes the expected reward/minimizes the expected cost. Such OS problems find an array of applications across areas including management science (MS) / operations management (OM) / operations research (OR), finance, economics, statistics and computer science (CS).
\\\indent Real-world stopping tasks are becoming increasingly challenging, with complex problem structure including long horizon and non-Markovian/path-dependent dynamics.  These features call for new methods and tools. We provide a few examples below.
\begin{itemize}
    \item \textbf{Statin Choice Decisions.\ } The decision of initiating statin therapies for diabetes patients can be modelled as an OS problem (\cite{denton2009optimizing}).  Here the underlying information is the evolving state of each patient, including their gender, height, weight, smoking status, blood pressure, cholesterol, HbA1c, and history of stroke and CHD. The decision window of such problems spans the patient's entire lifetime, with an abundance of decision epochs. The underlying state is both high-dimensional and with possibly complicated, history-dependent transition dynamics.
    
    \item \textbf{Seizure Detection.\ } OS is also used to model and solve the detection of epileptic seizures (\cite{santaniello2012optimal, santaniello2011quickest}). Here the decision relies on the electroencephalogram(EEG) data stream, which is generated by multiple electrode ``sensors'' (in some cases more than 20) and evolves in a non-stationary manner. Seizure detection needs to be prompt, which results in a large number of decision epochs. For example, in the model proposed by \cite{santaniello2012optimal} the number of decision epochs is typically on the order of $1000.$
    
    \item \textbf{Non-Markovian Bermudan Option Pricing.\ } An important application of OS is the pricing of Bermudan options.  In recent years, there has been a growing interest in capturing the long-range or short-range dependency in assets' volatility.  The OS problem associated to the options written on such assets has non-Markovian / path-dependent dynamics (\cite{leao2019discrete, sturt2021nonparametric}).
\end{itemize}

Such complicated structure can lead to severe challenges.  Here the scale of the OS problems ``\textit{grow especially large because of the need to incorporate some history-dependence into the state space}'' (Section 6.6. of \cite{boucherie2017markov}). Such problems are ``not easy to solve'' (\cite{leao2019discrete}), essentially because ``\textit{the enlarged state space will be high-dimensional when the OS problem has many periods, and the difficulty of solving a Markovian OS problem explodes in the dimensionality of the state space}'' (\cite{sturt2021nonparametric}). 

\indent Although a large number of heuristic approaches have been proposed to cope with these challenges, and many of them do exhibit excellent numerical performance on certain test benchmarks, much less is known regarding the theoretical computational complexity of high-dimensional OS. \cite{chalasani1999approximate} was the first work to put OS (and options pricing) in the formal framework of computational complexity theory and polynomial-time approximability as commonly defined in the OR/CS communities (\cite{arora2009computational,williamson2011design}).  The authors first prove an impossibility result, showing that under a standard complexity-theoretic assumption, no deterministic algorithm can exactly price even very simple one-dimensional options in polynomial time.  The same work, as well as later works \cite{aingworth2000accurate} and \cite{schalekamp2007some}, then provide several positive results.  These positive results show that for any fixed accuracy $\epsilon$, there are deterministic algorithms with polynomial runtime that can determine an $\epsilon$-approximation to the value of certain simple options.  However, none of these works address the question of whether complex high-dimensional OS problems can be approximated to fixed accuracy $\epsilon$ in polynomial time.
\\\indent In contrast to the few aforementioned complexity results, the rigorous performance analysis of most existing methods (e.g. ADP, pathwise optimization, decision-tree and deep neural network approximations) rely fundamentally on the underlying approximate basis functions/martingales at one's disposal, implicitly requiring them to be expressive enough to capture the continuation values/optimal martingale. However, there is no known polynomial-time method to generate such functions with provable guarantees in complicated, long-horizon, path-dependent instances.  This situation is analogous to that of other significant problems in OR/CS, such as the Euclidean Traveling Salesman Problem (TSP) and Linear Programming (LP), where many heuristics with strong empirical performance predate the first polynomial-time algorithms (see Appendix \ref{sec:ecptas} for a more detailed discussion).  In those cases the development of the first polynomial-time algorithms, although often themselves impractical, ultimately helped lead to new practical approaches (\cite{arora1996polynomial, kisfaludi2022gap, wright2005interior}).
\\\indent The above motivates us to study the theoretical computational complexity of high-dimensional complex OS, in particular whether these problems can be approximated to fixed accuracy $\epsilon$ in polynomial time.  In this paper, we propose such an algorithm (under appropriate assumptions).  We next highlight our main contributions in detail, followed by a discussion of the limitations of our methodology and an additional review of relevant literature, as well as an outline of the organization of the rest of the paper.
\subsection{Main Contributions}
\subsubsection{Polynomial time algorithm for OS with fixed accuracy.} We prove that given access to an efficient simulator of the underlying information process, and fixed accuracy $\epsilon$, our algorithm returns an $\epsilon$-optimal solution (both stopping policies and approximate optimal values) with computational complexity scaling polynomially in the time horizon $T$ and underlying dimension $D$ (under appropriate assumptions on the rewards).  We note that in many cases of interest (e.g. under the Black-Scholes model) the underlying information process will simply be a high-dimensional geometric Brownian motion (B.m.), for which such a simulation is straightforward.  
\\\indent More precisely, our theory requires that either $(1)$ reward functions are uniformly bounded (independent of T and D) or $(2)$ the reward process satisfies a mild concentration / moment condition.  Under either of these assumptions, and with no additional assumptions of e.g. independence, continuity, or Markovian dynamics, our proposed algorithm can implement $\epsilon-$approximate OS policies (in expectation) and estimates of the optimal value (with high probability), with a computational complexity scaling polynomially in the number of decision epochs and linearly in the unit simulation and reward evaluation costs, for any fixed accuracy $\epsilon \in (0, 1)$. In cases where the unit simulation and reward evaluation costs can be appropriately bounded (in terms of e.g. the horizon and/or underlying dimension), our algorithm achieves the desired approximation guarantee with computational time scaling polynomially in the time horizon and underlying dimension.  Thus in theory one can do (exponentially) better than either DP (whose runtime depends on the number of states, which may be exponentially large) or naive fully nested simulation (whose runtime scales exponentially in the time horizon), although due to the impractical constants appearing in our polynomial complexity bounds the time horizon and number of states have to be enormous for our method to exhibit a provably superior performance.
\\\indent Furthermore, we prove that the Bermudan max call (a standard benchmark in option pricing) satisfies the assumptions required for our methodology, and thus our algorithmic results are applicable in this setting.  The precise computational model we use for our analysis falls within the framework of information complexity (in line with past work on related problems), and we refer the reader to Section\ \ref{complexsec1} for a detailed discussion.
\subsubsection{Novel expansion theory.}  Our algorithm is based on a novel expansion representation of the optimal value.  This expansion is by nature different from both DP and naive fully nested simulation (which rely on backward inductive logic). In particular, the first term of our expansion coincides with the expected pathwise minimum/maximum (i.e. the hindsight optimal value). The subsequent terms are derived in a recursive manner, themselves expectations of certain minima/maxima conceptually corresponding to various notions of regret (and regret of regret, etc.).  Furthermore, by truncating our expansion after $O(\frac{1}{\epsilon})$ terms, one can achieve an $\epsilon$-approximation to the OS value with a level of nesting (of conditional expectations) scaling as $O(\frac{1}{\epsilon})$ (independent of the time horizon).  We believe this novel representation contributes to the theoretical understanding of optimal stopping and may be of independent interest.
\subsection{Limitations of our methodology and additional literature review}\label{limitssec}
To add additional context for our results, here we discuss several limitations of our methodology.  This includes: (1) the impractical theoretical complexity; (2) the limitations of our fixed $\epsilon$ framework; (3) the restrictions imposed by our technical assumptions.  We also provide an expanded discussion of relevant literature (to further provide additional context).
\subsubsection{Impractical theoretical complexity.} The theoretical runtime of our algorithm (to achieve error $\epsilon$) is $\exp\big( O(\epsilon^{-2}) \big) T^{O(\frac{1}{\epsilon})}$ in the case that rewards are uniformly bounded (independent of T and D), and $\exp\big( O( \gamma_0^8 \epsilon^{-6}) \big) T^{O(\gamma_0^{4} \epsilon^{-3})}$ in the case of unbounded rewards for which the squared coefficient of variation (s.c.v.) of the maximum reward is bounded by $\gamma_0 - 1$ (independent of $T$ and $D$).  The O-notation above hides very large absolute constants, making our theoretical guarantees impractical, especially in settings where one is interested in very accurate solutions (which includes many problems in financial engineering).  
\subsubsection{Fixed $\epsilon$ framework.} The notion of polynomial-time efficiency for any fixed accuracy $\epsilon$ is not standard in the study of OS in mathematical finance.  Alternative notions of efficiency (e.g. that of semi-tractability introduced in \cite{belomestny2019semi}) will in general be more appropriate for many problems in finance and options pricing, in which either: (1) the problem arises as a discretization of a continuous problem (in which case the horizon length cannot be treated as being conceptually separate from the desired accuracy); or (2) a very accurate estimate is needed.
\\\indent \emph{Additional motivation for the fixed $\epsilon$ framework.} This notion of efficiency coincides with that of Polynomial Time Approximation Scheme (PTAS), and is standard in the study of approximation algorithms in OR/CS, and consistent with past work on the theoretical computational complexity of options pricing (\cite{chalasani1999approximate, aingworth2000accurate}).  It has also recently been applied to the study of several problems in OM/MS (\cite{alaei2022revenue, aouad2021display, derakhshan2022product}).  Furthermore, there are several application domains of OS in which the problems do not arise as discretizations of continuous problems, and may not require very accurate estimates.  For example, in certain healthcare applications the time horizon may arise from either the limitations of certain medical technologies or common medical practices (e.g. annual screenings, see \cite{zhang2012optimization, denton2009optimizing, santaniello2012optimal}).  Similarly, the pricing of certain Bermudan options (e.g. Bermudan swaptions) may have annual, semi-annual, or monthly exercise dates (\cite{beinker2017accurate}).  Regarding problems which may not require very accurate estimates, it has been noted in past works on pricing Bermudan options that a less accurate solution may suffice when the associated bid-ask spread is large (see \cite{hin2015analysis,hamida2005recovering,kou2004option,
longstaff2001valuing,glasserman2003cap}).
\subsubsection{Assumptions required for theoretical analysis.} Our main theoretical results require either rewards to be uniformly bounded (independent of $T$ and $D$), or for the s.c.v. of the maximum reward to be bounded (independent of $T$ and $D$), as well as access to an efficient simulator of the underlying information process.  There may be many applications for which the relevant OS problems would not satisfy these regularity conditions, or for which the information process may be difficult to simulate.
\\\indent \emph{Additional motivation for these assumptions.} The boundedness assumption may be relevant to many applications outside of finance, including those building on the original ``secretary problem" in which the goal is to maximize a suitable probability (see e.g \cite{freeman1983secretary,emmerling2021accept,krotofil2014vulnerabilities,langrene2017field}).  The alternate assumption regarding the s.c.v. of the maximum arises in the unbounded setting due to our analysis of the error induced by truncating the rewards (which requires an application of the Cauchy-Schwarz inequality, see e.g. the proof of Lemma\ \ref{maxb1} in Section\ \ref{sec:ecmaximization}).  Furthermore, we require that this error is not too much larger than the optimal value of the stopping problem, and we show in Lemma\ \ref{maxb2} that this will hold as long as the s.c.v. of the maximum reward is not too large.  This assumption is also naturally motivated by the i.i.d. setting (i.e. secretary problem), since for many unbounded distributions the associated i.i.d. sequence of random variables (r.v.s) $\lbrace X_1,\ldots,X_T \rbrace$ will have the property that the s.c.v. of $\max(X_1,\ldots,X_T)$ is bounded independent of $T$, making our results applicable.  In particular, by well-known results in extreme value theory, this includes any distribution belonging to the Gumbel domain of attraction, such as the exponential distribution (\cite{embrechts2013modelling}).  Let us also note that even when these assumptions do not hold, our methodology can still be implemented as a heuristic.
\\\indent Regarding our assumption about access to an efficient simulator, we note that in many discrete-time options pricing problems, the information process will be distributed as the discretization of a multi-dimensional geometric Brownian motion (B.m).  This will be the case (for example) under the standard assumptions of the Black-Scholes model.  In any such model, one can efficiently implement the required simulator by generating trajectories of a multi-dimensional geometric B.m. conditioned on the values observed up to time $t$.  In an even broader family of problems (including those with long range dependencies, see e.g. \cite{bayer2016pricing}), one can take the information process to be a multi-dimensional Gaussian random vector, and thus one can implement the required simulator in time polynomial in $T$ and $D$ by using standard linear-algebra methods (\cite{sabino2007monte}).  Furthermore, in many applications the information process will evolve as a high-dimensional Markov chain.  In any such model, the complexity of implementing the required simulator is no greater than that of generating unconditioned trajectories of the information process, where the ability to implement such a simulation is a common assumption across the OS literature.
\\\indent Let us also note that the well-studied Bermudan max call indeed satisfies the assumptions required for our theoretical analysis, as we show in Section\ \ref{sec:maxcall} (and Section\ \ref{sec:ecmaxcall} of the Appendix).
\subsubsection{Additional literature review.}\label{litsec}
\indent There is a vast literature on computational OS, and one of the most popular approaches is regression/ADP.  Here one fixes a family of basis functions to approximate the DP value functions.  Seminal papers in this area include \cite{tsitsiklis2001regression} and \cite{longstaff2001valuing}.  There was much subsequent work on e.g. nonparametric regression, smoothing spline regression, recursive kernel regression, integer programming, reinforcement learning, multilevel Monte Carlo, and (deep) neural networks (see \cite{lai2004valuation}, \cite{egloff2005monte}, \cite{kohler2008regression}, \cite{belomestny2015multilevel}, \cite{belomestny2018optimal}, \cite{bayer2018dynamic}, \cite{becker2019solving}, \cite{ciocan2020interpretable}, \cite{sturt2021nonparametric}, \cite{zhou2021unbiased}). These methods mostly focus on achieving empirical success rather than ensuring theoretical guarantees, with their performance relying on how well the choice of basis functions can approximate the true value functions.  Those works with theoretical guarantees often require additional continuity assumptions, and/or scale exponentially in the time horizon.  Most existing theoretical analyses for these methods focus on their convergence to the best approximation within the fixed set of basis functions (\cite{clement2002analysis}, \cite{glasserman2004number}, \cite{bouchard2012monte}, \cite{bezerra2018discrete}, \cite{belomestny2011rates}, \cite{stentoft2004convergence}), and do not have guarantees comparable to our own.
\\\indent Building on the seminal work of \cite{davis1994deterministic}, significant progress was made simultaneously by \cite{haugh2004pricing} and \cite{rogers2002monte} in their formulation of a dual methodology for OS. Instead of finding an OS time, it formulates and solves a dual optimization problem over the set of adapted martingales. Other dual representations were subsequently discovered (e.g. \cite{jamshidian2007duality}), and the methodology was extended to more general control problems (\cite{brown2010information}).  We note that closely related dual formulations and methodologies have also appeared in the study of multi-stage stochastic optimization (\cite{rockafellar1991scenarios}).  Simulation approaches via dual formulation have since led to substantial algorithmic progress (see \cite{andersen2004primal, chen2007additive,kolodko2006iterative,belomestny2013multilevel, ibanez2020recursive, kolodko2004efficient, belomestny2006monte, desai2012pathwise,christensen2014method, belomestny2019optimal,lelong2018dual}).  However, the quality of the bounds derived using these approaches typically depends on the expressiveness of the set of basis functions/martingales, or the quality of the initial estimation of the value functions.  In addition, those works which do provide rigorous approximation  guarantees typically require a level of nesting (of conditional expectations) scaling with the time horizon (see \cite{chen2007additive, kolodko2004efficient}).
\\\indent More recently, new machine learning techniques such as decision-trees and deep neural networks, and ideas from robust optimization were introduced to solve the problem of high-dimensional OS and American option pricing. Relevant papers include  \cite{kohler2010pricing}, \cite{ciocan2020interpretable}, \cite{becker2019deep}, \cite{becker2019solving}, \cite{sturt2021nonparametric}. These approaches exhibit excellent performance on a number of test benchmarks, but generally don't have strong theoretical performance guarantees independent of the quality of the basis functions/free from additional continuity assumptions. 

\subsection{Organization}
\indent The rest of the paper is organized as follows. In Section \ref{sec:definition}, we define a general OS problem and specify the assumptions required for our main results.  In Section \ref{sec:mainresults}, we describe the novel expansion (and associated error bounds) which provides the foundation of our approach.  In Section \ref{sec:algoresults}, we introduce our expansion-inspired algorithm, and state our main algorithmic results.  In Section \ref{sec:maxcall}, we show that our results imply a polynomial-time algorithm (in our computational model) for the Bermudan max call.  We present some closing thoughts and directions for future research in Section \ref{conc}.  We include a broader discussion and motivation of the notion of polynomial-time efficiency for a fixed accuracy $\epsilon$, additional intuition for our main results, several auxiliary proofs, as well as the formal descriptions and analysis of all algorithms, in Appendix\ \ref{sec:ecptas} - \ref{sec:ecnonneg}.
\section{Problem Formulation} \label{sec:definition}
\subsection{Formulation of problem and notation}
Let $\mathbf{Y} = \lbrace Y_t, t = 1,\ldots,T \rbrace$ be a (possibly high-dimensonal and non-Markovian) discrete-time stochastic process. For $t \in \lbrace 1,\ldots,T \rbrace,$ let $Y_{[t]} \stackrel{\Delta}{=} (Y_1,\ldots,Y_t)$ be the time $t$ partial trajectory, where we assume that (for all $t$) $Y_t$ is a $D$-dimensional vector for some $D \geq 1$.  Let $\lbrace g_t: \mathbb{R}^{D \times t} \to \mathbb{R} , t = 1,\ldots,T \rbrace$ be a sequence of general (possibly time-dependent) measurable reward functions mapping trajectories $Y_{[t]}$ to reals.  The OS problem is that of computing $\textsc{opt} \stackrel{\Delta}{=} \sup_{\tau \in {\mathcal T}} E[g_\tau(Y_{[\tau]})],$ with ${\mathcal T}$ the set of all stopping times adapted to the natural filtration ${\mathcal F} = \lbrace {\mathcal F}_t, t = 1,\ldots,T \rbrace$ generated by the process $\lbrace Y_t, t = 1,\ldots,T \rbrace$. For notational simplicity, we let $Z_t \stackrel{\Delta}{=} g_t(Y_{[t]}),$ where we write $Z_t(Y_{[t]})$ if we wish to make the dependence explicit, and denote the corresponding process by $\mathbf{Z}$. Throughout the paper we assume $Z_t$ is integrable and non-negative for all $t$, and note this is without loss of generality (w.l.o.g.), as we formally show in Appendix\ \ref{sec:ecnonneg} (a similar observation was recently made in \citet{baas2023sampling}).  The problem of interest is to compute $\textsc{opt} = \sup_{\tau \in {\mathcal T}} E[Z_{\tau}],$ which we assume to be strictly positive (to rule out the degenerate case in which all rewards are zero w.p.1). 
\subsection{Assumptions}
Our theoretical and algorithmic results require additional assumptions of two flavors: assumptions on the regularity of the reward process $\mathbf{Z}$; and assumptions on one's ability to efficiently simulate certain processes.  We formally prove that our theoretical guarantees hold under either of the following regularity assumptions: 

\begin{assump}\label{assumptionA1}
Assume $Z_t$ is uniformly bounded with $0 \leq Z_t \leq U$ w.p.1 for all $t$, with $U$ a constant independent of both the time horizon $T$ and the dimension $D$.
\end{assump}
\begin{assump}\label{assumptionA2}
Let $M \ \stackrel{\Delta}{=} \max_{t = 1,\ldots,T} Z_t.$  Assume the squared coefficient of variation of $M$, $\frac{E[M^2]}{(E[M])^2} - 1$, is upper bounded by a constant independent of both the time horizon $T$ and the dimension $D.$  We let $\gamma_0 \stackrel{\Delta}{=} \frac{E[M^2]}{(E[M])^2}$.
\end{assump}

Both of these assumptions impose a certain level of regularity of the underlying reward process $\mathbf{Z},$ although are quite different in nature from assumptions appearing in previous literature (typically requiring certain strong continuity properties or knowledge of certain basis functions).  Although here we focus on these two regularity assumptions, we note that our approach is flexible enough to yield similar results across a variety of related conditions in both the minimization and maximization settings.  As our expansion theory is most clearly and concisely explained under Assumption \ref{assumptionA1}, we build our theory in Section \ref{sec:mainresults} under Assumption \ref{assumptionA1}.  We then combine with some additional auxiliary results and a truncation argument to conduct our analysis under Assumption \ref{assumptionA2}.  

We next formalize our assumption about one's ability to efficiently simulate the relevant processes.  We note that efficient Monte Carlo algorithms (with rigorous guarantees for high-dimensional problems) typically require related assumptions (at least the ability to generate unconditioned trajectories).  

\begin{assump}\label{assumptionB}
Assume that for any $t$ and partial trajectory $\gamma$ (of the information process up to time $t$), generating one trajectory of the underlying process $Y_{[T]}$ (conditioned on the event $\lbrace Y_{[t]} = \gamma \rbrace$) takes $C$ units of time, and evaluating any given stopping reward (i.e. the function $g_t$) on any given (partial) trajectory takes $G$ units of time, both scaling (at most) polynomially in the time horizon $T$ and dimension $D.$
\end{assump}

\section{Theory}\label{sec:mainresults}

\subsection{Novel expansion representation for \textsc{OPT}}\label{sec:mainresultsexpansion}

\subsubsection{Simple intuition.}
\indent We begin by giving the simple intuition behind the expansion.  We wish to compute $\textsc{opt} = \sup_{\tau \in \mathcal{T}} E[Z_\tau].$ Trivially, $\sup_{\tau \in \mathcal{T}} E[Z_\tau] \leq E[\max_{i = 1,\ldots,T} Z_i].$ We next turn this straightforward bound into an equality by compensating with a remainder term, which is characterized by a new OS problem.  We introduce a specific martingale $\big\lbrace E\big[ \max_{i = 1,\ldots,T} Z_i | \mathcal{F}_t\big], \ t = 1,\ldots,T \big\rbrace,$ i.e. the Doob martingale of $\max_{i = 1,\ldots,T} Z_i,$ adapted to filtration $\mathcal{F}.$  By the optional stopping theorem, for any $\tau \in \mathcal{T}$, it holds that $E\big[E\big[ \max_{i = 1,\ldots,T} Z_i | \mathcal{F}_\tau\big]\big] = E[\max_{i = 1,\ldots,T } Z_i],$ and thus
\begin{equation*}
    E[Z_\tau] = E\bigg[\max_{i = 1,\ldots,T } Z_i\bigg] + E\bigg[  Z_\tau - E\big[ \max_{i = 1,\ldots,T} Z_i | \mathcal{F}_\tau\big]  \bigg].
\end{equation*}
Taking supremum from both sides, we conclude that 
\begin{equation*}
    \textsc{opt} = E\bigg[\max_{i = 1,\ldots,T} Z_i\bigg] + \sup_{\tau \in \mathcal{T}}E\bigg[Z_\tau - E\big[ \max_{i = 1,\ldots,T} Z_i | \mathcal{F}_\tau\big] \bigg].
\end{equation*}
There are many bounds for the gap $E\bigg[\max_{i = 1,\ldots,T} Z_i\bigg] - \textsc{opt}$ appearing in the literature on so-called prophet inequalities, see e.g. \cite{hill1983stop,hill1992survey}.  The novelty of our approach will be to ``recurse" this logic, as follows.  For $t \in \lbrace 1,\ldots,T \rbrace,$ let $Z^1_t \stackrel{\Delta}{=} Z_t$ and $Z^2_t \stackrel{\Delta}{=} Z^1_t - E\big[ \max_{i = 1,\ldots,T} Z^1_i | \mathcal{F}_t\big]$ (w.p.1. non-positive).  The above equation becomes
\begin{equation*}
    \textsc{opt} = E\bigg[\max_{i = 1,\ldots,T } Z^1_i\bigg] + \sup_{\tau \in \mathcal{T}}E\big[Z^2_\tau \big].
\end{equation*}
Now, we simply observe that we may recursively repeat this process on the problem $\sup_{\tau \in \mathcal{T}}E\big[Z^2_\tau \big]$, and then all subsequent problems.  This yields an explicit expansion for the optimal value which is amenable to simulation.
\subsubsection{Expansion for \textsc{OPT}.}
 \indent We now formalize the above intuition and provide our expansion for \textsc{opt}.  For $k \geq 1$ and $t \in \lbrace 1,\ldots,T \rbrace$, let $Z^{k+1}_t \stackrel{\Delta}{=} Z^k_t - E\big[ \max_{i = 1,\ldots,T} Z^k_i | {\mathcal F}_t \big]$ (denoting the corresponding stochastic process by $\mathbf{Z}^k$).  It follows from the basic properties of conditional expectation and a straightforward induction that $Z^k_t \leq 0$ for all $t \in \lbrace 1,\ldots,T \rbrace$ and $k \geq 2$.  Then our main result is as follows.
\begin{theorem}\label{thm:main1}
$\textsc{opt} = \sum_{k=1}^{\infty} E[\max_{t = 1,\ldots,T} Z^k_t].$
\end{theorem}
Let us more explicitly give the first few terms.  For $k \geq 1$, let $L_k \stackrel{\Delta}{=} E[\max_{t = 1,\ldots,T} Z^k_t]$.  Then $L_1,L_2$ are as follows. 
\[
L_1 = E\big[\max_{t = 1,\ldots,T} g_t(Y_{[t]}) \big]\ \ \ ;\ \ \ L_2 = E\bigg[\max_{t = 1,\ldots,T}\bigg( g_t(Y_{[t]}) - E\big[ \max_{i = 1,\ldots,T} g_i(Y_{[i]}) | {\mathcal F}_t \big] \bigg) \bigg].
\]
Note that the first term, $L_1$, is the only positive term in the expansion, which corresponds to the obvious upper bound.  Later terms are the expectations of explicit suprema with interpretations in terms of certain notions of (negative) ``higher order" regret, each of which can be computed by simulation.  We defer a rigorous proof of Theorem \ref{thm:main1} to Appendix \ref{sec:ectheoryproof}. 

\subsection{Approximation guarantees when truncating the expansion}
 The power of Theorem\ \ref{thm:main1} is that it allows for rigorous approximation guarantees when the infinite sum is truncated.  Let $E_k \stackrel{\Delta}{=} \sum_{i=1}^k L_i$.  

\begin{theorem}\label{thm:main2}
Suppose Assumption \ref{assumptionA1} holds with $U = 1.$ Then for all $k \geq 1$, $0 \leq E_k - \textsc{opt}  \leq \frac{1}{k}$.
\end{theorem}

Thus truncating our expansion after $k$ terms yields an absolute error at most $\frac{1}{k}$, regardless of the distribution of the underlying stochastic processes $\mathbf{Y}$ and $\mathbf{Z}$ or length of the time horizon, requiring only that the reward functions are bounded in $[0,1].$   Before proving Theorem\ \ref{thm:main2}, we first state some key properties of our expansion.  These follow essentially immediately from our earlier simple intuition from Section \ref{sec:mainresultsexpansion}, i.e. recursively applying the optional stopping theorem to the appropriate remainder term, combined with definitions and a straightforward induction.  For completeness, we include a proof in Appendix \ref{sec:ectheoryproof}. 
\begin{lemma}\label{lem:intuit1}
For all $k \geq 1$, $\textsc{opt} = \sum_{i=1}^k E[\max_{t = 1,\ldots,T} Z^i_t] + \sup_{\tau \in {\mathcal T}} E[Z^{k+1}_{\tau}].$
In addition: (1) w.p.1 $Z^k_t$ is non-positive and integrable for all $k \geq 2$ and $t \in \lbrace 1,\ldots,T \rbrace$; (2) for each $t \in \lbrace 1,\ldots,T \rbrace $, $\lbrace Z^k_t , k \geq 2 \rbrace$ is a monotone increasing sequence of random variables; and (3) for each $k \geq 1$, $\mathbf{Z}^k$ is adapted to ${\mathcal F}$.
\end{lemma}
Next, we state and prove a pathwise rate of convergence guarantee, essentially that after $k$ iterations of our expansion, the negative maximum of every sample path is at most $\frac{1}{k}$.  This is a much stronger result than that stated in Theorem\ \ref{thm:main2}, which only regards expectations.  This pathwise convergence enables us to construct provably accurate and polynomial time policies. 
\begin{lemma}\label{converge1}
Suppose that Assumption \ref{assumptionA1} holds with $U = 1.$   For $k \geq 2$, w.p.1, $ 0 \ \leq \ - \max_{t = 1,\ldots,T} Z^k_t \ \leq\ \frac{1}{k - 1}$.
\end{lemma}
\proof{Proof}
That $0 \leq - \max_{t = 1, \ldots, T} Z^k_t $ follows from Lemma \ref{lem:intuit1}.  It thus suffices to show the upper bound.  It will be helpful to first get a better handle on what happens in the final period $T$, since the conditional expectations degenerate in the last period.  In particular, for $k \geq 1$, $Z^{k + 1}_T = Z^k_T  - \max_{t = 1, \ldots, T} Z^k_t. $ Summing over $k$ then yields a telescoping sum, leading to the equation $Z^{k+1}_T = Z_T - \sum_{i=1}^k \max_{t = 1,\ldots,T} Z^i_t.$ By Lemma \ref{lem:intuit1}, $Z^{k+1}_T$ is non-positive, hence $Z_T - \sum_{i=1}^k \max_{t = 1,\ldots,T} Z^i_t \leq 0$\ w.p.1. We may rewrite this inequality as $\sum_{i=2}^{k} \max_{t = 1,\ldots,T} Z^i_t \geq Z_T - \max_{t = 1,\ldots,T} Z_t. $ Again by Lemma \ref{lem:intuit1}, we know that $Z^i_t$ is monotone increasing in $i$ for $i \geq 2$ (and all $t \in 1, \ldots, T$).  As a result, $\max_{t = 1, \ldots, T} Z^i_t$ is also monotone increasing in $i$ for $i \geq 2.$ Thus $\sum_{i = 2}^k \max_{t = 1, \ldots, T} Z^i_t \leq (k - 1)\times \max_{t = 1, \ldots, T} Z^k_t.$ Combining the above arguments, we conclude that $ Z_T - \max_{t = 1,\ldots,T} Z_t \leq  \sum_{i = 2}^k \max_{t = 1, \ldots, T} Z^i_t \leq (k - 1)\times \max_{t = 1, \ldots, T} Z^k_t.$ Or equivalently, $(k - 1)\times (- \max_{t = 1, \ldots, T} Z^k_t) \leq \max_{t = 1,\ldots,T} Z_t - Z_T.$ Since we assume $Z_t \in [0, 1]$ for $t = 1, \ldots, T,$ it must hold that $\max_{t = 1,\ldots,T} Z_t - Z_T \leq 1.$  Combining the above, we conclude that $- \max_{t = 1, \ldots, T} Z^k_t \leq \frac{1}{k - 1},$ completing the proof.  $\qed$
\endproof
\proof{Proof of Theorem\ \ref{thm:main2}: }
By Lemma\ \ref{lem:intuit1}, $\textsc{opt} = \sum_{i=1}^k E[\max_{t = 1,\ldots,T} Z^i_t] + \sup_{\tau \in {\mathcal T}} E[Z^{k+1}_{\tau}] = E_k + \sup_{\tau \in {\mathcal T}} E[Z^{k+1}_{\tau}] $.  In light of Lemma\ \ref{converge1}, $\tau_k \stackrel{\Delta}{=} \min\lbrace t = 1,\ldots,T: Z^k_t \geq - \frac{1}{k - 1} \rbrace$ is a well-defined stopping time satisfying $E[Z^k_{{\tau_k}}] \geq - \frac{1}{k - 1}$. Therefore $E_k - \textsc{opt} = - \sup_{\tau \in {\mathcal T}} E[Z^{k+1}_{\tau}] \leq - E[Z^{k+1}_{{\tau_{k+1}}}] \leq \frac{1}{k},$ completing the proof.  $\qed$
\endproof
\begin{remark}
We will later prove that $\tau_k$, as a valid stopping time for the original stopping problem, achieves an optimality gap of at most $\frac{1}{k}.$  Our polynomial-time algorithm for implementing a stopping policy is based upon this insight.\end{remark}
\ \\ For the setting in which Assumption \ref{assumptionA1} holds with $U = 1$, we now present a lower bound showing that Theorem\  \ref{thm:main2} is tight up to an absolute multiplicative constant factor of $\frac{1}{e}$ (asymptotically in $k$). First let us give the intuition behind this result. Our hard instance is a simple two-period problem (different for each $k$) in which the process is highly structured with the following properties: (1) the process is a martingale; (2) random variable $Z_1$ is a constant; (3) random variable $Z_2$ takes only two different values, one of which is zero. A straightforward induction implies that all of these properties will also hold for $\mathbf{Z}^j$ for all $j \geq 2$ (although they become non-positive processes). Furthermore, the martingale property and the discussion in Section\ \ref{sec:mainresultsexpansion} imply that the truncation error (after $j$ terms in our expansion) will simply equal $Z^{j+1}_1$ (itself the optimal stopping value associated with $\mathbf{Z}^{j+1}$ by the martingale property).  For the hard instance associated with $k$, we deliberately choose a specific probability assignment (in particular, the probability that $Z_2$ takes value zero is set to $\frac{1}{k + 1}$), such that it holds true that $Z^{j+1}_1 = (1 - \frac{1}{k+1}) Z^j_1$ for all $j \geq 2$. In other words, the truncation error after $j+1$ terms is a $(1 - \frac{1}{k + 1})$ fraction of the truncation error after $j$ terms for all $j \geq 2.$ Applying inductively, we conclude that the error if one truncates after $k$ terms, $-\sup_{\tau \in {\mathcal T}} E[Z^{k+1}_{\tau}],$ equals $-(1 - \frac{1}{k+1})^{k-1} \times Z^2_1$.  Since in our example it will hold that $Z^2_1 = -\frac{1}{k+1} \times (1 - \frac{1}{k+1})$, we conclude the following (deferring the proof to Appendix \ref{sec:ectheoryproof}).

\begin{lemma}[Lower bound]\label{lem:hardinstance}
For any given $k \geq 2$, there exists a 2-period OS problem with $P(Z_t \in [0,1]) = 1$ for all $t \in \lbrace 1,2 \rbrace$, yet $E_k - \textsc{opt} \geq \big( \frac{k}{k+1} (1 - \frac{1}{k+1})^k \big) \times \frac{1}{k}$.  Since $\lim_{k \rightarrow \infty} \big( \frac{k}{k+1} (1 - \frac{1}{k+1})^k \big) = e^{-1}$, we conclude that asymptotically Theorem\ \ref{thm:main2} is tight to within a multiplicative factor $e^{-1}$.
\end{lemma}

Namely, there is a precise sense in which the rate of convergence of the stated expansion cannot be substantially improved (in the worst case).  

\section{Polynomial-time Algorithms and their analysis}\label{sec:algoresults}
\indent We now describe our main algorithms that turn our expansion into computer executable instructions via Monte-Carlo simulation. Our first algorithm (Algorithm $\mathcal{A}$) outputs an approximate optimal stopping value. This value is achievable by the stopping policy computed by our second algorithm (Algorithm $\mathcal{A}'$).  Both algorithms enjoy explicit theoretical guarantees that relate the precision of the output with their runtime and sample complexity. Our algorithmic results, as stated in this section, should not be perceived as practical methods that could be implemented in actual applications.  Instead, they serve as a proof-of-concept that approximation algorithms for OS with runtime scaling polynomially in the time horizon and dimension exist.  The remainder of this section proceeds as follows.
\begin{itemize}
\item In Section\ \ref{subsec:algdescrip}, we provide very high-level descriptions of our algorithms.
\item In Section\ \ref{complexsec1}, we state our main algorithmic and complexity results and provide some additional intuition for our complexity bounds.
\item In Section\ \ref{sec: ecproofofmainalgo2}, we sketch the proofs of and provide in-depth intuition for our main algorithmic results.
\end{itemize}
\ \\For clarity of exposition (and out of space considerations), we defer the formal description of our algorithms and details of their analysis to Appendix\ \ref{sec:ecalgo}.

\subsection{High-level Description of Algorithms}\label{subsec:algdescrip}
Before diving into the description of the two algorithms (${\mathcal A}$ and ${\mathcal A}'$), we begin by formalizing our ``simulator" $\mathcal{B}$ of the underlying information process $\mathbf{Y}$, whose existence is guaranteed by Assumption\ \ref{assumptionB}, and which will act as a building block of all subsequent algorithms.  For $t \in \lbrace 1,\ldots,T \rbrace$ and $\gamma \in \mathcal{R}^{D \times t}$, let $\mathbf{Y}(\gamma)$ denote a random $D$ by $T$ matrix distributed as $\mathbf{Y}$, conditioned on the event $\lbrace \mathbf{Y}_{[t]} = \gamma \rbrace$.
 \\\hrule
\vspace{2pt}
 Simulation Oracle ${\mathcal B}$:
 \vspace{2pt}
 \hrule
 \vspace{4pt}
\indent\indent Input: $\ t\in \lbrace 0,\ldots,T \rbrace$ and $\gamma \in \mathcal{R}^{D \times t}$
\\\indent\indent \indent\indent Return an independent sample $\mathbf{Y}(\gamma)$ 
\vspace{2pt}
\hrule
\vspace{18pt}

Algorithm $\mathcal{A}$, which captures the main essence of our approach, is based on using an appropriate truncation of our expansion, with all terms estimated by simulation. Algorithm $\mathcal{A}'$ takes as input a trajectory $Y_{[T]},$ and outputs a stopping decision $t \in \lbrace 1,\ldots,T \rbrace$ associated to $Y_{[T]}.$ The decision is made by implementing a proper threshold policy on a certain higher order term of the expansion.  Below, we present a simplified rendition of the two algorithms. All references to $k$ and $\epsilon'$ in our informal description will ultimately be chosen as an appropriate function of $\epsilon$ and $\delta$.  Here we only introduce ``sketches" of our algorithms to convey the main ideas, deferring a formal description to Appendix\ \ref{sec:ecalgo}.
\ \\\hrule
\vspace{2pt}
 Algorithm ${\mathcal A}$ (\textsc{opt} approximation):
 \vspace{2pt}
 \hrule
 \vspace{4pt}
 \indent\indent Input: desired accuracy $\epsilon$ and confidence $\delta$
\\\indent\indent Select $k$
\\\indent\indent For $i = 1$ to $k$
\\\indent\indent\indent\indent Repeatedly call $\mathcal{B}$ to compute $\hat{L}_i$ (a Monte Carlo estimate of $L_i$)
\\\indent\indent Return $\sum_{i = 1}^k \hat{L}_i$ 
\vspace{2pt}
\hrule
\vspace{18pt}
\hrule
\vspace{2pt}
 Algorithm ${\mathcal A'}$ (stopping policy):
 \vspace{2pt}
 \hrule
 \vspace{4pt}
\indent\indent Input: a trajectory $Y_{[T]}$ and desired accuracy $\epsilon$
\\\indent\indent Select $k, \epsilon'$
\\\indent\indent For $t = 1$ to $T$
\\\indent\indent\indent\indent Repeatedly call $\mathcal{B}$ to compute $\hat{Z}^k_t$ (a Monte Carlo estimate of $Z^k_t$)
\\\indent\indent\indent\indent Return \textsc{stop}  if $(- \hat{Z}^k_t < \epsilon'\ \ \textrm{or}\ \   t = T)$
\vspace{2pt}
\hrule
\vspace{18pt}

Let us elaborate a bit on how the auxiliary terms $\hat{L}_i$ and $\hat{Z}^k_t$ are computed. Recall that $L_i = E\bigg[\max_{j = 1,\dots, T} Z^i_j \bigg],$ where $Z^i_t$ is recursively defined by $Z^1_t = Z_t$ and $ Z^i_t = Z^{i - 1}_t -  E\big[\max_{j = 1,\dots, T}Z^{i-1}_j \big|\mathcal{F}_t\big]$ for $i \geq 2.$ The expansion terms $\hat{L}_i$ are computed using Monte Carlo simulation as follows.  First, generate a sufficient number of independent trajectories of the information process.  Second, for each trajectory, compute $\hat{Z}^i_t$ (the approximation to $Z^i_t$) to a desired precision for $t = 1,\ldots,T$.  Third, for each trajectory, compute $\max_t \hat{Z}^i_t$.  Finally, average these estimates over all trajectories and return the output sample-average approximation (SAA).  Let us point out that since $Z^i_t$ is defined recursively in terms of $Z^{i-1}_t$ and $E\big[\max_{j = 1,\dots, T}Z^{i-1}_j \big|\mathcal{F}_t\big]$, the computation of $\hat{Z}^i_t$ will itself require a recursive/nested simulation of many approximate estimates of terms $\hat{Z}^{i-1}_{t'}$ for different $t'$ (drawn from their own trajectories of the information process), themselves requiring the simulation of many estimates of terms $\hat{Z}^{i-2}_{t''}$, etc.  Thus to compute $\hat{Z}^i_t$, one needs a nested Monte Carlo simulation with $i$ levels of nesting.  In general nested simulation is very costly, but the key point is that we can take $i$ to be $O(\frac{1}{\epsilon})$, thus yielding a constant level of nesting (not depending on the length of the horizon or the dimension).

\subsection{Algorithmic analysis and complexity guarantees}\label{complexsec1}
In this section we state various complexity guarantees of our algorithms (under our different sets of assumptions).  We begin by formalizing our computational model for analyzing algorithms.  Although the original theory for analyzing the complexity of algorithms was based on the Turing-machine/bit complexity model, it was later observed that different complexity models were more appropriate for the study of certain stochastic optimization problems in which one's algorithm is given oracle access to a simulator (\cite{de2008complexity}).  As high-dimensional OS falls into this category, we adopt a model in line with other such works (\cite{shapiro2005complexity,swamy2012sampling}), in which one is given access to an appropriate (simulation) oracle and can perform basic arithmetic operations on real numbers (regardless of the number of bits).  Such a model falls within the framework of information complexity (\cite{traub1998complexity}), and we refer the reader to \cite{de2008complexity} for a more detailed comparison with the Turing-machine/bit complexity model.
\\\indent Consistent with the assumptions in these works, we suppose that addition, subtraction, multiplication, division, maximum, and minimum of any two numbers can be done in one unit of time, as can the raising of one number to the power of another, irregardless of the values of those numbers.  We also suppose that the median of $n$ numbers can be computed in $n$ units of computational time, using e.g. the celebrated linear-time method for computing the median (\cite{cormen2009introduction}), and ignoring the absolute constant governing the associated linear runtime.  We ignore all computational costs associated with reading, writing, and storing numbers in memory, as well as inputting numbers to functions. In general our model allows for computation over real numbers.  We also assume that drawing a sample from standard Bernoulli, uniform ($\mathcal{U}[0, 1]$) or normal ($\mathcal{N}(0,1)$) r.v.s takes one unit of time.  Most importantly, we assume access to an oracle which yields (conditioned) trajectories of the information process, and evaluates the relevant reward functions, as per Assumption\ \ref{assumptionB}.

We now state our main algorithmic results. Our results are twofold: $(1)$ computing an $\epsilon$-approximation for the optimal value \textsc{opt}; and $(2)$ providing an $\epsilon-$approximate OS strategy.  We prove that under Assumption \ref{assumptionB}, and either of Assumptions \ref{assumptionA1} or \ref{assumptionA2}, for any given error parameter $\epsilon$, our algorithms can achieve these two goals (in expectation or with high probability).  Furthermore, our algorithm's runtime is only polynomial in $T$ (for a given $\epsilon$), and depends on the dimension (and state-space more generally) only through the cost of simulating individual sample paths, where only a polynomial number of such simulations are needed.  We defer the proofs to Appendix\ \ref{sec:ecalgo}.

\begin{theorem}[Optimal value approximation]\label{thm:mainalgo2}
Suppose that Assumptions \ref{assumptionA1} and \ref{assumptionB} hold, with $U = 1$ in Assumption\ \ref{assumptionA1} (by normalizing this is w.l.o.g.). Then Algorithm ${\mathcal A},$ with a proper choice of parameters, takes as input any $\epsilon,\delta \in (0,1)$, and achieves the following.  In total computational time at most $(C + G + 1) \times ( 1  +  \log(\delta^{-1})) \times \exp(100 \epsilon^{-2})\times  T^{8 \epsilon^{-1}}$, returns a random number $X$ satisfying $P\big( |X - \textsc{opt}| \leq \epsilon \big) \geq 1 - \delta$.
\end{theorem}

\begin{theorem}[Optimal stopping policy]\label{thm:mainpolicy1}
Under the same assumptions as Theorem \ref{thm:mainalgo2}, Algorithm $\mathcal{A}'$ implements a randomized stopping time $\tau_{\epsilon}$ s.t. $  E[Z_{\tau_{\epsilon}}] \geq \textsc{opt} - \epsilon$.  In addition, at each time step, the decision of whether to stop (if one has not yet stopped) requires total computational time at most $(C + G + 1) \times \exp(250 \epsilon^{-2}) \times T^{12\epsilon^{-1}}.$
\end{theorem}

\begin{theorem}[Optimal value approximation for unbounded rewards]\label{thm:mainmax}Suppose that 
\\Assumptions \ref{assumptionA2} and \ref{assumptionB} hold.  Then Algorithm $\mathcal{A}$, with properly chosen parameters and an additional truncation step, takes any $\epsilon,\delta \in (0,1)$ and achieves the following.  In total computational time at most 
$(C+ G+ 4) \times ( 1+  \log(\delta^{-1}) ) \times \exp\bigg(6.1 \times 10^{8} \times \gamma_0^8 \times \epsilon^{-6}\bigg) \times T^{ 2 \times 10^4 \times \gamma_0^{4} \times \epsilon^{-3}},$ returns a random number $X$ satisfying $P\bigg(\frac{\big| X - \textsc{opt}\big|}{\textsc{opt}} \leq \epsilon \bigg) \geq 1 - \delta.$
\end{theorem}

\begin{theorem}[Optimal stopping policy for unbounded rewards]\label{thm:mainmaxpolicy} Under the \\same assumptions as Theorem \ref{thm:mainmax}, Algorithm $\mathcal{A}'$ (with properly chosen parameters and an additional truncation step) implements a randomized stopping time $\tau_{\epsilon}$ s.t. $\frac{E[Z_{\tau_{\epsilon}}]}{\textsc{opt}} \geq 1 - \epsilon$.  In addition, at each time step, the decision of whether to stop (if one has not yet stopped) requires total computational time at most $(C+ G+ 4) \times ( 1+  \log(\delta^{-1}) ) \times \exp\bigg(1.6 \times 10^{9} \times \gamma_0^8 \times \epsilon^{-6}\bigg) \times T^{ 3 \times 10^4 \times \gamma_0^{4} \times \epsilon^{-3}}$.
\end{theorem}

\subsubsection{Additional intuition for complexity bounds.}\label{moreintuitmain} Here we provide some additional intuition for our complexity bounds, and provide a more detailed discussion of the relevant intuition in Appendix\ \ref{sec:ecmoreintuit}.  

For any particular $i$, $t$, and trajectory, the computation of $\hat{Z}^i_t$ (on that trajectory) involves first generating sufficiently many (say $N$) trajectories and then computing $\hat{Z}^{i-1}_1, \dots, \hat{Z}^{i - 1}_T$ along each trajectory (to estimate $E[\max_{j = 1,\ldots,T} Z^{i-1}_j | {\mathcal F}_t]$ using SAA).  Roughly speaking, one needs to compute $N \times T$ terms of the form $\hat{Z}^{i - 1}_{t'}$ in order to compute $\hat{Z}^i_t$, where $N$ depends on the desired accuracy.  Applying this logic inductively for $i, i-1, \dots, 1$, we conclude that computing terms in the $i$-th level of the expansion roughly requires computational effort $(N T)^i.$  By Theorem \ref{thm:main2}, to attain an $\epsilon$ approximation of \textsc{opt}, one only needs to compute out to level $i = O(\frac{1}{\epsilon})$ in the expansion.  This intuitively results in a complexity of $(N T)^{O(\frac{1}{\epsilon})}$, and this is precisely where the $T^{O(\frac{1}{\epsilon})}$ term comes from in our complexity bounds.  The term involving $C + G$ comes from the cost of each individual simulation and reward function evaluation, and the $\log(\delta^{-1})$ term comes from standard results in probability and the analysis of algorithms for deriving high-confidence estimators.  

We now comment on $N$ (which has so far been unspecified) and the term $\exp\big( O(\epsilon^{-2}) \big)$ in our complexity bounds.  To estimate $Z^i_t$, we must estimate both $Z^{i-1}_t$ and $E[\max_{j = 1,\ldots,T} Z^{i-1}_j | {\mathcal F}_t]$, and then take the difference of these estimates.  Thus to yield an estimate with $\epsilon$ error for $Z^i_t$, a union bound indicates we must estimate each of $Z^{i-1}_t$ and  $E[\max_{j = 1,\ldots,T} Z^{i-1}_j | {\mathcal F}_t]$ to error $\frac{\epsilon}{2}$.  Applying this logic inductively, to estimate terms at recursive depth $O(\frac{1}{\epsilon})$ to error $\epsilon$ our algorithm must estimate certain terms (deep in the recursion) to error $\frac{1}{2^{\frac{1}{\epsilon}}}$.  The exact expression $\exp(100 \epsilon^{-2}),$ which is even larger than $2^{\frac{1}{\epsilon}}$, arises from roughly the same intuition (albeit from a more sophisticated recursive analysis accounting precisely for how errors accumulate through the different levels of nesting, see the discussion in Appendix\ \ref{sec:ecmoreintuit}).  Furthermore, we note that this logic (of needing greater accuracy at different depths) implies that in our actual algorithm $N$ is selected to be a different function of $\epsilon$ (and other parameters) at different recursive depths.  

The logic behind the complexity bound appearing in Theorem\ \ref{thm:mainmax} (namely, the case with possibly unbounded rewards) is very similar, the main difference being that the algorithm must perform an additional truncation step on all rewards.  The truncation threshold is chosen (as a function of $\frac{1}{\epsilon}$) to trade off truncation error with computational complexity.  As a result, the algorithm must go out to greater depth in the expansion (albeit still only depending on $\epsilon$ / independent of $T$ and $D$) to yield the same error $\epsilon$, which by our intuition leads to a larger exponent of $T$ compared to Theorem 3 and 4.  

\subsection{Proof sketches and in-depth intuition for algorithmic results}\label{sec: ecproofofmainalgo2}
In this section we sketch the proofs of and provide in-depth intuition for our main algorithmic results.  We shall impose Assumption \ref{assumptionA1} with $U = 1$ unless otherwise specified. We defer the formal description of our algorithms and details of their analysis to Appendix\ \ref{sec:ecalgo}.
 \subsubsection{Polynomial time algorithms for computing $Z^k_t$}
Our main algorithmic theorems follow from the expansion established in earlier sections. First, we need to obtain polynomial time algorithms for computing $Z^k_t.$  To do so, we devise a randomized algorithm ${\mathcal B}^k$ that takes as inputs $t \in \lbrace 1,\ldots,T \rbrace, \gamma \in \mathcal{R}^{D \times t}, \epsilon \in (0,1), \delta \in (0,1)$, and returns a value which with probability at least $1 - \delta$ is additively within $\epsilon$ of $Z^k_t(\gamma)$.  Here $Z^k_t(\gamma)$ denotes the r.v. $Z^k_t$ conditioned on the event $\lbrace Y_{[t]} = \gamma \rbrace$ (where we note that $Z^k_t$ is a deterministic function of $Y_{[t]}$ by measurability).  We use a Monte Carlo approach to algorithmically implement our expansion, with algorithm ${\mathcal B}^{k+1}$ recursively calling algorithm ${\mathcal  B}^k$ in the natural manner suggested by definitions and our expansion.  

We now sketch in greater detail the intuition and inductive logic behind ${\mathcal B}^k$.
For the base case $k = 1,$ ${\mathcal B}^1$ simply makes an evaluation call to $\mathcal{B}$ with any input $\gamma \in \mathcal{R}^{D \times t},$ outputting the exact value of $Z^1_t(\gamma)$ in a runtime of $G$ units of time.  For $k \geq 1$, and supposing we have (inductively) defined algorithm ${\mathcal  B}^k$ which can approximately simulate $Z^k_t$ as desired, we construct algorithm ${\mathcal B}^{k+1}$ based on the following logic.  By definition, $Z^{k +1}_t(\gamma) = Z^k_t(\gamma) - E[\max_{1 \leq i \leq T} Z^k_i |Y_{[t]} = \gamma].$  We thus compute $Z^{k+1}_t(\gamma)$ by a recursive call of ${\mathcal B}^k$ to compute $Z^k_t(\gamma)$, and an appropriate number of additional recursive calls of ${\mathcal B}^k$ to estimate the conditional expectation $E[\max_{1 \leq i \leq T} Z^k_i |Y_{[t]} = \gamma]$ using Monte Carlo simulation.  We note that two fundamental types of error must be controlled: (1) error introduced from the fact that each call to ${\mathcal B}^k$ only gurantees an approximation to $Z^k_t$ (depending on the parameters it is called with); and (2) error introduced from using Monte Carlo to estimate the conditional expectation. A careful analysis yields the following result.
 \begin{lemma}\label{almostthere1}
For all $k \geq 1$, $t \in \lbrace 1,\ldots,T \rbrace$, $\gamma \in \mathcal{R}^{D \times t}$, $\epsilon,\delta \in (0,1)$, there exists algorithm ${\mathcal B}^k$ that achieves the following when evaluated on $t,\gamma,\epsilon,\delta.$  In total computational time at most 
$$(C + G + 1)\times \log(2\delta^{-1})\times 10^{2 k^2} \times \epsilon^{- 2 k} \times (T+2)^k \times \big(1 + \log(\frac{1}{\epsilon}) + \log(T) \big)^k,$$
${\mathcal B}^k$ returns a random number $X$ satisfying $P\big( |X - Z^k_t(\gamma)| > \epsilon \big) \leq \delta.$  
\end{lemma}
We defer the formal description of algorithm ${\mathcal B}^k$ to Appendix\ \ref{ecBk}, and the formal analysis and proof of Lemma\ \ref{almostthere1} to Appendix\ \ref{eclemma3}.

\subsubsection{Optimal value approximation and proof sketch of Theorem \ref{thm:mainalgo2}}
Our main result for approximating the optimal value, Theorem \ref{thm:mainalgo2}, follows directly from Lemma\ \ref{almostthere1} and Theorem \ref{thm:main2}.  Indeed, by Theorem \ref{thm:main2}, truncating the expansion after the $n$-th term yields $E_n = \sum_{k = 1}^n L_k$, a $\frac{1}{n}-$approximation of \textsc{opt}.  Thus to compute an $\epsilon$-approximation of \textsc{opt} with high probability, it suffices to compute an $\frac{\epsilon}{2}-$approximation of $E_{\lceil\frac{2}{\epsilon}\rceil}$ with high probability.  By a simple union bound, this may be achieved by computing each $L_k$ (for $1\leq k \leq \lceil\frac{2}{\epsilon}\rceil$) to within additive error $\frac{\epsilon}{2}\big(\lceil\frac{2}{\epsilon}\rceil\big)^{-1}$ with sufficiently high probability.  Since by definition $L_k = E[\max_{1 \leq i \leq T} Z^k_i]$, this can be achieved by combining Lemma\ \ref{almostthere1} (to simulate the $Z^k_t$ terms) with a standard Monte Carlo approach to estimating the relevant expectation.  Combining with an argument similar to that used in the proof of Lemma\ \ref{almostthere1} can again control all relevant errors and achieve the desired runtime guarantees and sample complexity.   We defer the formal description of the associated algorithm ${\mathcal A}$, its analysis, and the proof of Theorem\ \ref{thm:mainalgo2}, to Appendix\ \ref{ecthm3}.

\subsubsection{Approximately optimal stopping times and proof sketch of Theorem \ref{thm:mainpolicy1}}\label{goodpolicy}
We now discuss how to use our simulation-based approach to implement approximately optimal stopping strategies, sketching the proof of Theorem\ \ref{thm:mainpolicy1}.  We begin with the following lemma relating the value achieved by a single stopping policy across different stopping problems (defined by $\mathbf{Z}^k$).  The logic is essentially the same as that of Lemma\ \ref{lem:intuit1}, with the result following from our definitions, the optional stopping theorem, and a straightforward induction, and we omit the details.  Out of considerations of space and clarity of exposition, we do not formalize the notion of randomized stopping time here.  We instead refer the reader to e.g. \cite{solan2012random}, and note that all relevant concepts can be formalized in our computational model by considering an appropriate augmented probability space, and we omit the details. 

\begin{lemma}\label{good1}
For all (possibly randomized) integer-valued stopping times $\tau$ adapted to ${\mathcal F}$ which w.p.1 belonging to $\lbrace 1,\ldots,T \rbrace$, and all $k \geq 1$, $E[Z_{\tau}] = E[Z^k_{\tau}] + \sum_{i=1}^{k-1} E[\max_{t = 1,\ldots,T} Z^i_t]$.
\end{lemma}
Combining Lemma\ \ref{good1} with Lemma\ \ref{lem:intuit1}, Lemma\ \ref{converge1}, and Theorem\ \ref{thm:main1}, we are led to the following corollary.

\begin{corollary}\label{good2}
For $k \geq 1$, let $\tau_k$ denote the stopping time that stops the first time that $-Z^{k+1}_t \leq \frac{1}{k}$, where by Lemma\ \ref{converge1} such a time exists w.p.1.  Then $\textsc{opt} - E[Z_{\tau_k}]  \leq \frac{1}{k}$.
\end{corollary}

A suitable approximation to the threshold policy $\tau_k$ can be found  by approximately computing $Z^{k+1}_t$ at each time step.  By Lemma \ref{almostthere1}, such a computation can be carried out by calling $\mathcal{B}^k$ an appropriate number of times.  The desired efficiently computable randomized stopping time $\tau_{\epsilon}$ is simply a natural approximation to $\tau_k$ for $k = O(\frac{1}{\epsilon})$, and the formal proof of Theorem\ \ref{thm:mainpolicy1} proceeds by making this completely precise and carefully bounding all sources of errors.

\subsubsection{Unbounded setting and proof sketch of Theorems\ \ref{thm:mainmax} and \ref{thm:mainmaxpolicy}}\label{maxsubsec}

Our main results for approximating the optimal value and implementing an approximately optimal stopping policy under Assumption \ref{assumptionA2}, Theorems \ref{thm:mainmax} and \ref{thm:mainmaxpolicy}, follow by combining our algorithms $({\mathcal A}$ and ${\mathcal A}'$) and results (Theorems\ \ref{thm:mainalgo2} and \ref{thm:mainpolicy1}) for the bounded case with an extra truncation step.  In particular, we select (as a function of the desired accuracy $\epsilon$) a truncation level $U^{\epsilon}$, and truncate all rewards at this level, allowing us to apply our results for the bounded case.  In addition, we prove results which allow us to control the associated truncation error.  We defer the formal proofs to Appendix \ref{sec:ecmaximization}.

\section{Bermudan Max Call}\label{sec:maxcall}
In this section we prove that the Bermudan max call, a popular test benchmark for OS algorithms (\cite{andersen2004primal}), satisfies Assumptions \ref{assumptionA2} and \ref{assumptionB} (for appropriate parameters).  We note that for the Bermudan max call, the unit simulation cost $C$ and $G$ can be explicitly bounded as simple functions of $T$ and $D$, and our simulation-based algorithm is thus able to return (with high probability) a value within $\epsilon$ of the true option price (in a relative error metric) with computational complexity scaling polynomially in $T$ and $D$.

\subsection{Basic setup}

Following \cite{andersen1999simple}, we assume that the underlying stochastic asset process evolves according to the dynamics of a multi-dimensional geometric B.m., i.e.
$\ d G^k_t = (r -\varrho) G^k_t d t + \sigma G^k_t d W^k_t,$ for all $t \geq 0 $ and $k = 1,...,D,$\ with $G^k_0 = y_0$ for all $k$,
where $\lbrace W^k_t, t \geq 0 ,  k = 1,...,D \rbrace$ are independent one-dimensional B.m. and $r,\varrho,\sigma$ are constants capturing the interest rate, the dividend rate and the volatility respectively.  We assume that the option lives in $[0, \mathcal{J}]$ for some horizon length ${\mathcal J}$,and that the exercising dates are restricted to $\mathcal{S} = \{t_0,..., t_T\}$ where $t_i \ \stackrel{\Delta}{ =}\ \frac{i}{T} \mathcal{J}$ for $i \in \lbrace 0,\ldots,T \rbrace$ (assumed to start from $0$ in line with standard notation in this setting). We also assume that $\mathcal{J}, T$ are both positive integers, and $T$ is divisible by $\mathcal{J}.$ One should think of $\mathcal{J}$ as, for example, the number of months, and $\frac{T}{\mathcal{J}}$ the number of exercising opportunities within one month. The payoff is given by
$g_t(G_t) =  \big(\max(G^1_t,...,G^D_t) -  \kappa\big)^+, $
with $\kappa > 0$ some positive constant and $x^+ = \max(x,0)$ for any scalar $x$. The OS problem associated to Bermudan max call pricing is thus 
$$
\textsc{opt} \ =\ \sup_{\tau \in \lbrace 0,\ldots,T \rbrace} E\bigg[e^{- r t_\tau} g_{t_\tau}(G_{t_\tau})\bigg], 
$$
with $\tau$ a stopping time adapted to the natural filtration generated by the asset process.  Since here the information process is simply a $D$-dimensional geometric B.m. discretized on the mesh $\lbrace t_i, i = 0,\ldots,T \rbrace$, the existence of an efficient simulator (as required by Assumption\ \ref{assumptionB}) is easily verified.

\subsection{Main result}
Our main result for Bermudan max call is summarized in the next corollary.
\begin{corollary}\label{cormaxcall}
Suppose $\varrho > 1.5 \sigma^2.$  Then there exists an algorithm (namely Algorithm $\mathcal{A}$, with properly chosen parameters and an additional truncation step) that takes as input $\epsilon,\zeta \in (0,1)$ and achieves the following for a Bermudan max call instance (under our computational model).  In total computational time at most $(10 T D + 4) \times ( 1 +  \log(\zeta^{-1}) ) \times \exp\bigg(6.1 \times 10^{8} \times \gamma_0^8 \times \epsilon^{-6}\bigg) \times T^{ 2 \times 10^4 \times \gamma_0^{4} \times \epsilon^{-3}},$ returns a random number $X$ satisfying $P\bigg(\frac{\big| X - \textsc{opt}\big|}{\textsc{opt} + \kappa} \leq \epsilon \bigg) \geq 1 - \zeta$, where $\gamma_0 =  10^4 \exp\big(2 (\frac{\varrho + \frac{\sigma^2}{2}}{\sigma})^2\big) (1 + \frac{20 \sigma^2}{\varrho - 1.5 \sigma^2}).$
\end{corollary}

We defer the proof to Appendix\ \ref{sec:ecmaxcall}.  To achieve the above result we only require that $\varrho > 1.5 \sigma^2,$ and the mild assumptions (for technical convenience) that $T, {\mathcal J}$ are integers with $T$ divisible by $\mathcal{J}.$  We impose no additional constraints on parameters $D,\mathcal{J}, T, r, \kappa $ and $y_0.$ In particular, $D, T$ and $\mathcal{J}$ are allowed to take arbitrarily large values, covering all parameter regimes.  We note that so long as $\varrho,\sigma,\frac{1}{\varrho - 1.5 \sigma^2},$ and $\frac{\kappa}{\textsc{opt}}$ are bounded independent of $T$ and $D$, the above will yield a polynomial time algorithm with relative error $\epsilon$ (for any fixed $\epsilon$, in our computational model) in the setting of the Bermudan max call.  Regarding the assumption $\varrho > 1.5 \sigma^2,$ let us point out that previous studies using Bermudan max call as their numerical benchmark (e.g. \cite{andersen2004primal}, \cite{becker2019solving}, \cite{rogers2015bermudan}, \cite{belomestny2011pricing}) all have parameter choices satisfying the condition.  Furthermore, such an assumption is essentially inevitable if we wish to deal with arbitrarily large $\mathcal{J}$ (e.g. perpetual options).  Indeed, when $\mathcal{J}$ is large and the condition $\varrho > 1.5 \sigma^2$ is violated, the path-wise maximum discounted payoff will be well-approximated by $\sup_{t \geq 0} e^{-r t} g_t(Y_t),$ which will have infinite variance, leading to severe challenges for any Monte-Carlo based method.

\section{Conclusion}\label{conc}
In this work, we proved the existence of a polynomial time (in the horizon $T$ and dimension $D$) algorithm for complex OS problems with fixed accuracy and access to an efficient simulator.  The algorithms come with provable performance guarantees under either boundedness or mild concentration assumptions.  Although their complexity is polynomial, it is impractical, as was the case for the first polynomial-time algorithms for other problems in the literature.  Our work leaves many interesting directions for future research. 
\\\indent \textbf{The design of practical algorithms with polynomial runtime guarantees.}  A natural question is whether there exists a practical algorithm with polynomial runtime guarantees.  One possible approach is to develop different expansions for the OS value which converge more rapidly, and/or to identify natural assumptions under which such expansions converge more rapidly.  Such approaches could also be combined with techniques from reinforcement learning, ADP, simulation, and parallel computing.
\\\indent\textbf{Generalization to high-dimensional online decision-making broadly and other applications.}  We believe that our methodology can be extended to a broad family of high-dimensional online decision-making problems in economics, mechanism design / prophet inequalities, finance, statistics, machine learning, operations, and computer science.  Preliminary results along these lines for certain problems in multiple stopping, dynamic pricing, and online combinatorial optimization appear in \cite{chen2021efficient}.  
\\\indent\textbf{Lower bounds and computational complexity.}  An interesting set of questions revolve around proving lower bounds on the computational and sample complexity for the problems studied, e.g. path-dependent OS.  There has been much interesting recent work laying out a theory of computational complexity (with positive and negative results) in the settings of stochastic control and reinforcement learning (\cite{shapiro2005complexity,halman2014fully, halman2015computationally, chen2017lower, sidford2018near, du2019good, wang2021exponential}).  Better understanding the connection between our approach and those works, as well as the many existing models for online optimization, remains an interesting direction for future research.  Furthermore, the question of what quality of approximation can be (efficiently) achieved with a given depth of nesting remains a very interesting open question.  Such questions relate to recent studies on the power of adaptivity in optimization (e.g. \cite{balkanski2017limitations}), as well as the literature on parallel and quantum computing, (multistage stochastic) convex optimization, and prophet inequalities.

\section*{Acknowledgements} 
The authors gratefully acknowledge Sid Banerjee, Erhan Bayraktar, Denis Belomestny, Thomas Bruss, Agostino Capponi, Jim Dai, Damek Davis, Mark Davis, Vivek Farias, Paul Glasserman, Nir Halman, Shane Henderson, Saul Jacka, Bobby Kleinberg, Ozalp Ozer, Philip Protter, Chris Rogers, John Schoenmakers, Timur Tankayev, John Tsitsiklis, Alberto Vera, and David Williamson for several helpful conversations and insights. Special thanks go to Kaidong Zhang and Di Wu, for their tremendous help with numerical implementations, and to Bobby Kleinberg, for pointing out that the median trick could be applied in our setting. The authors thank Jim Dai for his organization of a seminar series on Reinforcement Learning for Processing Networks, as well as all participants of that seminar series.  The authors also thank the organizers and participants of the 2018 Symposium on optimal stopping in honor of Larry Shepp held at Rice University.

\medskip
\bibliographystyle{informs2014} 
\bibliography{Opt_stop_Orig_Bib_8_6_2023} 

\begin{thebibliography}{104}
\providecommand{\natexlab}[1]{#1}
\providecommand{\url}[1]{\texttt{#1}}
\providecommand{\urlprefix}{URL }

\bibitem[{Abraham et~al.(2010)Abraham, Fiat, Goldberg, \protect\BIBand{}
  Werneck}]{abraham2010highway}
Abraham I, Fiat A, Goldberg AV, Werneck RF (2010) Highway dimension, shortest
  paths, and provably efficient algorithms. \emph{Proceedings of the
  twenty-first annual ACM-SIAM symposium on Discrete Algorithms}, 782--793
  (SIAM).

\bibitem[{Aingworth et~al.(2000)Aingworth, Motwani, \protect\BIBand{}
  Oldham}]{aingworth2000accurate}
Aingworth D, Motwani R, Oldham JD (2000) Accurate approximations for asian
  options. \emph{Proceedings of the eleventh annual ACM-SIAM symposium on
  Discrete algorithms}, 891--900.

\bibitem[{Alaei et~al.(2022)Alaei, Makhdoumi, Malekian, \protect\BIBand{}
  Peke{\v{c}}}]{alaei2022revenue}
Alaei S, Makhdoumi A, Malekian A, Peke{\v{c}} S (2022) Revenue-sharing
  allocation strategies for two-sided media platforms: Pro-rata vs.
  user-centric. \emph{Management Science} .

\bibitem[{Andersen \protect\BIBand{} Broadie(2004)}]{andersen2004primal}
Andersen L, Broadie M (2004) Primal-dual simulation algorithm for pricing
  multidimensional american options. \emph{Management Science}
  50(9):1222--1234.

\bibitem[{Andersen(1999)}]{andersen1999simple}
Andersen LB (1999) A simple approach to the pricing of bermudan swaptions in
  the multi-factor libor market model. \emph{Available at SSRN 155208} .

\bibitem[{Aouad \protect\BIBand{} Segev(2021)}]{aouad2021display}
Aouad A, Segev D (2021) Display optimization for vertically differentiated
  locations under multinomial logit preferences. \emph{Management Science}
  67(6):3519--3550.

\bibitem[{Arora(1996)}]{arora1996polynomial}
Arora S (1996) Polynomial time approximation schemes for euclidean tsp and
  other geometric problems. \emph{Proceedings of 37th Conference on Foundations
  of Computer Science}, 2--11 (IEEE).

\bibitem[{Arora(1997)}]{arora1997nearly}
Arora S (1997) Nearly linear time approximation schemes for euclidean tsp and
  other geometric problems. \emph{Proceedings 38th Annual Symposium on
  Foundations of Computer Science}, 554--563 (IEEE).

\bibitem[{Arora \protect\BIBand{} Barak(2009)}]{arora2009computational}
Arora S, Barak B (2009) \emph{Computational complexity: a modern approach}
  (Cambridge University Press).

\bibitem[{Baas et~al.(2023)Baas, Boucherie, \protect\BIBand{}
  Braaksma}]{baas2023sampling}
Baas S, Boucherie RJ, Braaksma A (2023) A sampling-based method for gittins
  index approximation. \emph{arXiv preprint arXiv:2307.11713} .

\bibitem[{Balkanski et~al.(2017)Balkanski, Rubinstein, \protect\BIBand{}
  Singer}]{balkanski2017limitations}
Balkanski E, Rubinstein A, Singer Y (2017) The limitations of optimization from
  samples. \emph{Proceedings of the 49th Annual ACM SIGACT Symposium on Theory
  of Computing}, 1016--1027.

\bibitem[{Baricz(2008)}]{baricz2008mills}
Baricz {\'A} (2008) Mills' ratio: Monotonicity patterns and functional
  inequalities. \emph{Journal of Mathematical Analysis and Applications}
  340(2):1362--1370.

\bibitem[{Bayer et~al.(2016)Bayer, Friz, \protect\BIBand{}
  Gatheral}]{bayer2016pricing}
Bayer C, Friz P, Gatheral J (2016) Pricing under rough volatility.
  \emph{Quantitative Finance} 16(6):887--904.

\bibitem[{Bayer et~al.(2021)Bayer, Redmann, \protect\BIBand{}
  Schoenmakers}]{bayer2018dynamic}
Bayer C, Redmann M, Schoenmakers J (2021) Dynamic programming for optimal
  stopping via pseudo-regression. \emph{Quantitative Finance} 21(1):29--44.

\bibitem[{Becker et~al.(2019)Becker, Cheridito, \protect\BIBand{}
  Jentzen}]{becker2019deep}
Becker S, Cheridito P, Jentzen A (2019) Deep optimal stopping. \emph{Journal of
  Machine Learning Research} 20:74.

\bibitem[{Becker et~al.(2021)Becker, Cheridito, Jentzen, \protect\BIBand{}
  Welti}]{becker2019solving}
Becker S, Cheridito P, Jentzen A, Welti T (2021) Solving high-dimensional
  optimal stopping problems using deep learning. \emph{European Journal of
  Applied Mathematics} 32(3):470--514.

\bibitem[{Beier \protect\BIBand{} V{\"o}cking(2004)}]{beier2004random}
Beier R, V{\"o}cking B (2004) Random knapsack in expected polynomial time.
  \emph{Journal of Computer and System Sciences} 69(3):306--329.

\bibitem[{Beinker \protect\BIBand{} Schlenkrich(2017)}]{beinker2017accurate}
Beinker M, Schlenkrich S (2017) Accurate vega calculation for bermudan
  swaptions. \emph{Novel Methods in Computational Finance} 65--82.

\bibitem[{Belomestny(2011)}]{belomestny2011pricing}
Belomestny D (2011) Pricing bermudan options by nonparametric regression:
  optimal rates of convergence for lower estimates. \emph{Finance and
  Stochastics} 15(4):655--683.

\bibitem[{Belomestny et~al.(2019)Belomestny, Hildebrand, \protect\BIBand{}
  Schoenmakers}]{belomestny2019optimal}
Belomestny D, Hildebrand R, Schoenmakers J (2019) Optimal stopping via pathwise
  dual empirical maximisation. \emph{Applied Mathematics \& Optimization}
  79(3):715--741.

\bibitem[{Belomestny et~al.(2020)Belomestny, Kaledin, \protect\BIBand{}
  Schoenmakers}]{belomestny2019semi}
Belomestny D, Kaledin M, Schoenmakers J (2020) Semitractability of optimal
  stopping problems via a weighted stochastic mesh algorithm.
  \emph{Mathematical Finance} 30(4):1591--1616.

\bibitem[{Belomestny et~al.(2015)Belomestny, Ladkau, \protect\BIBand{}
  Schoenmakers}]{belomestny2015multilevel}
Belomestny D, Ladkau M, Schoenmakers J (2015) Multilevel simulation based
  policy iteration for optimal stopping--convergence and complexity.
  \emph{SIAM/ASA Journal on Uncertainty Quantification} 3(1):460--483.

\bibitem[{Belomestny \protect\BIBand{} Milstein(2006)}]{belomestny2006monte}
Belomestny D, Milstein GN (2006) Monte carlo evaluation of american options
  using consumption processes. \emph{International Journal of theoretical and
  applied finance} 9(04):455--481.

\bibitem[{Belomestny et~al.(2013)Belomestny, Schoenmakers, \protect\BIBand{}
  Dickmann}]{belomestny2013multilevel}
Belomestny D, Schoenmakers J, Dickmann F (2013) Multilevel dual approach for
  pricing american style derivatives. \emph{Finance and Stochastics}
  17(4):717--742.

\bibitem[{Belomestny et~al.(2018)Belomestny, Schoenmakers, Spokoiny,
  \protect\BIBand{} Tavyrikov}]{belomestny2018optimal}
Belomestny D, Schoenmakers J, Spokoiny V, Tavyrikov Y (2018) Optimal stopping
  via deeply boosted backward regression. \emph{arXiv preprint
  arXiv:1808.02341} .

\bibitem[{Belomestny et~al.(2011)}]{belomestny2011rates}
Belomestny D, et~al. (2011) On the rates of convergence of simulation-based
  optimization algorithms for optimal stopping problems. \emph{The Annals of
  Applied Probability} 21(1):215--239.

\bibitem[{Bezerra et~al.(2018)Bezerra, Ohashi, Russo, \protect\BIBand{}
  de~Souza}]{bezerra2018discrete}
Bezerra SC, Ohashi A, Russo F, de~Souza F (2018) Discrete-type approximations
  for non-markovian optimal stopping problems: Part ii. \emph{Methodology and
  Computing in Applied Probability} 1--35.

\bibitem[{Bouchard \protect\BIBand{} Warin(2012)}]{bouchard2012monte}
Bouchard B, Warin X (2012) Monte-carlo valuation of american options: facts and
  new algorithms to improve existing methods. \emph{Numerical methods in
  finance}, 215--255 (Springer).

\bibitem[{Boucherie \protect\BIBand{} Van~Dijk(2017)}]{boucherie2017markov}
Boucherie RJ, Van~Dijk NM (2017) \emph{Markov decision processes in practice}
  (Springer).

\bibitem[{Boukai(1990)}]{boukai1990explicit}
Boukai B (1990) An explicit expression for the distribution of the supremum of
  brownian motion with a change point. \emph{Communications in
  Statistics-Theory and Methods} 19(1):31--40.

\bibitem[{Brown et~al.(2010)Brown, Smith, \protect\BIBand{}
  Sun}]{brown2010information}
Brown DB, Smith JE, Sun P (2010) Information relaxations and duality in
  stochastic dynamic programs. \emph{Operations research}
  58(4-part-1):785--801.

\bibitem[{Chalasani et~al.(1999)Chalasani, Jha, \protect\BIBand{}
  Saias}]{chalasani1999approximate}
Chalasani P, Jha S, Saias I (1999) Approximate option pricing.
  \emph{Algorithmica} 25:2--21.

\bibitem[{Chawla et~al.(2020)Chawla, Gergatsouli, Teng, Tzamos,
  \protect\BIBand{} Zhang}]{chawla2020pandora}
Chawla S, Gergatsouli E, Teng Y, Tzamos C, Zhang R (2020) Pandora's box with
  correlations: Learning and approximation. \emph{2020 IEEE 61st Annual
  Symposium on Foundations of Computer Science (FOCS)}, 1214--1225 (IEEE).

\bibitem[{Chen \protect\BIBand{} Glasserman(2007)}]{chen2007additive}
Chen N, Glasserman P (2007) Additive and multiplicative duals for american
  option pricing. \emph{Finance and Stochastics} 11(2):153--179.

\bibitem[{Chen(2021)}]{chen2021efficient}
Chen Y (2021) \emph{Efficient Algorithms for High-Dimensional Data-Driven
  Sequential Decision-Making}. Ph.D. thesis, Cornell University.

\bibitem[{Chen \protect\BIBand{} Wang(2017)}]{chen2017lower}
Chen Y, Wang M (2017) Lower bound on the computational complexity of discounted
  markov decision problems. \emph{arXiv preprint arXiv:1705.07312} .

\bibitem[{Chow et~al.(1971)Chow, Robbins, \protect\BIBand{}
  Siegmund}]{chowgreat}
Chow Y, Robbins H, Siegmund D (1971) \emph{Great Expectations: The Theory of
  Optimal Stopping} (Houghton Mifflin, Boston).

\bibitem[{Christensen(2014)}]{christensen2014method}
Christensen S (2014) A method for pricing american options using semi-infinite
  linear programming. \emph{Mathematical Finance: An International Journal of
  Mathematics, Statistics and Financial Economics} 24(1):156--172.

\bibitem[{Ciocan \protect\BIBand{}
  Mi{\v{s}}i{\'c}(2022)}]{ciocan2020interpretable}
Ciocan DF, Mi{\v{s}}i{\'c} VV (2022) Interpretable optimal stopping.
  \emph{Management Science} 68(3):1616--1638.

\bibitem[{Cl{\'e}ment et~al.(2002)Cl{\'e}ment, Lamberton, \protect\BIBand{}
  Protter}]{clement2002analysis}
Cl{\'e}ment E, Lamberton D, Protter P (2002) An analysis of a least squares
  regression method for american option pricing. \emph{Finance and Stochastics}
  6(4):449--471.

\bibitem[{Cormen et~al.(2009)Cormen, Leiserson, Rivest, \protect\BIBand{}
  Stein}]{cormen2009introduction}
Cormen TH, Leiserson CE, Rivest RL, Stein C (2009) \emph{Introduction to
  algorithms} (MIT press).

\bibitem[{Correa et~al.(2019)Correa, Foncea, Hoeksma, Oosterwijk,
  \protect\BIBand{} Vredeveld}]{correa2019recent}
Correa J, Foncea P, Hoeksma R, Oosterwijk T, Vredeveld T (2019) Recent
  developments in prophet inequalities. \emph{ACM SIGecom Exchanges}
  17(1):61--70.

\bibitem[{Davis \protect\BIBand{} Karatzas(1994)}]{davis1994deterministic}
Davis MH, Karatzas I (1994) A deterministic approach to optimal stopping.
  \emph{Probability, Statistics and Optimisation (ed. FP Kelly). NewYork
  Chichester: John Wiley \& Sons Ltd} 455--466.

\bibitem[{De~Klerk(2008)}]{de2008complexity}
De~Klerk E (2008) The complexity of optimizing over a simplex, hypercube or
  sphere: a short survey. \emph{Central European Journal of Operations
  Research} 16(2):111--125.

\bibitem[{Denton et~al.(2009)Denton, Kurt, Shah, Bryant, \protect\BIBand{}
  Smith}]{denton2009optimizing}
Denton BT, Kurt M, Shah ND, Bryant SC, Smith SA (2009) Optimizing the start
  time of statin therapy for patients with diabetes. \emph{Medical Decision
  Making} 29(3):351--367.

\bibitem[{Derakhshan et~al.(2022)Derakhshan, Golrezaei, Manshadi,
  \protect\BIBand{} Mirrokni}]{derakhshan2022product}
Derakhshan M, Golrezaei N, Manshadi V, Mirrokni V (2022) Product ranking on
  online platforms. \emph{Management Science} 68(6):4024--4041.

\bibitem[{Desai et~al.(2012)Desai, Farias, \protect\BIBand{}
  Moallemi}]{desai2012pathwise}
Desai VV, Farias VF, Moallemi CC (2012) Pathwise optimization for optimal
  stopping problems. \emph{Management Science} 58(12):2292--2308.

\bibitem[{Du et~al.(2020)Du, Kakade, Wang, \protect\BIBand{} Yang}]{du2019good}
Du SS, Kakade SM, Wang R, Yang LF (2020) Is a good representation sufficient
  for sample efficient reinforcement learning? \emph{International Conference
  on Learning Representations}.

\bibitem[{Egloff et~al.(2005)}]{egloff2005monte}
Egloff D, et~al. (2005) Monte carlo algorithms for optimal stopping and
  statistical learning. \emph{The Annals of Applied Probability}
  15(2):1396--1432.

\bibitem[{Embrechts et~al.(2013)Embrechts, Kl{\"u}ppelberg, \protect\BIBand{}
  Mikosch}]{embrechts2013modelling}
Embrechts P, Kl{\"u}ppelberg C, Mikosch T (2013) \emph{Modelling extremal
  events: for insurance and finance}, volume~33 (Springer Science \& Business
  Media).

\bibitem[{Emmerling et~al.(2021)Emmerling, Yavas, \protect\BIBand{}
  Yildirim}]{emmerling2021accept}
Emmerling TJ, Yavas A, Yildirim Y (2021) To accept or not to accept: Optimal
  strategy for sellers in real estate. \emph{Real Estate Economics}
  49(S1):268--296.

\bibitem[{Freeman(1983)}]{freeman1983secretary}
Freeman P (1983) The secretary problem and its extensions: A review.
  \emph{International Statistical Review/Revue Internationale de Statistique}
  189--206.

\bibitem[{Gasull \protect\BIBand{} Utzet(2014)}]{gasull2014approximating}
Gasull A, Utzet F (2014) Approximating mills ratio. \emph{Journal of
  Mathematical Analysis and Applications} 420(2):1832--1853.

\bibitem[{Glasserman \protect\BIBand{} Merener(2003)}]{glasserman2003cap}
Glasserman P, Merener N (2003) Cap and swaption approximations in libor market
  models with jumps. \emph{Journal of Computational Finance} 7(1):1--36.

\bibitem[{Glasserman et~al.(2004)Glasserman, Yu et~al.}]{glasserman2004number}
Glasserman P, Yu B, et~al. (2004) Number of paths versus number of basis
  functions in american option pricing. \emph{The Annals of Applied
  Probability} 14(4):2090--2119.

\bibitem[{Gordon(1941)}]{gordon1941values}
Gordon RD (1941) Values of mills' ratio of area to bounding ordinate and of the
  normal probability integral for large values of the argument. \emph{The
  Annals of Mathematical Statistics} 12(3):364--366.

\bibitem[{Halman et~al.(2014)Halman, Klabjan, Li, Orlin, \protect\BIBand{}
  Simchi-Levi}]{halman2014fully}
Halman N, Klabjan D, Li CL, Orlin J, Simchi-Levi D (2014) Fully polynomial time
  approximation schemes for stochastic dynamic programs. \emph{SIAM Journal on
  Discrete Mathematics} 28(4):1725--1796.

\bibitem[{Halman \protect\BIBand{} Nannicini(2019)}]{halman2019toward}
Halman N, Nannicini G (2019) Toward breaking the curse of dimensionality: an
  fptas for stochastic dynamic programs with multidimensional actions and
  scalar states. \emph{SIAM Journal on Optimization} 29(2):1131--1163.

\bibitem[{Halman et~al.(2015)Halman, Nannicini, \protect\BIBand{}
  Orlin}]{halman2015computationally}
Halman N, Nannicini G, Orlin J (2015) A computationally efficient fptas for
  convex stochastic dynamic programs. \emph{SIAM Journal on Optimization}
  25(1):317--350.

\bibitem[{Hamida \protect\BIBand{} Cont(2005)}]{hamida2005recovering}
Hamida SB, Cont R (2005) Recovering volatility from option prices by
  evolutionary optimization. \emph{The Journal of Computational Finance} .

\bibitem[{Haugh \protect\BIBand{} Kogan(2004)}]{haugh2004pricing}
Haugh MB, Kogan L (2004) Pricing american options: a duality approach.
  \emph{Operations Research} 52(2):258--270.

\bibitem[{Hill \protect\BIBand{} Kertz(1983)}]{hill1983stop}
Hill TP, Kertz RP (1983) Stop rule inequalities for uniformly bounded sequences
  of random variables. \emph{Transactions of the American Mathematical Society}
  278(1):197--207.

\bibitem[{Hill \protect\BIBand{} Kertz(1992)}]{hill1992survey}
Hill TP, Kertz RP (1992) A survey of prophet inequalities in optimal stopping
  theory. \emph{Contemp. Math} 125:191--207.

\bibitem[{Hin(2015)}]{hin2015analysis}
Hin LY (2015) \emph{Analysis and Modelling of Implied Market Parameters}. Ph.D.
  thesis, Curtin University.

\bibitem[{Ib{\'a}{\~n}ez \protect\BIBand{} Velasco(2020)}]{ibanez2020recursive}
Ib{\'a}{\~n}ez A, Velasco C (2020) Recursive lower and dual upper bounds for
  bermudan-style options. \emph{European Journal of Operational Research}
  280(2):730--740.

\bibitem[{Jamshidian(2007)}]{jamshidian2007duality}
Jamshidian F (2007) The duality of optimal exercise and domineering claims: A
  doob--meyer decomposition approach to the snell envelope. \emph{Stochastics
  An International Journal of Probability and Stochastic Processes}
  79(1-2):27--60.

\bibitem[{Kakade et~al.(2003)}]{kakade2003sample}
Kakade SM, et~al. (2003) \emph{On the sample complexity of reinforcement
  learning}. Ph.D. thesis, University of London London, England.

\bibitem[{Kisfaludi-Bak et~al.(2022)Kisfaludi-Bak, Nederlof, \protect\BIBand{}
  Wegrzycki}]{kisfaludi2022gap}
Kisfaludi-Bak S, Nederlof J, Wegrzycki K (2022) A gap-eth-tight approximation
  scheme for euclidean tsp. \emph{2021 IEEE 62nd Annual Symposium on
  Foundations of Computer Science (FOCS)}, 351--362 (IEEE).

\bibitem[{Kohler(2008)}]{kohler2008regression}
Kohler M (2008) A regression-based smoothing spline monte carlo algorithm for
  pricing american options in discrete time. \emph{AStA Advances in Statistical
  Analysis} 92(2):153--178.

\bibitem[{Kohler et~al.(2010)Kohler, Krzy{\.z}ak, \protect\BIBand{}
  Todorovic}]{kohler2010pricing}
Kohler M, Krzy{\.z}ak A, Todorovic N (2010) Pricing of high-dimensional
  american options by neural networks. \emph{Mathematical Finance: An
  International Journal of Mathematics, Statistics and Financial Economics}
  20(3):383--410.

\bibitem[{Kolodko \protect\BIBand{} Schoenmakers(2004)}]{kolodko2004efficient}
Kolodko A, Schoenmakers J (2004) An efficient dual monte carlo upper bound for
  bermudan style derivative .

\bibitem[{Kolodko \protect\BIBand{} Schoenmakers(2006)}]{kolodko2006iterative}
Kolodko A, Schoenmakers J (2006) Iterative construction of the optimal bermudan
  stopping time. \emph{Finance and Stochastics} 10(1):27--49.

\bibitem[{Kou \protect\BIBand{} Wang(2004)}]{kou2004option}
Kou SG, Wang H (2004) Option pricing under a double exponential jump diffusion
  model. \emph{Management science} 50(9):1178--1192.

\bibitem[{Krotofil et~al.(2014)Krotofil, Cardenas, Larsen, \protect\BIBand{}
  Gollmann}]{krotofil2014vulnerabilities}
Krotofil M, Cardenas A, Larsen J, Gollmann D (2014) Vulnerabilities of
  cyber-physical systems to stale data—determining the optimal time to launch
  attacks. \emph{International journal of critical infrastructure protection}
  7(4):213--232.

\bibitem[{Lai \protect\BIBand{} Wong(2004)}]{lai2004valuation}
Lai TL, Wong SS (2004) Valuation of american options via basis functions.
  \emph{IEEE transactions on automatic control} 49(3):374--385.

\bibitem[{Langren{\'e} et~al.(2017)Langren{\'e}, Chen, \protect\BIBand{}
  Zhu}]{langrene2017field}
Langren{\'e} N, Chen W, Zhu Z (2017) Field exploration: when to start
  extracting? \emph{Proceedings of the 22nd International Congress on Modelling
  and Simulation (MODSIM 2017), Hobart, Australia}.

\bibitem[{Le{\~a}o et~al.(2019)Le{\~a}o, Ohashi, \protect\BIBand{}
  Russo}]{leao2019discrete}
Le{\~a}o D, Ohashi A, Russo F (2019) Discrete-type approximations for
  non-markovian optimal stopping problems: Part i. \emph{Journal of Applied
  Probability} 56(4):981--1005.

\bibitem[{Lelong(2018)}]{lelong2018dual}
Lelong J (2018) Dual pricing of american options by wiener chaos expansion.
  \emph{SIAM Journal on Financial Mathematics} 9(2):493--519.

\bibitem[{Longstaff \protect\BIBand{} Schwartz(2001)}]{longstaff2001valuing}
Longstaff FA, Schwartz ES (2001) Valuing american options by simulation: a
  simple least-squares approach. \emph{The review of financial studies}
  14(1):113--147.

\bibitem[{Nemirovski \protect\BIBand{} Yudin(1983)}]{nemirovskij1983problem}
Nemirovski AS, Yudin DB (1983) Problem complexity and method efficiency in
  optimization .

\bibitem[{Rockafellar \protect\BIBand{} Wets(1991)}]{rockafellar1991scenarios}
Rockafellar RT, Wets RJB (1991) Scenarios and policy aggregation in
  optimization under uncertainty. \emph{Mathematics of operations research}
  16(1):119--147.

\bibitem[{Rogers(2015)}]{rogers2015bermudan}
Rogers L (2015) Bermudan options by simulation. \emph{arXiv preprint
  arXiv:1508.06117} .

\bibitem[{Rogers(2002)}]{rogers2002monte}
Rogers LC (2002) Monte carlo valuation of american options. \emph{Mathematical
  Finance} 12(3):271--286.

\bibitem[{Rubinstein et~al.(2020)Rubinstein, Wang, \protect\BIBand{}
  Weinberg}]{rubinstein2020optimal}
Rubinstein A, Wang JZ, Weinberg SM (2020) Optimal single-choice prophet
  inequalities from samples. \emph{Innovations in Theoretical Computer Science}
  .

\bibitem[{Sabino(2007)}]{sabino2007monte}
Sabino P (2007) Monte carlo methods and path-generation techniques for pricing
  multi-asset path-dependent options. \emph{arXiv preprint arXiv:0710.0850} .

\bibitem[{Santaniello et~al.(2011)Santaniello, Burns, Golby, Singer, Anderson,
  \protect\BIBand{} Sarma}]{santaniello2011quickest}
Santaniello S, Burns SP, Golby AJ, Singer JM, Anderson WS, Sarma SV (2011)
  Quickest detection of drug-resistant seizures: An optimal control approach.
  \emph{Epilepsy \& Behavior} 22:S49--S60.

\bibitem[{Santaniello et~al.(2012)Santaniello, Sherman, Thakor, Eskandar,
  \protect\BIBand{} Sarma}]{santaniello2012optimal}
Santaniello S, Sherman DL, Thakor NV, Eskandar EN, Sarma SV (2012) Optimal
  control-based bayesian detection of clinical and behavioral state
  transitions. \emph{IEEE transactions on neural systems and rehabilitation
  engineering} 20(5):708--719.

\bibitem[{Schalekamp(2007)}]{schalekamp2007some}
Schalekamp F (2007) \emph{Some Results In Universal And A Priori Optimization}.
  Ph.D. thesis, Cornell University, Ph.D. Thesis.

\bibitem[{Segev \protect\BIBand{} Singla(2021)}]{segev2021efficient}
Segev D, Singla S (2021) Efficient approximation schemes for stochastic probing
  and prophet problems. \emph{Proceedings of the 22nd ACM Conference on
  Economics and Computation}, 793--794.

\bibitem[{Shapiro \protect\BIBand{} Nemirovski(2005)}]{shapiro2005complexity}
Shapiro A, Nemirovski A (2005) On complexity of stochastic programming
  problems. \emph{Continuous optimization}, 111--146 (Springer).

\bibitem[{Sidford et~al.(2018)Sidford, Wang, Wu, Yang, \protect\BIBand{}
  Ye}]{sidford2018near}
Sidford A, Wang M, Wu X, Yang L, Ye Y (2018) Near-optimal time and sample
  complexities for solving markov decision processes with a generative model.
  \emph{Advances in Neural Information Processing Systems}, 5186--5196.

\bibitem[{Solan et~al.(2012)Solan, Tsirelson, \protect\BIBand{}
  Vieille}]{solan2012random}
Solan E, Tsirelson B, Vieille N (2012) Random stopping times in stopping
  problems and stopping games. \emph{arXiv preprint arXiv:1211.5802} .

\bibitem[{Spielman \protect\BIBand{} Teng(2004)}]{spielman2004smoothed}
Spielman DA, Teng SH (2004) Smoothed analysis of algorithms: Why the simplex
  algorithm usually takes polynomial time. \emph{Journal of the ACM (JACM)}
  51(3):385--463.

\bibitem[{Stentoft(2004)}]{stentoft2004convergence}
Stentoft L (2004) Convergence of the least squares monte carlo approach to
  american option valuation. \emph{Management Science} 50(9):1193--1203.

\bibitem[{Sturt(2023)}]{sturt2021nonparametric}
Sturt B (2023) A nonparametric algorithm for optimal stopping based on robust
  optimization. \emph{Operations Research} 71(5):1530--1557.

\bibitem[{Swamy \protect\BIBand{} Shmoys(2012)}]{swamy2012sampling}
Swamy C, Shmoys DB (2012) Sampling-based approximation algorithms for
  multistage stochastic optimization. \emph{SIAM Journal on Computing}
  41(4):975--1004.

\bibitem[{Traub \protect\BIBand{} Werschulz(1998)}]{traub1998complexity}
Traub JF, Werschulz AG (1998) \emph{Complexity and information}, volume 26862
  (Cambridge University Press).

\bibitem[{Tsitsiklis \protect\BIBand{}
  Van~Roy(2001)}]{tsitsiklis2001regression}
Tsitsiklis JN, Van~Roy B (2001) Regression methods for pricing complex
  american-style options. \emph{IEEE Transactions on Neural Networks}
  12(4):694--703.

\bibitem[{Wang \protect\BIBand{} Han(2015)}]{wang2015basics}
Wang D, Han Z (2015) Basics for sublinear algorithms. \emph{Sublinear
  Algorithms for Big Data Applications}, 9--21 (Springer).

\bibitem[{Wang et~al.(2021)Wang, Wang, \protect\BIBand{}
  Kakade}]{wang2021exponential}
Wang Y, Wang R, Kakade S (2021) An exponential lower bound for linearly
  realizable mdp with constant suboptimality gap. \emph{Advances in Neural
  Information Processing Systems} 34:9521--9533.

\bibitem[{Williamson \protect\BIBand{} Shmoys(2011)}]{williamson2011design}
Williamson DP, Shmoys DB (2011) \emph{The design of approximation algorithms}
  (Cambridge university press).

\bibitem[{Wright(2005)}]{wright2005interior}
Wright M (2005) The interior-point revolution in optimization: history, recent
  developments, and lasting consequences. \emph{Bulletin of the American
  mathematical society} 42(1):39--56.

\bibitem[{Zhang et~al.(2012)Zhang, Denton, Balasubramanian, Shah,
  \protect\BIBand{} Inman}]{zhang2012optimization}
Zhang J, Denton BT, Balasubramanian H, Shah ND, Inman BA (2012) Optimization of
  prostate biopsy referral decisions. \emph{Manufacturing \& Service Operations
  Management} 14(4):529--547.

\bibitem[{Zhou et~al.(2023)Zhou, Wang, Blanchet, \protect\BIBand{}
  Glynn}]{zhou2021unbiased}
Zhou Z, Wang G, Blanchet JH, Glynn PW (2023) Unbiased optimal stopping via the
  muse. \emph{Stochastic Processes and their Applications} 166:104088.

\end{thebibliography}

\ECSwitch


\ECHead{Online Appendix}
In this appendix, we include a broader discussion and motivation of the notion of polynomial-time efficiency for a fixed accuracy $\epsilon$, additional intuition for our main results, several auxiliary proofs, as well as the formal descriptions and analysis of all algorithms.  We now provide a more detailed outline.  In Appendix\ \ref{sec:ecptas}, we provide a more detailed discussion and broader context for algorithms with polynomial complexity for fixed accuracy (and e.g. the notion of a polynomial-time approximation scheme).  In Appendix \ref{sec:ecmoreintuit}, we provide additional intuition regarding the complexity bounds for our algorithms.  In Appendix\ \ref{sec:ectheoryproof}, we provide the formal proofs of several of our theoretical results.  In Appendix\ \ref{sec:ecalgo}, we provide formal descriptions of our main algorithms, and the proofs of our main algorithmic results Theorems\ \ref{thm:mainalgo2}, \ref{thm:mainpolicy1}, \ref{thm:mainmax}, \ref{thm:mainmaxpolicy}.  In Appendix\ \ref{sec:ecmaxcall}, we prove Corollary\ \ref{cormaxcall}, showing that the Bermudan max call satisfies our assumptions.  In Appendix\ \ref{sec:ecnonneg}, we include the short proof of the observation that assuming the reward process to be non-negative is w.l.o.g. (as also observed recently in \citet{baas2023sampling}).

\section{Broader discussion of algorithms with polynomial complexity for fixed accuracy $\epsilon$}\label{sec:ecptas}
Be it in the traditional bit complexity model or the information complexity model, it is reasonable to ask why it makes sense to first fix the accuracy $\epsilon$ and then ask how the complexity scales in terms of the other parameters.  In either framework, an algorithm which has polynomial runtime for any fixed accuracy $\epsilon$ is typically called a polynomial-time approximation scheme (PTAS), where the counterpart which is given access to its own random coin tosses (e.g. to implement a randomized algorithm) is called a PRAS.  Such algorithms play a key role in the formal study of many problems in operations research, computer science, and operations management (\cite{williamson2011design}).  The question of whether or not a PTAS exists for any given optimization problem has become central to understanding the complexity of these problems across multiple academic communities, and plays a significant role in the broader subject of computational complexity (\cite{arora2009computational}).  
\\\indent For several recent operations management / management science papers utilizing this concept, see for example \cite{alaei2022revenue, aouad2021display, derakhshan2022product}.  Analogous concepts in the setting of stochastic dynamic optimization have also been widely studied across operations research and computer science, see e.g. \cite{swamy2012sampling, halman2019toward,kakade2003sample}, including in problems closely related to optimal stopping (\cite{segev2021efficient}).  For such stochastic problems, the even weaker notion of an algorithm which produces a solution within a constant factor of optimal in polynomial time is also commonly studied, with many relevant works in the context of optimal stopping and prophet inequalities (\cite{rubinstein2020optimal,correa2019recent}).  
\\\indent That said, the concept of PTAS remains largely unstudied in stochastic dynamic settings with complex and high-dimensional dependencies over time.  Indeed, this gap in the literature was recently pointed out in \cite{chawla2020pandora}, which analyzed constant-factor approximation algorithms for such problems related to the celebrated ``Pandora's box" problem.  Similarly, the concept of PTAS has not yet been formally studied in the context of high-dimensional path-dependent optimal stopping.  Our work takes a first step along these lines, and we view bridging the formalism of PTAS with fundamental problems in high-dimensional stochastic control arising in operations management and finance as another contribution of our work.  
\\\indent Let us also note that although PTAS is a standard framework for assessing approximation algorithms in computer science and operations research, other frameworks have been proposed by various academic communities.  In certain settings PTAS may not be the most appropriate measure of the complexity of computing approximately optimal solutions.  For example, \cite{belomestny2019semi} introduced the concept of semi-tractability, which requires an algorithm to return an $\epsilon$ approximation with a sample complexity growing slower than $\epsilon^{-D}$ as $\epsilon \to 0$ and $D \to \infty,$ and exhibits stochastic mesh methods with such a property under additional continuity assumptions.  
\\\indent Such alternative notions of tractability may be especially relevant in settings where the discrete-time problem arises as a discretization of an underlying continuous problem (as in several financial models for options pricing).  In general, the runtime of PTAS typically degrades as $\epsilon \downarrow 0$, and thus such methods may not be appropriate in settings where high accuracy is required.  None-the-less, as discussed previously, the question of whether or not a PTAS exists for any given optimization problem has become central to understanding the complexity of these problems across multiple academic communities, and plays a significant role in the broader subject of computational complexity (\cite{arora2009computational}).   Furthermore, let us point out that in the historical development of algorithms in operations research and computer science, there are several cases of important problems (perhaps most notably the Euclidean Traveling Salesman Problem) for which the first PTAS was not practically efficient  (especially for small $\epsilon$), but whose development prompted subsequent work leading to much more efficient algorithms (\cite{arora1996polynomial, arora1997nearly}).  Similarly, there are several cases of well-known problems for which algorithms were first developed, but only later were a ``nice" set of assumptions balancing practicality with theoretical tractability discovered, leading to a fast runtime for the known algorithm under reasonable conditions relevant to practice (\cite{spielman2004smoothed,abraham2010highway,beier2004random}).  Also, more broadly in the theory of algorithms, there are several high-profile examples (e.g. in linear programming) where a first polynomial-time algorithm for a problem was not practically efficient, but spurred the development of more efficient methods (\cite{wright2005interior}).

\section{Further discussion and interpretation of complexity bounds}\label{sec:ecmoreintuit}
Here we provide additional intuition behind the complexity bound appearing in Theorem\ \ref{thm:mainalgo2}.  For concreteness, let us first describe in some depth how our algorithm computes $\hat{L}_2$, our estimate of $L_2$ (the second term in our expansion).  To simulate $\hat{L}_2$, we repeat the following many (say $N$) times.  First, generate a trajectory $Y^i_{[T]}$ of the information process.  Second, for each of the $T$ values $Y^i_1,Y^i_2,\ldots,Y^i_T$ of the information process (on that trajectory), compute an approximation $\hat{Z}^2_t$ to $Z^2_t$ as follows.  Step A: compute $Z^1_t$ exactly (for trajectory $Y^i_{[T]}$), as this is simply the value of the reward process on that trajectory at time $t$.  Step B: use ${\mathcal B}$ to generate many trajectories (say again $N$) of the information process (conditioned on the information up to time $t$ being exactly $Y^i_{[t]}$) to accurately estimate (by simulation) the value of the conditional expectation $E\big[\max_{j = 1,\dots, T} Z^1_j \big|\mathcal{F}_t\big]$ (for trajectory $Y^i_{[T]}$).  For each of these $N$ ``inner trajectories", one must compute the maximum of $T$ terms (which, for each inner trajectory, takes effort roughly $T$).  Step C: set $\hat{Z}^2_t$ (for the single trajectory $Y^i_{[T]}$) to be the difference between the exactly computed value $Z^1_t$ and the estimated (by simulation) value for $E\big[\max_{j = 1,\dots, T} Z^1_j \big|\mathcal{F}_t\big].$  Thus the total effort needed to estimate each $\hat{Z}^2_t$ term is (roughly) $N \times T$.  
\\\indent Note that for each single trajectory $Y^i_{[T]}$, we must conduct this procedure to estimate all $T$ of the values $\hat{Z}^2_1, \hat{Z}^2_2,\ldots,\hat{Z}^2_T$ (thus requiring roughly $N \times T^2$ effort). Then for each single trajectory $Y^i_{[T]}$, the algorithm computes $\max_{t=1,\ldots,T} \hat{Z}^2_t$.  Finally, the algorithm averages this over all $N$ trajectories to yield its estimate $\hat{L}_2$.  As for each of these $N$ trajectories the algorithm must use roughly $N \times T^2$ effort, the overall effort to estimate $\hat{L}_2$ is roughly $(N T)^2$.  The key insight of the above is that computing $\hat{L}_2$ required roughly $(N T)^2$ effort.  Extending a very similar logic recursively to the computation of $\hat{L}_i$ for $i \geq 3$, one concludes that roughly $(N T)^i$ effort is required.  Let us note that for $i \geq 3$ each inner simulation itself requires further levels of inner simulation etc. (until one gets to the ``base level" in which the actual reward values are used), but the same general logic still applies to yield an overall effort of roughly $(N T)^i$. As to yield an $\epsilon$-approximation we must truncate our expansion after $O(\frac{1}{\epsilon})$ terms, this would yield a rough estimated runtime of $(N T)^{O(\frac{1}{\epsilon})}$.  This is indeed precisely where the $T^{O(\frac{1}{\epsilon})}$ term comes from in our complexity analysis.  The term involving $C + G$ comes from the cost of each individual call to the simulator and evaluation of a reward function, and the $\log(\delta^{-1})$ term comes from standard results in probability and the analysis of algorithms which show that (under appropriate assumptions) re-computing an estimator roughly $\log(\delta^{-1})$ independent times is sufficient to yield an estimator with confidence $1 - \delta$.  
\\\indent Finally, let us comment in greater depth on the expression $\exp(100 \epsilon^{-2}).$  As noted in the discussion in Section\ \ref{moreintuitmain}, intuitively, the accuracy required in simulations at nested depth $k+1$ is double that required at depth $k$.  Furthermore, by standard results in probability, achieving an accuracy of $\epsilon$ requires roughly $\epsilon^{-2}$ replications (ignoring other considerations such as the desired confidence).  Informally letting ${\mathfrak f}_k(\epsilon)$ denote a rough estimate for the complexity requried to carry out the required simulations at depth $k$ of recursion to within accuracy $\epsilon$, one is thus led to the recursion ${\mathfrak f}_{k+1}(\epsilon) = \epsilon^{-2} {\mathfrak f}_{k}(\frac{\epsilon}{2}).$  As our algorithm must go out to recursive depth $O(\frac{1}{\epsilon}),$ one is led to the rough complexity bound 
$$ (\frac{1}{\epsilon})^2 \times (\frac{2}{\epsilon})^2 \times (\frac{4}{\epsilon})^2 \ldots \times (\frac{2^{\frac{1}{\epsilon}}}{\epsilon})^2\ \ =\ \ (\frac{1}{\epsilon})^{\frac{2}{\epsilon}} \times 2^{2 \times \sum_{i=1}^{\epsilon^{-1}} i}\ \ = 2^{O(\epsilon^{-2})}.
$$
For further discussion of $N$ (and its dependence on $\epsilon$), and the complexity bound appearing in Theorem\ \ref{thm:mainmax}, we refer the reader to the discussion in Section\ \ref{moreintuitmain}.  For the formal derivation of our algorithm's complexity, we refer the reader to Appendix\ \ref{sec:ecalgo}.

\section{Omitted proofs of theoretical results}\label{sec:ectheoryproof}
In this section, we provide the proofs of several of our theoretical results, namely Theorem \ref{thm:main1} and Lemma\ \ref{lem:hardinstance}.
\subsection{Expansion and Proof of Lemma\ \ref{lem:intuit1} and Theorem\ \ref{thm:main1}}
Recall that for $k \geq 1$ and $t \in \lbrace 1,\ldots,T \rbrace$, $Z^{k+1}_t = Z^k_t - E\big[ \max_{i = 1,\ldots,T} Z^k_i | {\mathcal F}_t \big]$. We begin by proving Lemma\ \ref{lem:intuit1}, which formalizes some of the key properties of our expansion.
\begin{proof}{Proof of Lemma\ \ref{lem:intuit1}}
First, let us prove integrability.  We proceed by induction.  The base case ($k = 1$) follows from our assumption that the reward process is integrable.  Suppose the induction is true for all $j \leq k$ for some $k \geq 1$.  Recall that by definition $Z^{k+1}_t = Z^k_t - E[\max_{i = 1, \dots, T} Z^k_i | \mathcal{F}_t].$  It follows that w.p.1, $|Z^{k+1}_t| \leq |Z^{k}_t| + E[ \sum_{i=1}^T |Z^{k}_i| \big| {\mathcal F}_t] =   |Z^{k}_t| + \sum_{i=1}^T E\big[ |Z^{k}_i| \big| {\mathcal F}_t \big]$, where we note that all conditional expectations are finite and well-defined by the integrability of $Z^k_t$ for all $t$ and the basic properties of conditional expectation (e.g. the Radon-Nikodym Theorem).  By the tower property of conditional expectation, it follows that $E[|Z^{k+1}_t|] \leq E[|Z^{k}_t|] + \sum_{i=1}^T E[|Z^{k}_i|]$, and the desired integrability (and induction) then follows from the induction hypothesis.
\\\\Next, let us prove non-positivity.  First we treat the base case $k = 2$.  Note that $Z^2_t = Z^1_t - E[\max_{i=1,\ldots,T} Z^1_i | {\mathcal F}_t]$.  Since $\max_{i=1,\ldots,T} Z^1_i \geq Z^1_t$, we conclude from the basic properties of conditional expectation that $Z^2_t \leq Z^1_t - E[Z^1_t | {\mathcal F}_t]$.  As $Z^1_t$ is ${\mathcal F}_t$-measurable, the basic properties of conditional expectation imply that $E[Z^1_t | {\mathcal F}_t] = Z^1_t$.  We conclude that $Z^2_t \leq Z^1_t - Z^1_t = 0$, completing the proof.  Suppose the induction is true for all $j \in \lbrace 2,\ldots,k \rbrace$ for some $k \geq 2$.  The induction step follows essentially identically to that of the base case, and we omit the details.  Combining the above completes the proof.
\\\\Next, let us prove monotonicity.  By definition, $Z^{k+1}_t = Z^k_t - E[\max_{i = 1, \dots, T} Z^k_i | \mathcal{F}_t].$  By the demonstrated non-positivity and the basic properties of conditional expectation, $E[\max_{i = 1, \dots, T} Z^k_i | \mathcal{F}_t] \leq 0$.  It follows that $Z^{k+1}_t \geq Z^k_t$ for all $k \geq 2$.
\\\\Next, let us prove adaptedness.  We proceed by induction.  The base case $(k = 1$) follows from our assumption that the reward process is adapted.  Suppose the induction is true for all $j \leq k$ for some $k \geq 1$.  Recall that by definition $Z^{k+1}_t = Z^k_t - E[\max_{i = 1, \dots, T} Z^k_i | \mathcal{F}_t].$  $Z^k_t$ is adapted by the induction hypothesis, and $E[\max_{i = 1, \dots, T} Z^k_i | \mathcal{F}_t]$ is adapted by the basic properties of conditional expectation.  As the difference of adapted functions is adapted, this completes the induction, and the proof.
\\\\Finally, let us prove that $\textsc{opt} = \sum_{i=1}^k E[\max_{t = 1,\ldots,T} Z^i_t] + \sup_{\tau \in {\mathcal T}} E[Z^{k+1}_{\tau}].$  We proceed by induction.  The base case $k = 0$ follows from the definition of the OS problem.  Suppose the induction is true for all $j \leq k$ for some $k \geq 1$.  Thus $\textsc{opt} = \sum_{i=1}^k E[\max_{t = 1,\ldots,T} Z^i_t] + \sup_{\tau \in {\mathcal T}} E[Z^{k+1}_{\tau}].$  By the optional stopping theorem, for all $\tau \in {\mathcal T}$, it holds that 
\begin{eqnarray*}
E[Z^{k+1}_{\tau}] &=& E\big[ Z^{k+1}_{\tau} - E[\max_{i=1,\ldots,T} Z^{k+1}_i | {\mathcal F}_{\tau}] \big] + E[\max_{i=1,\ldots,T} Z^{k+1}_i]
\\&=& E[ Z^{k+2}_{\tau} ] + E[\max_{i=1,\ldots,T} Z^{k+1}_i].
\end{eqnarray*}
We conclude that 
\begin{eqnarray*}
\sup_{\tau \in {\mathcal T}} E[Z^{k+1}_{\tau}] &=& \sup_{\tau \in {\mathcal T}} \big( E[ Z^{k+2}_{\tau} ] + E[\max_{i=1,\ldots,T} Z^{k+1}_i] \big)
\\&=& E[\max_{i=1,\ldots,T} Z^{k+1}_i] + \sup_{\tau \in {\mathcal T}} E[ Z^{k+2}_{\tau} ].
\end{eqnarray*}
Applying the induction hypothesis, we conclude that 
\begin{eqnarray*}
\textsc{opt} &=& \sum_{i=1}^k E[\max_{t = 1,\ldots,T} Z^i_t] + \sup_{\tau \in {\mathcal T}} E[Z^{k+1}_{\tau}]
\\&=& \sum_{i=1}^k E[\max_{t = 1,\ldots,T} Z^i_t] + E[\max_{i=1,\ldots,T} Z^{k+1}_i] + \sup_{\tau \in {\mathcal T}} E[ Z^{k+2}_{\tau} ]
\\&=& \sum_{i=1}^{k+1} E[\max_{t = 1,\ldots,T} Z^i_t] + \sup_{\tau \in {\mathcal T}} E[ Z^{k+2}_{\tau} ],
\end{eqnarray*}
completing the proof. \qed
\end{proof}
With Lemma \ref{lem:intuit1}, we would be done (at least with the proof of Theorem\ \ref{thm:main1}) if $\lim_{k \rightarrow \infty} \sup_{\tau \in {\mathcal T}} E[Z^{k+1}_{\tau}] = 0.$ We now prove that this is indeed the case.   
\proof{Proof of Theorem\ \ref{thm:main1}:}
It follows from Lemma \ref{lem:intuit1} and monotone convergence that $\lbrace \mathbf{Z}^k, k \geq 1 \rbrace$ converges almost surely (a.s.), and thus $\lbrace Z^{k+1}_T - Z^k_T, k \geq 1 \rbrace$ converges a.s. to 0.  By definition, $Z^{k+1}_T = Z^k_T - \max_{i = 1,\ldots,T} Z^k_i,$ and it follows that $\lbrace \max_{i = 1,\ldots,T} Z^k_i , k \geq 1 \rbrace$ converges a.s. to 0, i.e. $\lim_{k \to \infty} \max_{t = 1,\ldots,T} Z^k_t  = 0 \ a.s.$  We now argue that, combined with the monotonicity and non-positivity of $\lbrace Z^k_t, k \geq 2 \rbrace$, we may interchange limit and max to conclude that $\max_{t = 1,\ldots,T} \lim_{k \to \infty} Z^k_t = 0$ a.s.  Indeed, suppose for contradiction that there is some positive measure event $A$ on which 
$\max_{t = 1,\ldots,T} \lim_{k \to \infty} Z^k_t \neq 0$.  By non-positivity, it must be that on the event $A$ it holds that $\max_{t = 1,\ldots,T} \lim_{k \to \infty} Z^k_t < 0,$ and thus $\lim_{k \rightarrow \infty} Z^k_t < 0$ for all $t \in \lbrace 1,\ldots,T \rbrace$ (on this event $A$).  By the monotonicity of $\lbrace Z^k_t, k \geq 2 \rbrace$ (proven in Lemma\ \ref{lem:intuit1}), we conclude that (on this event $A$) there exist finite $K$ and $\epsilon > 0$ s.t. $Z^k_t < - \epsilon$ for all $k \geq K$ and all $t \in \lbrace 1,\ldots,T \rbrace$, and thus 
$\max_{t = 1,\ldots,T} Z^k_t < - \epsilon$ for all $k \geq K$.  But this is a contradiction, since we have already proven that $\lim_{k \to \infty} \max_{t = 1,\ldots,T} Z^k_t  = 0 \ a.s.$  We conclude that such an interchange is indeed justified, and thus $\max_{t = 1,\ldots,T} \lim_{k \to \infty} Z^k_t = 0$ a.s.  In other words, \begin{equation}\label{existneed1}
P\big( \exists\ t\ \in \lbrace 1,\ldots,T \rbrace\ \textrm{such that} \lim_{k \rightarrow \infty} Z^k_t = 0 \big) = 1.
\end{equation}
Now we introduce a stopping time $\tau_\infty$. Let $\tau_\infty \stackrel{\Delta}{=} \min\{t \in \lbrace 1,\ldots,T \rbrace: \lim_{k \to \infty} Z^k_t = 0\}$, where (\ref{existneed1}) ensures $\tau_{\infty}$ is well-defined, and the adaptedness of $\mathbf{Z}^k$ ensures $\tau_{\infty}$ is a stopping time.  It follows that 
$\lbrace Z^k_{{\tau_{\infty}}}, k \geq 2 \rbrace$ is a sequence of non-positive r.v.s converging a.s. to 0.  Indeed, this follows from the non-positivity of $Z^k_t$ for $k \geq 2$, and the fact that $Z^k_{{\tau_{\infty}}} = \sum_{t=1}^T Z^k_t I(\tau_{\infty} = t)$, and 
\begin{eqnarray*}
\lim_{k \rightarrow \infty} \sum_{t=1}^T Z^k_t I(\tau_{\infty} = t) &=& \sum_{t=1}^T I(\tau_{\infty} = t) \lim_{k \rightarrow \infty} Z^k_t 
\\&=& \sum_{t=1}^T I(\tau_{\infty} = t) \times 0 \ \ \textrm{since}\ \ I(\tau_{\infty} = t) = 1 \rightarrow \lim_{k \rightarrow \infty} Z^k_t = 0
\\&=& 0.
\end{eqnarray*}
Since $|Z^k_{{\tau_{\infty}}}| \leq \sum_{t=1}^T |Z^k_t| \leq \sum_{t=1}^T |Z^2_t|$ (by the non-positivity and monotonicity shown in Lemma\ \ref{lem:intuit1}), 
and $\sum_{t=1}^T |Z^2_t|$ is integrable, we may thus apply dominated convergence to conclude that $\lim_{k \to \infty} E[Z^k_{\tau_\infty}] = E\big[\lim_{k \to \infty} Z^k_{\tau_\infty}\big] = 0$.  Now, for each $k \geq 2$, note that $0 \geq \sup_{\tau \in \mathcal{T}} E[Z^k_{\tau}]$ by the non-positivity of $Z^k_t, t = 1,\ldots,T$, and thus $0 \geq \limsup_{k \rightarrow \infty} \sup_{\tau \in \mathcal{T}} E[Z^k_{\tau}].$  But as the supremum over a set of stopping times is at least what one achieves with any given stopping time, we also have that (for all $k \geq 2$) $\sup_{\tau \in \mathcal{T}} E[Z^k_{\tau}] \geq E[Z^k_{\tau_\infty}]$, and thus $\liminf_{k \rightarrow \infty} \sup_{\tau \in \mathcal{T}} E[Z^k_{\tau}] \geq \liminf_{k \rightarrow \infty} E[Z^k_{\tau_\infty}].$  But as we have already shown that $\lim_{k \rightarrow \infty} E[Z^k_{\tau_\infty}] = 0$, we conclude that $\liminf_{k \rightarrow \infty} E[Z^k_{\tau_\infty}] = 0$, and $\liminf_{k \rightarrow \infty} \sup_{\tau \in \mathcal{T}} E[Z^k_{\tau}] \geq 0$.  Combining the above, we find that
$$0 \geq \limsup_{k \rightarrow \infty} \sup_{\tau \in \mathcal{T}} E[Z^k_{\tau}] \geq \liminf_{k \rightarrow \infty} \sup_{\tau \in \mathcal{T}} E[Z^k_{\tau}] \geq 0.$$
Combining the above completes the proof of Theorem \ref{thm:main1}. $\qed$
\endproof

\subsection{Lower bound and proof of Lemma\ \ref{lem:hardinstance}}
This section contains a detailed description and analysis of the hard instance and the formal proof of Lemma \ref{lem:hardinstance}.
\proof{Proof of Lemma\ \ref{lem:hardinstance}:}
Let us fix $k \geq 1$, and consider the stopping problem s.t $T = 2$, $D = 1$, $Z_t = Y_t$ for $t \in \lbrace 1,2 \rbrace$, and $P(Z_1 = 1 - \frac{1}{k+1}) = 1$, $P(Z_2 = 1) = 1 - \frac{1}{k+1}, P(Z_2 = 0) = \frac{1}{k+1}$.  As $\mathbf{Z}$ is a martingale, it follows from optional stopping that $\textsc{opt} = 1 - \frac{1}{k+1}$.  
It then follows from definitions and the basic preservation properties of martingales that $\mathbf{Z}^j$ is a martingale for all $j \geq 1$, and thus $\sup_{\tau \in {\mathcal T}} E\big[Z^j_\tau\big] = E\big[Z^j_1\big]$ for all $j \geq 1$.  A similar argument yields that $Var[Z^j_1] = 0$ for all $j \geq 1$, and that $Z^j_2$ has 2-point support, with one of those points equal to 0.  We conclude that $Z^j_1 = E\big[Z^j_1\big]$, and by the martingale property $Z^j_2$ equals $(k+1) \times Z^j_1$ w.p. $\frac{1}{k+1}$ and equals $0$ w.p. $1 - \frac{1}{k+1}$, for all $j \geq 1$.
\\\indent We now prove by induction that $Z^j_1 = -\frac{1}{k+1}(1-\frac{1}{k+1})^{j-1}$ for all $j \geq 2.$  For the base case $j = 2$, note that 
\begin{eqnarray*}
E[\max(Z^1_1,Z^1_2) | Y_1] &=& \frac{1}{k+1} \times (1 - \frac{1}{k+1}) + (1 - \frac{1}{k+1}) \times 1 
\\&=& (1 - \frac{1}{k+1}) \times (1 + \frac{1}{k+1}).
\end{eqnarray*}
It then follows from the definition of $\mathbf{Z}^2$ that 
\begin{eqnarray*}
Z^2_1 &=& Z^1_1 - E[\max(Z^1_1,Z^1_2) | Y_1]
\\&=& 1 - \frac{1}{k+1} - \big( (1 - \frac{1}{k+1}) \times (1 + \frac{1}{k+1}) \big)
\\&=& - \frac{1}{k+1} \times \big( 1 - \frac{1}{k+1} \big),
\end{eqnarray*}
completing the proof of the base case.  Now, suppose the induction holds for all $j \leq i$ for some $i \geq 2$.  Using the previous observations and the inductive hypothesis, it follows that $Z^i_1 = -\frac{1}{k+1}(1-\frac{1}{k+1})^{i-1}$, and $Z^{i}_2$ equals 
$-(1-\frac{1}{k+1})^{i-1}$ w.p. $\frac{1}{k+1}$, and 0 w.p. $1 - \frac{1}{k+1}$.  It follows that
\begin{eqnarray*}
E[\max(Z^i_1,Z^i_2) | Y_1] &=& \frac{1}{k+1} \times -\frac{1}{k+1}(1-\frac{1}{k+1})^{i-1} + (1 - \frac{1}{k+1}) \times 0
\\&=& -(\frac{1}{k+1})^2(1-\frac{1}{k+1})^{i-1}.
\end{eqnarray*}
Thus
\begin{eqnarray*}
Z^{i+1}_1 &=& Z^i_1 - E[\max(Z^i_1,Z^i_2) | Y_1]
\\&=& -\frac{1}{k+1}(1-\frac{1}{k+1})^{i-1}  + (\frac{1}{k+1})^2(1-\frac{1}{k+1})^{i-1}
\\&=& \frac{1}{k+1} \times (1-\frac{1}{k+1})^{i-1} \times (\frac{1}{k+1} - 1)
\\&=& - \frac{1}{k+1} \times (1-\frac{1}{k+1})^{i},
\end{eqnarray*}
completing the induction.  We conclude that $\sup_{\tau \in {\mathcal T}} E\big[Z^{k+1}_\tau\big] = -\frac{1}{k+1} (1 - \frac{1}{k+1})^k$.  As it is easily verified that $\frac{k}{k+1} (1 - \frac{1}{k+1})^k$ converges to $\frac{1}{e}$ as $k \rightarrow \infty,$ the result then follows from Lemma\ \ref{lem:intuit1}, since $\sup_{\tau \in {\mathcal T}} E\big[Z^{k+1}_\tau\big]$ is the remainder term when the expansion is truncated after $k$ terms.  $\qed$
\endproof

\section{Formal descriptions of algorithms and proofs of algorithmic results}\label{sec:ecalgo}
In this section we provide formal descriptions of our algorithms, and formal proofs of our main algorithmic results.  Several arguments in this section refer back to the proof sketches and statements of results in Section\ \ref{sec: ecproofofmainalgo2}, to ultimately prove our main algorithmic results Theorems\ \ref{thm:mainalgo2}, \ref{thm:mainpolicy1}, \ref{thm:mainmax}, and \ref{thm:mainmaxpolicy}. 

\subsection{Additional notation and auxiliary lemmas}\label{sec:ecnotation}
Before presenting our main algorithms and proofs, it will be helpful to introduce several additional notations, definitions, and preliminary results that will be used repeatedly in later arguments. 
\subsubsection{Additional notations.}
For a D by T matrix $M$, and $t \in [1,T]$, let $M_{[t]}$ denote the submatrix consisting of the first t columns of M (consistent with our essentially identical definition regarding the information process). For $\epsilon, \eta \in (0,1)$, let $N(\epsilon,\eta) \stackrel{\Delta}{=} \lceil \frac{1}{2\epsilon^2} \log(\frac{2}{\eta}) \rceil,$ a quantity well-known to arise in standard concentration arguments, and which will appear frequently in several of our proofs. We also let 
$$
f_k(\epsilon,\delta) \stackrel{\Delta}{=} \log(2\delta^{-1})\times 10^{2 k^2} \times \epsilon^{- 2 k} \times (T+2)^k \times \big(1 + \log(\frac{1}{\epsilon}) + \log(T) \big)^k.
$$
This function will arise in the recursions needed to analyze the complexity of the nested recursions utilized by our algorithms, where (intuitively) $f_k(\epsilon,\delta)$ will coincide with the complexity of approximating $Z^k_t$ to accuracy $\epsilon$ with confidence $1 - \delta$, and matches the associated complexity bound appearing in Lemma\ \ref{almostthere1}.  We note that $f_k$ is the formalization of the rough intuitive complexity measure ${\mathfrak f}_k$ described in Section\ \ref{sec:ecmoreintuit}.
\\\\We state without proof several easily verified facts regarding $N$ and $f_k.$ For $\epsilon,\eta,\delta \in (0,1),$
\begin{align*}
& N(\epsilon,\eta)\ \textrm{and}\ f_k(\epsilon,\delta)\ \textrm{are non-decreasing in}\ \epsilon,\eta,\delta; 
\\& N(\epsilon,\eta) \leq \frac{1}{2\epsilon^2} \log(\frac{2}{\eta})+ 1;
\\& \lceil \log(2 \delta^{-1}) \rceil \leq 2 \log(2 \delta^{-1}).
\end{align*}

\subsubsection{Auxiliary lemmas.}
We recall a standard result from probability theory used often to prove concentration for estimators. 
\begin{lemma}[Hoeffding's inequality]\label{Hofflem}
Suppose that for some $n \geq 1$ and $U > 0$, $\lbrace X_i, i =1,\ldots,n \rbrace$ are i.i.d., and either $P(X_1 \in [ 0 , U]) = 1$ or $P(X_1 \in [-U,0]) = 1$.  Then for all $\eta > 0$,
$ P\bigg(\bigg| n^{-1} \sum_{i=1}^n X_i - E[X_1] \bigg| \geq \eta\bigg) \leq 2 \exp\bigg(-\frac{2 \eta^2 n}{ U^2}\bigg).$
\end{lemma}
\ \indent We next recall the ``median trick'' -- a direct corollary of Hoeffding's inequality that is commonly used in speeding up algorithms for estimating a number, see e.g. \cite{wang2015basics}, and which goes back at least to the work of \cite{nemirovskij1983problem}.  For completeness, and as many variants of the core idea appear throughout the literature, below we state and prove a specific result along these lines.
\begin{lemma}[Median Trick]\label{median}
Suppose there is a given deterministic real number $X \in [0,U]$ (for some $U > 0$), and a randomized algorithm $A$ with the following property: the (random) output $\hat{X}$ of $A$ satisfies $P(|X - \hat{X}| \leq \epsilon) \geq \frac{3}{4}$.  Consider a new algorithm $A'$ constructed as follows: Repeat the algorithm $A$ for $m = \lceil 8\log(2\delta^{-1}) \rceil$ independent trials (with any random numbers used by $A$ independent across trials), yielding outputs $\lbrace \hat{X}_i, i = 1,\ldots,m \rbrace,$ and output the median $X'$ of $\lbrace \hat{X}_i, i = 1,\ldots,m \rbrace.$  \footnote{Here by default we define the median of an ordered sequence $x_1 \leq ... \leq x_n$ (with $n$ even or odd) as $\frac{1}{2}(x_{ \lfloor \frac{n+1}{2} \rfloor } + x_{ \lfloor \frac{n}{2} \rfloor + 1})$.}  Then $P(|X - X'| \leq \epsilon) \geq 1 - \delta$.
\end{lemma}
\proof{Proof of Lemma\ \ref{median}:} Suppose the $m$ outputs from the repetitions of $A$ (when running $A'$) are $q_1,...,q_m$ (possibly non-distinct) and that w.l.o.g. these are non-decreasing (i.e. $q_1 \leq \ldots \leq q_m$).  For $i \in \lbrace 1,\ldots,m \rbrace$, let $Z_i \stackrel{\Delta}{=} I(|q_i - X| \leq \epsilon).$  It may be easily verified from our definition of median that for any $j \in \lbrace 1,\ldots, \lfloor \frac{m + 1}{2}\rfloor \rbrace,$ the interval $[q_j, q_{j + \lfloor \frac{m}{2} \rfloor}]$ must contain the median of the sequence.  Also, the event $\lbrace \sum_{i = 1}^m Z_i > \frac{m}{2} \rbrace$ implies that there exists at least $\lceil \frac{m + 1}{2} \rceil$ elements in $\{q_1,.., q_m\}$ that satisfy $|q_i - x| \leq \epsilon.$  Conditional on this event, let $u$ denote the index of the smallest of these (in the order $q_1 \leq \ldots \leq q_m$) and $v$ the index of the largest.  It may be easily verified that (conditional on that event) $v - u \geq \lfloor \frac{m}{2} \rfloor,$ and $u \leq \lfloor \frac{m + 1}{2}\rfloor$.  Combining with our previous observation about the median, it follows that in this case $\textrm{median}(q_1,...,q_m) \in [q_u, q_v],$ and therefore $|\textrm{median}(q_1,...,q_m) - x| \leq \epsilon.$ Taking the contraposition of the above reasoning, we conclude that  $\lbrace |\textrm{median}(q_1,...,q_m) - x| > \epsilon \rbrace$ implies $\lbrace \sum_{i = 1}^m Z_i \leq \frac{m}{2} \rbrace,$ and thus  
\begin{eqnarray*}
P\big(\big|X' - X\big| > \epsilon\big) &\leq& P\bigg(\frac{1}{m}\sum_{i=1}^m Z_i \leq \frac{1}{2}\bigg) 
\\&=& P\bigg(\frac{1}{m}\sum_{i=1}^m Z_i - E[Z_1] \leq \frac{1}{2} - E[Z_1] \bigg) 
\\&\leq& P\bigg( \frac{1}{m}\sum_{i=1}^m Z_i - E[Z_1] \leq - \frac{1}{4} \bigg) \ \textrm{since}\ E[Z_1] \geq \frac{3}{4}
\\&\leq& P\bigg(\big| \frac{1}{m}\sum_{i=1}^m Z_i - E[Z_1] \big| \geq \frac{1}{4} \bigg).
\end{eqnarray*}
Note that $\lbrace Z_i, i = 1,\ldots,m \rbrace$ viewed as an unordered set has the same joint distribution as $m$ i.i.d. r.v.s each distributed as the indicator of whether a single output of $A$ is within $\epsilon$ of $X$ in absolute value.  Thus we may apply Hoeffding's inequality with $U = 1, n = m = \lceil 8\log(2\delta^{-1}) \rceil,$ and $\eta = \frac{1}{4}$ to conclude that
$$P\big(\big|X' - X\big| > \epsilon\big) \leq 2 \exp\bigg(-\frac{m}{8}\bigg) \leq \delta,$$
thus completing the proof. $\qed$
\endproof
We next prove that the functions $\lbrace f_k, k \geq 1 \rbrace$ satisfy a certain recursive relationship.  This fact will later be used in the proof of Lemma \ref{almostthere1} and Theorem \ref{thm:mainalgo2}.
\begin{lemma}\label{fkrecursebound}
For all $\epsilon, \delta \in (0,1)$ and $k \geq 1,$
$$
f_{k+1}(\epsilon,\delta) \geq 
16\lceil \log(2\delta^{-1}) \rceil\times \big( N(\frac{\epsilon}{4},\frac{1}{16}) + 1 \big) \times (T+2) \times f_k\big( \frac{\epsilon}{4} ,\frac{1}{16 N(\frac{\epsilon}{4},\frac{1}{16}) T} \big).
$$
\end{lemma}
\proof{Proof of Lemma\ \ref{fkrecursebound}: } Combining definitions, the monotonicty of both $f_k(\epsilon,\delta)$ and $N(\epsilon,\eta)$ in $\epsilon,\delta,\eta$ on $(0,1),$ and some straightforward algebra (including the facts that $\lceil x \rceil \leq 2 x$ for $x \geq \log(2), x + 1 \leq 2x$ for $x \geq 1$, and $8 x + 2 \leq 10 x $ for $x \geq 1$), we have that
\begin{eqnarray*}
& & 16\lceil \log(2\delta^{-1}) \rceil\times \big( N(\frac{\epsilon}{4},\frac{1}{16}) + 1 \big) \times (T+2) \times f_k\big( \frac{\epsilon}{4} ,\frac{1}{16 N(\frac{\epsilon}{4},\frac{1}{16}) T} \big) 
\\& & \leq  16\lceil \log(2\delta^{-1}) \rceil\times \big( 8 \log(32) \epsilon^{-2} + 2 \big) \times (T+2) \times f_k\big( \frac{\epsilon}{4} ,\frac{1}{16 N(\frac{\epsilon}{4},\frac{1}{16}) T} \big) 
\\& & \leq  \big( 32 \log(2 \delta^{-1}) \big) \times \big( 10 \log(32) \epsilon^{-2} \big) \times (T + 2) \times f_k\big( \frac{\epsilon}{4}, \frac{\epsilon^2}{256 \log(32) T} \big)
\\& & =  \big( 32 \log(2 \delta^{-1}) \big) \times \big( 10 \log(32) \epsilon^{-2} \big) \times (T + 2) 
\\& &\ \ \ \ \ \ \ \ \ \times \log\big( 512 \log(32) T \epsilon^{-2} \big) \times 10^{2 k^2} \times 4^{2 k} \times \epsilon^{- 2 k} \times (T+2)^k \times \big(1 + \log(\frac{4}{\epsilon}) + \log(T) \big)^k
\\ & & \leq \ \ \log(2\delta^{-1}) \times (T+2)^{k+ 1} \times \epsilon^{-2 (k + 1)} \times \bigg(1 + \log\big(\frac{1}{\epsilon}\big) + \log(T)\bigg)^{k + 1} \times 10^{2 k^2} \times 4^{2 k} \times 10^4\times 4^k
\\& &\ \ \ \ \ \textrm{since}\ \log\big( 512 \log(32) T \epsilon^{-2} \big) \leq  \log\big( 512 \log(32) \big) \big(1 + \log\big(\frac{1}{\epsilon}\big) + \log(T)\big), 
\\&\ &\ \ \ \ \ \ \ \ \ \ \ \ \big(1 + \log(\frac{4}{\epsilon}) + \log(T) \big) \leq 4 \big(1 + \log(\frac{1}{\epsilon}) + \log(T) \big),\ \textrm{and}
\\&\ &\ \ \ \ \ \ \ \ \ \ \ \ \ \ \ \ \ \ \ \ \ \ \ \  32 \times 10 \log(32) \times  \log\big( 512 \log(32) \big) \leq 10^4
\\ & & \leq \ \ \log(2\delta^{-1}) \times 10^{2 (k + 1)^2} \times (T+2)^{k + 1} \times \epsilon^{-2 (k + 1)} \times \bigg(1 + \log\big(\frac{1}{\epsilon}\big) + \log(T)\bigg)^{k + 1},
\end{eqnarray*}
the final inequality following from the fact that
\begin{eqnarray*}
\ &\ &\ 10^{2 k^2} \times 4^{2 k} \times 10^4\times 4^k
\\&=&  10^{2 k^2 + 4} \times 4^{3 k} \ \leq \ 10^{2 k^2 + 4} \times 10^{2 k}  = 10^{2(k + 1)^2 - 2 k + 2}
\\&\leq& 10^{2 (k + 1)^2 } \ \ \ \ \ \ \textrm{for $k \geq 1.$}
\end{eqnarray*}
\endproof

\indent To conclude this set of preliminary results, let us state several well-known results regarding the c.d.f. and hazard rate of the standard normal distribution, which follow from e.g. the results of \cite{gordon1941values}, \cite{baricz2008mills}, and \cite{gasull2014approximating}.  Let $\Phi(\cdot)$ denote the standard normal c.d.f., and $\Phi^c(\cdot) \stackrel{\Delta}{=} 1 - \Phi(\cdot)$.
\begin{lemma}\label{lmnormaltail}
For all  $x > 0,$ one has
\[
\frac{1}{\sqrt{2 \pi}}\bigg(\frac{1}{x} - \frac{1}{x^3}\bigg) e^{-\frac{x^2}{2}} \leq \Phi^c(x) \leq \frac{1}{\sqrt{2 \pi}}\frac{1}{x} e^{-\frac{x^2}{2}}.
\]
Also, $\Phi^c(x) \leq \sqrt{\frac{2}{\pi}} e^{-\frac{x^2}{2}} (x + \sqrt{x^2 + \frac{8}{\pi}})^{-1}$.  Furthermore, let $R(x, a) \triangleq \frac{\Phi^c(x + a)}{\Phi^c(x)}.$ Then $R(x,a)$ is monotone decreasing in $a$ for each fixed $x \in \mathbb{R}$, and monotone decreasing in $x$ for each fixed $a \geq 0$. 
\end{lemma}

\subsection{Description of algorithms $\mathcal{B}^k$}\label{ecBk}
We recursively define algorithms $\{{\mathcal B}^k\}_{k \geq 1}$, which takes inputs $(t,\gamma,\epsilon,\delta)$ and returns an (additive) $\epsilon$-approximation to $Z^k_t(\gamma)$ w.p. at least $1 - \delta$.  Here we recall that ${\mathcal B}$ is the simulation subroutine defined in Section\ \ref{sec:algoresults}.
\\\\\hrule
\vspace{2pt}
 Algorithm ${\mathcal B}^1(t,\gamma,\epsilon,\delta)$:
 \vspace{2pt}
 \hrule
 \vspace{2pt}
\indent\indent\indent\indent Return $g_t(\gamma)$
\vspace{2pt}
\hrule
\vspace{8pt}
\ \\ \hrule
\vspace{2pt}
Algorithm ${\mathcal B}^{k+1}(t,\gamma,\epsilon,\delta)$:
\vspace{2pt}
 \hrule
\ \\\indent\indent Create a length-$\lceil 8 \log(2\delta^{-1}) \rceil$ vector $\mathbf{A}^m$
\\\indent\indent For s = 1 to $\lceil 8 \log(2\delta^{-1}) \rceil$
\\\indent\indent\indent\indent Generate an ind. call to ${\mathcal B}^k(t,\gamma,\frac{\epsilon}{2}, \frac{1}{8})$ and store as variable $A^3$
\\\indent\indent\indent\indent\indent\indent\indent\indent\indent\indent\indent\indent\indent $\backslash\backslash$\textit{ compute an\ $\epsilon/2-$approx. of $Z^k_t(\gamma)$ w.p. $\frac{7}{8}$}
\\\indent\indent\indent\indent Create a length-$N(\frac{\epsilon}{4},\frac{1}{16})$ vector $\mathbf{A}^0$
\\\indent\indent\indent\indent For i = 1 to $N(\frac{\epsilon}{4},\frac{1}{16})$
\\\indent\indent\indent\indent\indent\indent Generate an ind. call to ${\mathcal B}(t,\gamma)$ and store in D by T matrix $\mathbf{A}^1$
\\\indent\indent\indent\indent\indent\indent\indent\indent\indent\indent\indent\indent\indent $\backslash\backslash$\textit{ draw an ind. sample path of $\mathbf{Y}$ cond. on $Y_t = \gamma$} 
\\\indent\indent\indent\indent\indent\indent Create a length-T vector $\mathbf{A}^2$
\\\indent\indent\indent\indent\indent\indent For j = 1 to T
\\\indent\indent\indent\indent\indent\indent\indent\indent Generate an ind. call to ${\mathcal B}^k\big(j,\mathbf{A}^1_{[j]},\frac{\epsilon}{4} ,\frac{1}{16 N(\frac{\epsilon}{4},\frac{1}{16}) T} \big)$ and store in $A^2_j$
\\\indent\indent\indent\indent\indent\indent Compute the maximum value of $\mathbf{A}^2$ and store in $A^0_i$
\\\indent\indent\indent\indent\indent\indent\indent\indent\indent\indent\indent\indent\indent $\backslash\backslash$\textit{ approx. $\max_{i = 1,\ldots,T]}Z^k_i$ on the given sample path} 
\\\indent\indent\indent\indent Compute $A^3 - \big(N(\frac{\epsilon}{4},\frac{1}{16})\big)^{-1} \sum_{i=1}^{N(\frac{\epsilon}{4},\frac{1}{16})} A^0_i$ and store in $A^m_s$
\\\indent\indent\indent\indent\indent\indent\indent\indent\indent\indent\indent\indent\indent $\backslash\backslash$\textit{ compute an $\epsilon-$approx. of $Z^{k+1}_t(\gamma)$ w.p. $\frac{3}{4}$}
\\\indent\indent Return the median of $\mathbf{A}^m$
\vspace{4pt}
\hrule
\vspace{12pt}
\subsection{Proof of Lemma \ref{almostthere1}}\label{eclemma3}
We now show the algorithms introduced in the last section satisfy the desired approximation and complexity guarantee, thereby proving Lemma\ \ref{almostthere1}. 
\proof{Proof of Lemma \ \ref{almostthere1}: } We proceed by induction. In the base case $k = 1,$ $\mathcal{B}^1$ simply makes one call to evaluate $g_t.$ The output is exact with $G$ units of computational time. Combining with the definition of $f_k$, we have that the base case trivially holds true.
\\\indent  Now suppose the induction is true for some $k \geq 1.$ We first show that ${\mathcal B}^{k+1}$ satisfies the desired high probability error bounds. By Lemma\ \ref{median}, it suffices to show that in each of the $\lceil 8 \log(2\delta^{-1}) \rceil$ outer loops, the algorithm ${\mathcal B}^{k+1}$ computes and stores a (random) number $R$ s.t. $P\big(|R - Z^{k+1}_t(\gamma)| > \epsilon\big) \leq \frac{1}{4},$ as each iteration of this loop uses independent randomness (with the common $\gamma$).
Let $\lbrace X_i, i = 1,\ldots,N(\frac{\epsilon}{4},\frac{1}{16}) \rbrace$ be a sequence of r.v.s constructed as follows.  For the fixed $\gamma,$ let $\big\lbrace \mathbf{Y}^i(\gamma), i = 1,\ldots,N(\frac{\epsilon}{4},\frac{1}{16}) \big\rbrace$ be a sequence of (conditionally) independent random matrices, each distributed as $\mathbf{Y}(\gamma)$.  For $i = 1,\ldots, N(\frac{\epsilon}{4},\frac{1}{16}),$ let $X_i = \max_{j = 1,\ldots,T} Z^k_j(\mathbf{Y}^i(\gamma)_{[j]}).$
Then it follows from our inductive hypothesis, the Lipschitz property of the max function, a union bound over all $i = 1,\ldots,N(\frac{\epsilon}{4},\frac{1}{16})$ and $j \in \lbrace 1,\ldots,T \rbrace$, and some straightforward algebra, that we may construct $\lbrace X_i, i = 1,\ldots, N(\frac{\epsilon}{4},\frac{1}{16}) \rbrace$ and $\lbrace A^0_i, i = 1,\ldots, N(\frac{\epsilon}{4},\frac{1}{16}) \rbrace$ on a common probability space s.t. with probability at least $1 - \frac{1}{16}$, $|X_i - A^0_i| < \frac{\epsilon}{4}$ for all $i = 1,\ldots,N(\frac{\epsilon}{4},\frac{1}{16})$.  Applying Lemma\ \ref{Hofflem} to $\lbrace X_i, i = 1,\ldots,N(\frac{\epsilon}{4},\frac{1}{16}) \rbrace$, with parameters $\eta = \frac{\epsilon}{4}, U = 1, n = N(\frac{\epsilon}{4},\frac{1}{16})$, we conclude (after some straightforward algebra) that on the same probability space, 
$$P\bigg( \big| \big(N(\frac{\epsilon}{4},\frac{1}{16})\big)^{-1} \sum_{i=1}^{N(\frac{\epsilon}{4},\frac{1}{16})} X_i - E[X_1] \big| < \frac{\epsilon}{4} \bigg) \geq 1 - \frac{1}{16}.$$  
Here we can apply Lemma\ \ref{Hofflem} with $U = 1$ since $X_i \in [-1,0]$ for all $i \geq 1.$  Furthermore, the event $\bigg\lbrace |X_i - A^0_i| < \frac{\epsilon}{4}\ \textrm{for all
}\ i = 1,\ldots,N(\frac{\epsilon}{4},\frac{1}{16}) \bigg\rbrace$, which occurs with probability at least $1 - \frac{1}{16},$ implies the event
$$\bigg\lbrace \bigg| \big(N(\frac{\epsilon}{4},\frac{1}{16})\big)^{-1} \sum_{i=1}^{N(\frac{\epsilon}{4},\frac{1}{16})} A^0_i - \big(N(\frac{\epsilon}{4},\frac{1}{16})\big)^{-1} \sum_{i=1}^{N(\frac{\epsilon}{4},\frac{1}{16})} X_i \bigg| < \frac{\epsilon}{4} \bigg\rbrace.$$
By combining the above with a union bound and the triangle inequality, we conclude that on the same probability space (and hence in general),
$$P\bigg( \bigg| \big(N(\frac{\epsilon}{4},\frac{1}{16})\big)^{-1} \sum_{i=1}^{N(\frac{\epsilon}{4},\frac{1}{16})} A^0_i - E[X_1] \bigg| < \frac{\epsilon}{2} \bigg) \geq 1 - \frac{1}{8}.$$
As the inductive hypothesis ensures that $P\big( \big| A^3 - Z^k_t(\gamma) \big| > \frac{\epsilon}{2} \big) \leq \frac{1}{8}$, we may again apply a union bound and the triangle inequality, along with the definition of $Z^k_t(\gamma)$ and $X_1$, to conclude that 
\begin{equation}\label{step1}
P\bigg( \bigg| A^3 - \big(N(\frac{\epsilon}{4},\frac{1}{16})\big)^{-1} \sum_{i=1}^{N(\frac{\epsilon}{4},\frac{1}{16})} A^0_i - Z^{k+1}_t(\gamma) \bigg| > \epsilon \bigg) \leq \frac{1}{4}\ \ \ \textrm{as desired.}
\end{equation}
\\\indent We next focus on computational costs.  For each fixed $s \in \lbrace 1,\ldots,\lceil 8 \log(2\delta^{-1}) \rceil \rbrace,$
in each of the $N(\frac{\epsilon}{4},\frac{1}{16})$ iterations (indexed by i), first one direct call is made to ${\mathcal B}(t,\gamma)$ (at computational cost C); then T calls are made to ${\mathcal B}^k\big(j,A^1_{[j]},\frac{\epsilon}{4} ,\frac{1}{16 N(\frac{\epsilon}{4},\frac{1}{16}) T} \big)$ (for different values of j), each at computational cost at most $(C + G+ 1) \times f_k\big( \frac{\epsilon}{4} ,\frac{1}{16 N(\frac{\epsilon}{4},\frac{1}{16}) T} \big)$; then the maximum of a length-T array is computed (at computational cost T).  One additional call is then made to ${\mathcal B}^k(t,\gamma,\frac{\epsilon}{2}, \frac{1}{8})$, at computational cost at most 
$(C + G + 1) \times f_k(\frac{\epsilon}{2}, \frac{1}{8});$ the average of $N(\frac{\epsilon}{4},\frac{1}{16})$ is computed and subtracted from $A^3$, at computational cost $N(\frac{\epsilon}{4},\frac{1}{16}) + 1;$ and finally the median of $\lceil8 \log(2\delta^{-1})\rceil$ numbers is computed. at a cost of $\lceil8 \log(2\delta^{-1})\rceil$.  It thus follows from the inductive hypothesis and some straightforward algebra that the computational cost of ${\mathcal B}^{k+1}(t,\gamma,\epsilon,\delta)$ is at most 
\begin{eqnarray}
\ &\ &\ 8\times \lceil \log(2\delta^{-1}) \rceil\times \bigg( N(\frac{\epsilon}{4},\frac{1}{16}) C + N(\frac{\epsilon}{4},\frac{1}{16}) T (C + G+ 1) f_k\big( \frac{\epsilon}{4} ,\frac{1}{16 N(\frac{\epsilon}{4},\frac{1}{16}) T} \big) \nonumber 
\\&\ &\ \ \ \ +\ \ N(\frac{\epsilon}{4},\frac{1}{16}) T + (C + G+ 1) f_k(\frac{\epsilon}{2}, \frac{1}{8}) + N(\frac{\epsilon}{4},\frac{1}{16}) + 1 + 1 \bigg) \nonumber
\\&\ &\ \ \ \ \ \ \ \ \leq\ \ \ 16\times \lceil \log(2\delta^{-1}) \rceil\times(C + G+ 1) \times \big( N(\frac{\epsilon}{4},\frac{1}{16}) + 1 \big) \times (T+2) \times f_k\big( \frac{\epsilon}{4} ,\frac{1}{16 N(\frac{\epsilon}{4},\frac{1}{16}) T} \big). \nonumber
\\&\ &\ \ \ \ \ \ \ \ \leq\ \ \ (C + G + 1) \times f_{k+1}(\epsilon, \delta). \nonumber
\end{eqnarray}
The last step follows from Lemma\ \ref{fkrecursebound}. By the definition of $f_k$, we then get the desired sampling and computational complexity bound for ${\mathcal B}^{k+1}$. Combining all of the above, we complete the proof by induction. 
\qed
\endproof

\subsection{Proof of Theorem \ref{thm:mainalgo2}}\label{ecthm3}
In this section, we complete our proof of Theorem \ref{thm:mainalgo2}. First, let us formally define algorithm $\mathcal{A}$, which uses ${\mathcal B}^k$ and simulation to approximate \textsc{opt}.
\ \\\\\hrule
\vspace{2pt}
Algorithm $\mathcal{A}(\epsilon,\delta)$:
\vspace{2pt}
\hrule
\vspace{4 pt}
\indent\indent Set $k = \lceil\frac{2}{\epsilon}\rceil, \alpha = \frac{\epsilon}{2} \big( \lceil \frac{2}{\epsilon} \rceil \big)^{-1}, \beta = \delta \big( \lceil \frac{2}{\epsilon} \rceil \big)^{-1}$, and create a length-$k$ vector $\hat{\mathbf{L}}$
\\\indent\indent For l = 1 to k
\\\indent\indent\indent\indent Create a 
length-$\lceil8\log(2 \beta^{-1})\rceil$ vector $\mathbf{A}^l$
\\\indent\indent\indent\indent For $s = 1$ to $\lceil8\log(2 \beta^{-1})\rceil$
\\\indent\indent\indent\indent\indent\indent Create a length-$N(\frac{\alpha}{2},\frac{1}{8})$ vector $\mathbf{A}^0$
\\\indent\indent\indent\indent\indent\indent For i = 1 to $N(\frac{\alpha}{2},\frac{1}{8})$
\\\indent\indent\indent\indent\indent\indent\indent\indent Generate an ind. call to ${\mathcal B}(0,\emptyset)$ and store in D by T matrix $\mathbf{A}^1$
\\\indent\indent\indent\indent\indent\indent\indent\indent Create a length-T array $\mathbf{A}^2$
\\\indent\indent\indent\indent\indent\indent\indent\indent For j = 1 to T
\\\indent\indent\indent\indent\indent\indent\indent\indent\indent\indent Generate an ind. call to ${\mathcal B}^l\big(j,\mathbf A^1_{[j]},\frac{\alpha}{2} ,\frac{1}{8 N(\frac{\alpha}{2},\frac{1}{8}) T} \big)$ and store in $ A^2_j$
\\\indent\indent\indent\indent\indent\indent\indent\indent Compute the maximum value of $\mathbf{A}^2$ and store in $ A^0_i$
\\\indent\indent\indent\indent\indent\indent Compute $\big(N(\frac{\alpha}{2},\frac{1}{8})\big)^{-1} \sum_{i=1}^{N(\frac{\alpha}{2},\frac{1}{8})}  A^0_i$ and store in $ A^l_s$ 
\\\indent\indent\indent\indent Take the median of $\mathbf{A}^l$ and store in $\hat L_l$
\\\indent\indent Return $\sum_{l=1}^k  \hat L_l$  
\vspace{4 pt}
\hrule
\vspace{8pt}
\ \proof{Proof of Theorem\ \ref{thm:mainalgo2}: }
We first show that $\mathcal{A}(\epsilon,\delta)$ satisfies the desired high probability error bounds.  In light of  Lemma\ \ref{almostthere1} and Theorem\ \ref{thm:main2}, by a union bound and the triangle inequality, to output an $\epsilon-$approximation of $\textsc{opt}$ with probability at least $1 - \delta,$ it suffices to individually approximate the first $\lceil\frac{2}{\epsilon}\rceil$ terms $(L_l)_{l = 1,\ldots,\lceil\frac{2}{\epsilon}\rceil}$ of the expansion each to within additive error $\alpha$ with probability $1- \beta.$  Here we recall that $\alpha = \frac{\epsilon}{2} \big( \lceil \frac{2}{\epsilon} \rceil \big)^{-1}, \beta = \delta \big( \lceil \frac{2}{\epsilon} \rceil \big)^{-1}.$  By Lemma \ref{median}, to get an $\alpha-$approximation of each $L_l$ with probability at least $1 - \beta,$ it suffices that for each $s \in \lbrace 1,\ldots,\lceil 8\log(2 \beta^{-1})\rceil \rbrace,$ the algorithm computes and stores a number $A^l_s$ that satisfies $P(|A^l_s - L_l| > \alpha) \leq \frac{1}{4}.$ The proof of this fact is nearly identical to that of Lemma \ \ref{almostthere1} (with several union bounds and an application of Lemma \ref{Hofflem}), and we omit the details.
\\\indent We now prove the computational complexity bounds.  By Lemma\ \ref{almostthere1} and a runtime analysis nearly identical to that used in the proof of Lemma\ \ref{almostthere1}, we find that for each one iteration of the outer loop $l \in \lbrace 1,\ldots,k \rbrace,$ the computational complexity is bounded by
\begin{equation}\label{interabound}
(C + G+ 1)\times \lceil 8 \log(2 \beta^{-1}) \rceil \times \Bigg( N(\frac{\alpha}{2}, \frac{1}{8}) \times  \bigg( T \times f_l(\frac{\alpha}{2}, \frac{1}{8 N(\frac{\alpha}{2},\frac{1}{8}) T}) + T + 1\bigg) + 1 \Bigg).
\end{equation}
It then follows from Lemma\ \ref{fkrecursebound}, the monotonicity of $f_l$ and related auxiliary functions, and some straightforward algebra that (\ref{interabound}) is at most $f_{l + 1}\big(\alpha, \beta\big).$  Accounting for the additional $k$ units of runtime needed to compute $\sum_{l=1}^k \hat L_l$, the total computational complexity (divided by $C+ G + 1$) is at most
\begin{eqnarray*}
& & \sum_{l = 1}^{\lceil\frac{2}{\epsilon}\rceil}  f_{l + 1}\big(\alpha, \beta\big) + \lceil\frac{2}{\epsilon}\rceil 
\\ & & \ \ \leq  \ \ \frac{3}{\epsilon} f_{\lceil\frac{2}{\epsilon}\rceil + 1}(\frac{\epsilon}{2}(\lceil\frac{2}{\epsilon}\rceil)^{-1}, \delta(\lceil\frac{2}{\epsilon}\rceil)^{-1}) + \frac{3}{\epsilon}\ \textrm{by monotonicity of}\ f_k
\\ & & \ \ \leq \ \ 6 \epsilon^{-1} f_{\lceil\frac{2}{\epsilon}\rceil + 1}\big(\frac{\epsilon^2}{6}, \frac{\delta\epsilon}{3}\big) \ \leq \ 6 \epsilon^{-1} f_{\frac{4}{\epsilon}}\big(\frac{\epsilon^2}{6}, \frac{\delta\epsilon}{3}\big)\ \textrm{since}\ \lceil \frac{2}{\epsilon} \rceil \leq \frac{3}{\epsilon}
\\ & & \ \ = \ \ 6 \epsilon^{-1} \times  \log\big(6\delta^{-1}\epsilon^{-1} \big)\times 10^{2(4\epsilon^{-1})^2} \times (6 \epsilon^{-2})^{2(4\epsilon^{-1})}\times (T+2)^{4\epsilon^{-1}} \times \big(1 + \log(6\epsilon^{-2}) + \log(T)\big)^{4\epsilon^{-1}}
\\ & & \ \ = \ \ 6\epsilon^{-1} \times (\log(6) + \log(\epsilon^{-1}) + \log(\delta^{-1}))
\\& & \ \indent\indent\indent\indent\indent\indent\indent\times\ \ 10^{32\epsilon^{-2}}  \times 6^{8\epsilon^{-1}} \times \epsilon^{-16\epsilon^{-1}} \times (T+2)^{4\epsilon^{-1}} \times \big(1 + \log(6) + 2\log(\epsilon^{-1}) + \log(T)\big)^{4\epsilon^{-1}}
\\ & & \ \ \leq \ \ \big(\log(\delta^{-1}) + 1\big) (1 + \log(6)) \big( \log(\epsilon^{-1}) + 1 \big) 10^{32\epsilon^{-2}} 6^{8\epsilon^{-1}+1}\epsilon^{-16\epsilon^{-1}- 1}(T+2)^{4\epsilon^{-1}} 
\\& & \ \indent\indent\indent\indent\indent\indent\indent\times\ \ \big(1 + \log(6)\big)^{4\epsilon^{-1}} (1 + \log(T))^{4\epsilon^{-1}} (1 +  2 \log(\epsilon^{-1}))^{4\epsilon^{-1}}
\\& & \ \indent\indent\indent\indent\indent\indent\indent\indent\indent\indent\indent\indent\indent\indent \ \ \textrm{since}\ 1 + x + y + z \leq (1 + x)(1 + y)(1 + z) \ \textrm{for}\ x,y,z \geq 0
\\ & & \ \ \leq \ \ \big(1 + \log(\delta^{-1})\big) \times 10^{32\epsilon^{-2}}\times 6^{14\epsilon^{-1}}\times \epsilon^{-26\epsilon^{-1}} \times (3 T)^{8\epsilon^{-1}}
\\& & \ \indent\indent\indent\indent\indent\indent\indent\indent\indent\indent\indent\indent\indent\indent \ \ \textrm{since}\ 1 +  a \log(b) \leq b^a 
\\ & & \ \ \leq \ \ (1 + \log(\delta^{-1}) ) \times \exp(100\epsilon^{-2})\times T^{8 \epsilon^{-1}}.
\end{eqnarray*}
Combining the above completes the proof.  $\qed$  
\endproof

\subsection{Proof of Theorem \ref{thm:mainpolicy1}}\label{goodpolicydetails}
In this section, we complete the proof of Theorem\ \ref{thm:mainpolicy1}.  First, we make the description of our efficiently computable randomized stopping time $\tau_{\epsilon}$ completely explicit.  At time 1, after seeing $Y_{[1]}$, make an independent call to ${\mathcal B}^{\lceil 4 \epsilon^{-1} \rceil + 1}\bigg( 1, Y_{[1]}, \frac{\epsilon}{4}, \frac{\epsilon}{4 T} \bigg)$.  If the value returned is at 
least $-\frac{1}{2} \epsilon$, stop.  If not, continue.  Suppose that for some $t \in \lbrace 1,\ldots,T-2 \rbrace$, we have not yet stopped by the end of period $t$.  At time $t + 1$, after observing $Y_{t+1}$, make an independent call to 
${\mathcal B}^{\lceil 4 \epsilon^{-1} \rceil + 1 }\bigg( t + 1, Y_{[t+1]}, \frac{\epsilon}{4}, \frac{\epsilon}{4 T} \bigg)$.  If the value returned is at least $-\frac{1}{2} \epsilon$, stop.  If not,  continue.  Finally, if we have not yet stopped by period $T$, stop in period $T$.  It is easily verified that for any $\epsilon \in (0,1)$, the implied stopping time $\tau_{\epsilon}$ is a well-defined, appropriately adapted, randomized stopping time.  We now use $\tau_{\epsilon}$ to complete the proof of Theorem\ \ref{thm:mainpolicy1}.
\proof{Proof of Theorem\ \ref{thm:mainpolicy1}: }
Let $k_{\epsilon} \stackrel{\Delta}{=} \lceil \frac{4}{\epsilon} \rceil + 1$.  Let ${\mathcal G}_{1,\epsilon}$ denote the event
$$
\Bigg\lbrace \bigg| {\mathcal B}^{k_{\epsilon}}\bigg(t, Y_{[t]}, \frac{\epsilon}{4}, \frac{\epsilon}{4 T} \bigg) - Z^{k_{\epsilon}}_t(Y_{[t]}) \bigg| \leq \frac{\epsilon}{4}\ \ \ \forall\ \ t \in \lbrace 1,\ldots,T \rbrace \Bigg\rbrace,$$
${\mathcal G}_{2,\epsilon}$ denote the event
$$
\Bigg\lbrace \exists\ t \in \lbrace 1,\ldots,T \rbrace\ \textrm{such that}\ - {\mathcal B}^{k_{\epsilon}}\bigg(t, Y_{[t]}, \frac{\epsilon}{4}, \frac{\epsilon}{4 T} \bigg) \leq \frac{1}{2} \epsilon \Bigg\rbrace,$$
and ${\mathcal G}_{3,\epsilon}$ denote the event $\big\lbrace Z^{k_{\epsilon}}_{\tau_{\epsilon}} \geq - \frac{3}{4} \epsilon \big\rbrace$.  Observe that Lemma\ \ref{converge1}, definitions, and several straightforward union bounds and applications of the triangle inequality ensure that: (1) $P({\mathcal G}_{1,\epsilon}) \geq 1 - \frac{\epsilon}{4};$ (2) $P({\mathcal G}_{2,\epsilon} | {\mathcal G}_{1,\epsilon}) = 1;$ (3) $P({\mathcal G}_{3,\epsilon} | {\mathcal G}_{1,\epsilon} \bigcap {\mathcal G}_{2,\epsilon}) = 1$.  It follows that $P\big({\mathcal G}^c_{3,\epsilon}\big) \leq \frac{\epsilon}{4}$, and thus since by our assumptions and monotonicity $P(Z^{k_{\epsilon}}_t \geq - 1) = 1$ for all $t \in \lbrace 1,\ldots,T \rbrace$,
\begin{eqnarray*}
E\big[ Z^{k_{\epsilon}}_{\tau_{\epsilon}} \big] &=& E\big[ Z^{k_{\epsilon}}_{\tau_{\epsilon}} I({\mathcal G}_{3,\epsilon}) \big] + E\big[ Z^{k_{\epsilon}}_{\tau_{\epsilon}} I({\mathcal G}^c_{3,\epsilon}) \big] 
\\&\geq& -\frac{3}{4} \epsilon - E\big[I({\mathcal G}^c_{3,\epsilon}) \big]\ \ \ \geq\ \ \ - \epsilon.
\end{eqnarray*}
Combining with Lemma\ \ref{good1} and Theorem\ \ref{thm:main1}, we complete the proof of the first part of the theorem.  As for the computational cost, the algorithm in each time step makes one call to $\mathcal{B}^{\lceil4\epsilon^{-1}\rceil + 1}$ to compute an $\frac{\epsilon}{4}-$approximation with probability at least $1 -  \frac{\epsilon}{4 T}$.  The desired result then follow from Lemma\ \ref{almostthere1} along with some straightforward algebra, and we omit the details.  $\qed$
\endproof

\subsection{Unbounded rewards and proofs of Theorems \ref{thm:mainmax} and \ref{thm:mainmaxpolicy}}\label{sec:ecmaximization}

In this section we present the proofs of Theorems\ \ref{thm:mainmax} and \ref{thm:mainmaxpolicy}, and all related auxiliary results.  Our proofs proceed by: (1) bounding the error induced by truncating the rewards; (2) bounding $\textsc{opt}$ in terms of the mean and s.c.v. of the maximum reward (which will allow us to relate the truncation error to $\textsc{opt}$); and (3) applying our previous algorithms and results for the bounded case to the truncated problem.  In addition, we prove an alternative bound not depending on the s.c.v. in Section\ \ref{altubsec}.  

\subsubsection{Additional Notation.}
For any $U > 0,$ define $Z^1_{U,t} \stackrel{\Delta}{=} U^{-1}\times \min\big(U, Z^1_t\big)$.  For $k \geq 1,$ let $Z^{k+1}_{U,t} \stackrel{\Delta}{=} Z^{k}_{U,t} - E\big[\max_{i = 1,\ldots,T}Z^{k}_{U,i} |\mathcal{F}_t\big]; L_{U,k} \stackrel{\Delta}{=} E\big[\max_{i = 1,\ldots,T}Z^{k}_{U,i}\big];$ and $E_{U,k} \stackrel{\Delta}{=} \sum_{i =1}^k L_{U,i}.$ 
Let $M_1 = E\big[\max_{t = 1,\ldots,T}Z^1_t\big]$, $M_2 = E\big[(\max_{t = 1,\ldots,T}Z^1_t)^2\big],$ and $\gamma_0 = \frac{M_2}{(M_1)^2}$, i.e. the constant specified in Assumption \ref{assumptionA2}. 
\subsubsection{Bounding the truncation error.}\label{proofofmainmax}
In this section, we bound the error induced (in the optimal stopping value) by truncating the rewards.  
\begin{lemma}\label{maxb1}
For all $U > 0$, $$0 \leq \textsc{opt} - U \times \sup_{\tau \in {\mathcal T}} E[Z^1_{U,\tau}] \leq (M_2)^{\frac{1}{2}} \times \big( \frac{M_1}{U} \big)^{\frac{1}{2}}.$$
\end{lemma}
\proof{Proof}
Non-negativity follows from the fact that w.p.1 $Z^1_t \geq U \times Z^1_{U,t}$ for all $t \in \lbrace 1,\ldots,T \rbrace$.  To prove the other direction, let $\tau_*$ denote an optimal stopping time for the problem 
$\sup_{\tau \in {\mathcal T}} E[Z^1_{\tau}]$, where existence follows from \cite{chowgreat}.  Then by a straightforward coupling and rescaling,
\begin{eqnarray*}
E[Z^1_{\tau_*}] - U \times E[Z^1_{U,\tau_*}] &\leq& E\big[Z^1_{\tau_*} I\big( Z^1_{\tau_*} > U\big)\big]
\\&\leq& E\big[Z^1_{\tau_*} I\big( \max_{t = 1,\ldots,T} Z^1_t > U\big)\big]
\\&\leq& E\big[ \big( \max_{t = 1,\ldots,T} Z^1_t \big) \times  I\big( \max_{t = 1,\ldots,T} Z^1_t > U\big)\big]
\\&\leq& (M_2)^{\frac{1}{2}} \times \bigg( P\big( \max_{t = 1,\ldots,T} Z^1_t > U \big) \bigg)^{\frac{1}{2}}\ \ \ \textrm{by Cauchy-Schwarz}
\\&\leq& (M_2)^{\frac{1}{2}} \times \big( \frac{M_1}{U} \big)^{\frac{1}{2}}\ \ \ \textrm{by Markov's inequality.}
\end{eqnarray*}
Noting that $\sup_{\tau \in {\mathcal T}} E[Z^1_{U,\tau}]  \geq E[Z^1_{U,\tau_*}]$ completes the proof.  $\qed$
\endproof

\subsubsection{Bounding $\textsc{opt}$ in terms of $M_1$ and $\gamma_0$.}
In this section, we show that $\textsc{opt}$ may be bounded (from below) as a simple function of $M_1$ and $\gamma_0$. 
\begin{lemma}\label{maxb2}
$\textsc{opt} \geq \frac{4}{27} \times \gamma_0^{-1} \times M_1.$
\end{lemma}
\proof{Proof}
Recall the celebrated Paley-Zygmund inequality, i.e. the fact that for any $\delta \in (0,1)$ and non-negative r.v. $X$,
\begin{equation}\label{paley1}
P\big( X > \delta E[X] \big) \geq (1 - \delta)^2 \times \frac{ (E[X])^2 }{E[X^2]}.
\end{equation}
Now, for $\delta \in (0,1)$, consider the stopping time $\tau_{\delta}$ which stops the first time that that $Z^1_t \geq \delta \times M_1$, and stops at time $T$ if no such time exists in $\lbrace 1,\ldots,T \rbrace$.  Then by non-negativity and (\ref{paley1}),
$$E[Z^1_{\tau_{\delta}}]\ \ \ \geq\ \ \ \delta \times (1 - \delta)^2 \times \frac{(M_1)^3}{ M_2 }.$$
Optimizing over $\delta$ (a straightforward exercise in calculus) then completes the proof.  $\qed$
\endproof

\subsubsection{Proof of Theorem\ \ref{thm:mainmax}.}
\proof{Proof of Theorem\ \ref{thm:mainmax} : } For $\epsilon > 0$, let $U^{\epsilon} \stackrel{\Delta}{=} (\frac{27}{2})^2 \gamma_0^3 \epsilon^{-2} M_1.$  It then follows from Lemma\ \ref{maxb1} and the definition of $\gamma_0$ that 
$$0 \leq \textsc{opt} - U^{\epsilon} \times \sup_{\tau \in {\mathcal T}} E[Z^1_{U^{\epsilon},\tau}] \leq \big( \frac{\gamma_0 M_1^3}{U^{\epsilon}} \big)^{\frac{1}{2}} = \frac{2}{27} \gamma_0^{-1} \epsilon M_1.$$
Combining with Lemma\ \ref{maxb2} (which implies that $M_1 \leq \frac{27}{4} \gamma_0 \textsc{opt}$), we conclude that 
\begin{equation}\label{altuberror1bb}
0 \leq \textsc{opt} - U^{\epsilon} \times \sup_{\tau \in {\mathcal T}} E[Z^1_{U^{\epsilon},\tau}] \leq \frac{\epsilon}{2} \textsc{opt}.
\end{equation}
Since $0 \leq Z^1_{U^{\epsilon},t} \leq 1$ for all $t$ (by definition of the truncated and normalized values), for any $\epsilon',\delta \in (0,1)$, we may use Algorithm ${\mathcal A}$ (as in Theorem\ \ref{thm:mainalgo2}) to compute an $\epsilon'$-approximation to $\sup_{\tau \in {\mathcal T}} E[Z^1_{U^{\epsilon},\tau}]$ in total computational time at most $(C + G + 1) \times ( 1  +  \log(\delta^{-1})) \times \exp(100 \epsilon'^{-2})\times  T^{8 \epsilon'^{-1}}$ with probability at least $1 - \delta$.  In light of (\ref{altuberror1bb}), it is natural to then return $U^{\epsilon}$ times this approximation, which will magnify the error multiplicatively by $U^{\epsilon}$.  Since by Lemma\ \ref{maxb2} we have $U^{\epsilon} \leq \overline{U^{\epsilon}} \stackrel{\Delta}{=} \frac{27^3}{16} \gamma_0^4 \epsilon^{-2} \textsc{opt}$, it follows that returning $U^{\epsilon}$ times this approximation magnifies the error multiplicatively by at most $\overline{U^{\epsilon}}.$  By setting $\epsilon' = \frac{\epsilon}{2} \times \textsc{opt} \times \frac{1}{\overline{U^{\epsilon}}} =  \epsilon^3 \times \frac{8}{27^3} \times \gamma_0^{-4}$, we compensate for this error magnification, and by an application of the triangle inequality derive an overall error of $\frac{\epsilon}{2} \times \textsc{opt} + \frac{\epsilon}{2} \times \textsc{opt} = \epsilon \times \textsc{opt}$.  Combining the above, and accounting for the extra time needed to truncate and normalize the rewards, completes the proof.  \qed
\endproof

\subsubsection{Proof of Theorem\ \ref{thm:mainmaxpolicy}.}
\proof{Proof of Theorem\ \ref{thm:mainmaxpolicy} : }
For any $U > 0$ and $\epsilon \in (0,1)$, we may use Algorithm ${\mathcal A}'$ (as in Theorem\  \ref{thm:mainpolicy1}) to implement a randomized stopping time $\tau_{U,\epsilon}$ s.t. $E[Z^1_{U,\tau_{U,\epsilon}}] \geq \sup_{\tau \in {\mathcal T}} E[Z^1_{U,\tau}] - \epsilon$.
Since $\tau_{U,\epsilon} \in {\mathcal T}$,  it also yields a stopping time for the original problem $\sup_{\tau \in {\mathcal T}} E[Z^1_{\tau}]$.  By construction, for all $t \geq 1$, w.p.1 $Z^1_t \geq U \times Z^1_{U,t}.$  Combining the above, it follows that for any $U > 0$ and $\epsilon \in (0,1)$, 
\begin{eqnarray*}
E[Z^1_{\tau_{U,\epsilon}}] &\geq& U \times E[Z^1_{U,\tau_{U,\epsilon}}]
\\&\geq& U \times \sup_{\tau \in {\mathcal T}} E[Z^1_{U,\tau}] - U \times \epsilon.
\end{eqnarray*}
Combining with Lemma\ \ref{maxb1} and the triangle inequality, we conclude that for all $U > 0$ and $\epsilon' \in (0,1)$, 
$$\textsc{opt} - E[Z^1_{\tau_{U,\epsilon'}}] \leq \big( \frac{\gamma_0 M_1^3}{U} \big)^{\frac{1}{2}} + U \times \epsilon'.$$
For $\epsilon > 0$, let $U^{\epsilon} \stackrel{\Delta}{=} (\frac{27}{2})^2 \gamma_0^3 \epsilon^{-2} M_1$ and $\epsilon' = \epsilon^3 \times \frac{8}{27^3} \times \gamma_0^{-4}.$
It then follows from an argument identical to that used in our proof of Theorem\ \ref{thm:mainmax} that for any $\epsilon \in (0,1)$, 
$$\textsc{opt} - E[Z^1_{\tau_{U^{\epsilon},\epsilon'}}] \leq \epsilon \times \textsc{opt}.$$
Combining with Theorem\ \ref{thm:mainpolicy1}, and accounting for the extra time needed to truncate and normalize the rewards, completes the proof.  \qed
\endproof

\subsubsection{An alternate bound not involving the s.c.v.}\label{altubsec} As noted in the introduction, the s.c.v. arises naturally in our analysis due to having to control the error induced by truncating the rewards (which introduces terms involving $M_1$ and $M_2$ due to an application of the Cauchy-Schwarz inequality in the proof of Lemma\ \ref{maxb1}), and our having to show that this error is not too large relative to the value of the optimal stopping problem (we show in Lemma\ \ref{maxb2} that this will hold as long as the s.c.v. is not too large).  However, in some settings it may be difficult to bound the ratio of $M_2$ and $M_1^2$, but easy to upper bound both $M_1$ and $M_2$ separately (e.g. in the case that it is difficult to derive tight lower bounds on $M_1$).  We now present an alternate analysis which yields a runtime depending on the product of $M_1$ and $M_2$, but not the s.c.v.
\begin{lemma}\label{altub1}
Suppose that only the assumption regarding the efficient simulator, Assumption\ \ref{assumptionB}, holds.  Then Algorithm ${\mathcal A},$ with properly chosen parameters and an additional truncation step, takes any $\epsilon,\delta \in (0,1)$ and achieves the following. In total computational time at most $(C + G + 1) \times ( 1  +  \log(\delta^{-1})) \times \exp\big(10^4 (M_1 M_2)^2 \epsilon^{-6}\big)\times  T^{64 M_1 M_2 \epsilon^{-3}}$, returns a random number $X$ satisfying $P\big( |X - \textsc{opt}| \leq \epsilon \big) \geq 1 - \delta$.
\end{lemma}
\proof{Proof} For $\epsilon > 0$, let $U^{\epsilon} \stackrel{\Delta}{=} 4 M_1 M_2 \epsilon^{-2}.$  By Lemma\ \ref{maxb1},
$$0 \leq \textsc{opt} - U^{\epsilon} \times \sup_{\tau \in {\mathcal T}} E[Z^1_{U^{\epsilon},\tau}] \leq (M_2)^{\frac{1}{2}} \times \big( \frac{M_1}{U^{\epsilon}} \big)^{\frac{1}{2}}.$$
Since by construction of $U^{\epsilon}$ it holds that $(M_2)^{\frac{1}{2}} \times \big( \frac{M_1}{U^{\epsilon}} \big)^{\frac{1}{2}} = \frac{\epsilon}{2}$, we conclude that
\begin{equation}\label{altuberror1}
0 \leq \textsc{opt} - U^{\epsilon} \times \sup_{\tau \in {\mathcal T}} E[Z^1_{U^{\epsilon},\tau}] \leq \frac{\epsilon}{2}.
\end{equation}
Since $0 \leq Z^1_{U^{\epsilon},t} \leq 1$ for all $t$ (by definition of the truncated and normalized values), for any $\epsilon',\delta \in (0,1)$, we may use Algorithm ${\mathcal A}$ (as in Theorem\ \ref{thm:mainalgo2}) to compute an $\epsilon'$-approximation to $\sup_{\tau \in {\mathcal T}} E[Z^1_{U^{\epsilon},\tau}]$ in total computational time at most $(C + G + 1) \times ( 1  +  \log(\delta^{-1})) \times \exp(100 \epsilon'^{-2})\times  T^{8 \epsilon'^{-1}}$ with probability at least $1 - \delta$.  In light of (\ref{altuberror1}), it is natural to then return $U^{\epsilon}$ times this approximation, which will magnify the error multiplicatively by $U^{\epsilon}$.  By setting $\epsilon' = \frac{\epsilon}{2 U^{\epsilon}} = \frac{\epsilon^3}{8 M_1 M_2}$, we compensate for this error magnification, and by an application of the triangle inequality derive an overall error of $\frac{\epsilon}{2} + \frac{\epsilon}{2} = \epsilon$.  Combining the above completes the proof.  \qed
\endproof
Let us note that in many settings, e.g. the case that $\lbrace Z_t, t \geq 1 \rbrace$ are i.i.d. exponentially distributed r.v.s, or the setting of the Bermudan Max Call analyzed in Section\ \ref{sec:maxcall}, the s.c.v. of the maximum reward can be bounded independent of the time horizon, while $M_1$ and $M_2$ themselves cannot be so bounded.  In such settings our bounds involving the s.c.v. will thus yield a much tighter analysis.

\section{Proof of Corollary\ \ref{cormaxcall}: Bermudan max call satisfies assumptions}\label{sec:ecmaxcall}
In this section, we provide a detailed analysis of our result regarding the pricing of the Bermudan max call, as stated in Corollary \ref{cormaxcall}. As we discussed before, the key is to verify Assumptions \ref{assumptionA2} and \ref{assumptionB}, especially Assumption \ref{assumptionA2} whose proof will be non-trivial.  Let ${\mathcal M} \stackrel{\Delta}{=} \max_{i = 1,..., D} \max_{t \in \{t_0, ... , t_T\}} e^{-r t} (G^i_t - \kappa)^+ + \kappa.$ 

\begin{lemma}[Bermudan max call satisfies Assumption \ref{assumptionA2}]\label{expA2}
 For any Bermudan max call instance satisfying $\varrho > 1.5 \sigma^2,$ 
 \[
 \frac{E[{\mathcal M}^2]}{E[{\mathcal M}]^2} \leq 10^4 e^{2 (\frac{\varrho + \frac{\sigma^2}{2}}{\sigma})^2} (1 + \frac{20 \sigma^2}{\varrho - 1.5 \sigma^2}).
 \]
\end{lemma}
\proof{Proof}
Note that
\begin{eqnarray*}
e^{-r t}(G^i_t - \kappa)^+ + \kappa &\geq& e^{-r t}(G^i_t - \kappa) + e^{- rt} \kappa
\\&=& e^{- r  t}G^i_t.
\end{eqnarray*}
It follows that 
\begin{align*}
    \frac{E[{\mathcal M}^2]}{E[{\mathcal M}]^2} & \ = \ 1 + \frac{Var({\mathcal M})}{E[{\mathcal M}]^2},
    \\ & \leq  \ 1 + \frac{Var\bigg( \max_{i = 1,..., D} \max_{t \in \{t_0, ... , t_T\}} e^{-r t} (G^i_t - \kappa)^+ + \kappa\bigg)}{E\bigg[ \max_{i = 1,..., D} \max_{t \in \{t_0, ... , t_T\}} e^{-r t} G^i_t \bigg]^2} 
    \\ & =  \ 1 + \frac{Var\bigg( \max_{i = 1,..., D} \max_{t \in \{t_0, ... , t_T\}} e^{-r t} (G^i_t - \kappa)^+ \bigg)}{E\bigg[ \max_{i = 1,..., D} \max_{t \in \{t_0, ... , t_T\}} e^{-r t} G^i_t \bigg]^2} 
    \\ & \leq  \ 1 +  \frac{E\bigg[ \bigg(\max_{i = 1,..., D} \max_{t \in \{t_0, ... , t_T\}} e^{-r t} (G^i_t - \kappa)^+ \bigg)^2\bigg]}{E\bigg[ \max_{i = 1,..., D} \max_{t \in \{t_0, ... , t_T\}} e^{-r t} G^i_t \bigg]^2} 
    \\ & \leq  \ 1 + \frac{E\bigg[ \bigg(\max_{i = 1,..., D} \max_{t \in \{t_0, ... , t_T\}} e^{-r t} G^i_t  \bigg)^2 \bigg]}{E\bigg[ \max_{i = 1,..., D} \max_{t \in \{t_0, ... , t_T\}} e^{-r t} G^i_t \bigg]^2}.
\end{align*}

Hence it suffices to upper bound the above expression.  Furthermore, as the above expression takes the same value for all values of $y_0$ (due to the multiplicative nature of geometric brownian motion), we may w.l.o.g. assume $y_0 = 1$.  We now bound the above ratio (where the maximum is over a discrete set of times in $[0,{\mathcal J}]$) in terms of the analogous ratio in which the suprema are over all times in $[0,{\mathcal J}]$.  Trivially, we may bound the numerator as follows: 
\[
E\bigg[ \bigg(\max_{i = 1,..., D} \max_{t \in \{t_0, ... , t_T\}} e^{-r t} G^i_t  \bigg)^2 \bigg] \leq E\bigg[ \bigg(\max_{i = 1,..., D} \sup_{t \in [0, \mathcal{J}]} e^{-r t} G^i_t  \bigg)^2 \bigg],
\]
where we replace the grids $\{t_0,...,t_T\}$ by the interval $[0, \mathcal{J}].$ We now provide a similar lower bound on the denominator.  Indeed, we show that 
\begin{equation}\label{ineq0}
    \frac{E\big[ \max_{i = 1,..., D} \max_{t \in \{t_0, ... , t_T\}} e^{-r t} G^i_t  \big]}{ E\big[ \max_{i = 1,..., D} \sup_{t \in [0, \mathcal{J}]} e^{-r t} G^i_t \big]} \ \geq\  \Phi^c\big(\frac{\varrho}{\sigma} + \frac{\sigma}{2}\big).
\end{equation}
We argue as follows. Let $\max_{t \in \{t_0,...,t_T\}} e^{-r t } G^i_t \stackrel{\Delta}{=} M_i'$ and $\sup_{t \in [0, \mathcal{J}]} e^{-r t } G^i_t \stackrel{\Delta}{=} M_i.$ Then $\{(M'_i, M_i), i = 1,...,D\}$ are $i.i.d.$ random vectors, where w.p.1 $M_i \geq M'_i$ for all $i$. Here we note that $M_i,M'_i$ have no connection to the constants $M_1,M_2$ used in the previous section.  By definition, process $\lbrace e^{-r t } G^i_t, t \in [0, {\mathcal J}] \rbrace$ can be explicitly represented as $\lbrace \exp\bigg((-\varrho- \frac{\sigma^2}{2}) t + \sigma W^i_t\bigg), t \in [0, {\mathcal J}] \rbrace$ namely a geometric B.m. with negative drift.  Here $\lbrace W^i_t, t \in [0, {\mathcal J}] \rbrace$ are independent standard B.m.  To prove bound \ref{ineq0}, we first show $M'_i$ is near $M_i$ in the following precise sense:
\begin{equation}\label{ineq1}
    P(M'_i  \geq z ) \geq \Phi^c\big(\frac{\varrho}{\sigma} + \frac{\sigma}{2}\big) P(M_i \geq z)\ \textrm{for all}\ z \geq 0\ \textrm{and all}\ i.
\end{equation}
As both discounted processes (continuous and discretized) start at 0 at time 0 (as we have assumed $y_0 = 0$), note that w.p.1 $M'_i \geq 1, M_i \geq 1$ for all $i$.  Thus it suffices to demonstrate (\ref{ineq1}) for $z > 1$.  Indeed, let $\tau_z \stackrel{\Delta}{=} \inf\lbrace t \in [0, \mathcal{J}]:  e^{-r t} G^i_t = z \rbrace$ denote the first time in $[0, \mathcal{J}]$ that the discounted payoff process equals value $z$ for $z > 1$ (here we are also implicitly using the continuity of the relevant processes).  We let $\tau_z = \mathcal{J}$ if the event does not happen. Then we have that for all $i$,
\begin{align*}
     & P(M'_i \geq z) \ = \ E[I(M'_i \geq z)]
    \\& \ = \  E[I(M'_i \geq z)\times I(\tau_z < \mathcal{J})] + E[I(M'_i \geq z)\times I(\tau_z = \mathcal{J})]
    \\& \ \geq \   E[I(e^{-r \lceil \tau_z \rceil} G^i_{\lceil \tau_z \rceil} \geq z)\times I(\tau_z < \mathcal{J})]  + E[I(M'_i \geq z)\times I(\tau_z = \mathcal{J})]
    \\&\ \ \ \ \ \ \ \ \ \textrm{since $\lceil \tau_z \rceil \in \{t_0,..., t_{T}\}$ by our assumption that $T$ is divisible by ${\mathcal J}$}
    \\& \ \geq \   E[I(e^{-r \lceil \tau_z\rceil} G^i_{\lceil \tau_z\rceil} \geq z)\times I(\tau_z < \mathcal{J})]  
    \\& \ = \   E[I(e^{-r \lceil \tau_z \rceil} G^i_{\lceil \tau_z \rceil} \geq e^{-r  \tau_z} G^i_{\tau_z})\times I(\tau_z < \mathcal{J})]
    \\& \ = \   E\bigg[I\bigg((-\varrho- \frac{\sigma^2}{2}) \lceil\tau_z\rceil + \sigma W^i_{\lceil\tau_z\rceil}\ \geq\ (-\varrho- \frac{\sigma^2}{2}) \tau_z + \sigma W^i_{\tau_z} \bigg)\times I(\tau_z < \mathcal{J})\bigg]
     \\& \ = \   E\bigg[I\bigg(W^i_{\lceil\tau_z\rceil} - \ W^i_{\tau_z} \geq (\frac{\varrho}{\sigma} + \frac{\sigma}{2}) (\lceil \tau_z \rceil - \tau_z) \bigg)\times I(\tau_z < \mathcal{J})\bigg]
     \\& \ \geq \ \inf_{t \in (0,1)} P\big( W^i_t \geq (\frac{\varrho}{\sigma} + \frac{\sigma}{2}) t \big)  P(\tau_z < \mathcal{J}) 
\\& \ \ \ \ \ \ \ \ \textrm{by the strong Markov property, since w.p.1\ $\lceil \tau_z \rceil - \tau_z \in (0,1)$}
\\& \ =\ \inf_{t \in (0,1)} \Phi^c\big( (\frac{\varrho}{\sigma} + \frac{\sigma}{2}) t^{\frac{1}{2}} \big)  P(\tau_z < \mathcal{J})
\\& \ =\ \Phi^c\big( \frac{\varrho}{\sigma} + \frac{\sigma}{2} \big)  P(M_i \geq z),
\end{align*}
completing the proof of (\ref{ineq1}).  With (\ref{ineq1}) in hand, we proceed to prove bound (\ref{ineq0}).
\begin{align*}
  & E\big[ \max_{i = 1,..., D} \max_{t \in \{t_0, ... , t_T\}} e^{-r t} G^i_t  \big]  \ =\ E\big[ \max_{i = 1,..., D} M'_i  \big] 
   \\& \ = \ \int_{z = 0}^\infty  P\bigg(\max_{i = 1,..., D} M'_i  \geq z\bigg) d  z
   \\& \ = \ \int_{z = 0}^\infty  E\bigg[I\bigg(\max_{i = 1,..., D} M'_i  \geq z\bigg)\bigg] d  z
   \\& \ = \ \int_{z = 0}^\infty  \sum_{j = 1}^D E\bigg[I\bigg(\max_{i = 1,..., D} M'_i  \geq z\bigg)\times I\bigg(M_1,...,M_{j - 1} < z, M_j \geq z\bigg)\bigg] d  z 
   \\& \ \ \ \ \ \ \textrm{since} \max_{i=1,\ldots,D} M'_i \geq z \rightarrow \max_{i=1,\ldots,D} M_i \geq z, \textrm{implying some}\ M_j\ \textrm{is the first to exceed}\ z
   \\& \ \geq \ \int_{z = 0}^\infty  \sum_{j = 1}^D E\bigg[I\bigg(M'_j \geq z\bigg)\times I\bigg(M_1,...,M_{j - 1} < z, M_j \geq z\bigg)\bigg] d  z 
   \\& \ = \ \int_{z = 0}^\infty  \sum_{j = 1}^D E\bigg[I\bigg(M'_j \geq z\bigg)\times I\bigg( M_j \geq z\bigg)\bigg] \times E\bigg[I\bigg(M_1,...,M_{j - 1} < z\bigg)\bigg] d  z 
   \\& \ \ \ \ \ \ \ \textrm{by independence}
    \\& \ = \ \sum_{j = 1}^D \int_{z = 0}^\infty   E\bigg[I\bigg(M'_j \geq z\bigg)\times I\bigg( M_j \geq z\bigg)\bigg] \times E\bigg[I\bigg(M_1,...,M_{j - 1} < z\bigg)\bigg] d  z
    \\& \ = \ \sum_{j = 1}^D \int_{z = 0}^\infty   E\bigg[I\bigg(M'_j \geq z\bigg)\bigg] \times E\bigg[I\bigg(M_1,...,M_{j - 1} < z\bigg)\bigg]  d  z
   \\& \ \ \ \ \ \ \ \textrm{since $M'_j \geq z$ implies $M_j \geq z$}
   \\& \ \geq \ \sum_{j = 1}^D \int_{z = 0}^\infty \Phi^c\big( \frac{\varrho}{\sigma} + \frac{\sigma}{2} \big) E\bigg[I\bigg(M_j \geq z\bigg)\bigg]\times E\bigg[I\bigg(M_1,...,M_{j - 1} < z\bigg)\bigg]  d  z
 \ \ \ \ \ \ \ \ \ \ \ \ \ \  \textrm{by (\ref{ineq1})}
   \\& \ = \ \Phi^c\big( \frac{\varrho}{\sigma} + \frac{\sigma}{2} \big) \int_{z = 0}^\infty  \sum_{j = 1}^D  E\bigg[I\bigg(M_1,..., M_{j-1} < z, M_j \geq z\bigg)\bigg]  d  z
    \\& \ = \ \Phi^c\big( \frac{\varrho}{\sigma} + \frac{\sigma}{2} \big) \int_{z = 0}^\infty  E\bigg[I\bigg(\max_{1 \leq i \leq D} M_i \geq z\bigg)\bigg]  d  z
    \\& \ = \ \Phi^c\big( \frac{\varrho}{\sigma} + \frac{\sigma}{2} \big) E[\max_{1 \leq i \leq D} M_i].
\end{align*}
Thus we complete the proof of bound \ref{ineq0}.  Recall again that $M_i = \sup_{t \in [0, \mathcal{J}]} e^{-r t} G^i_t = \sup_{t \in [0, \mathcal{J}]} \exp\big(-(\varrho + \sigma^2/2) t + \sigma W^i_t \big)$, and let us define $M^* \stackrel{\Delta}{=} \max_{i = 1,.., D} M_i.$  Then combining both the upper bound and the lower bound (\ref{ineq1}), we have that 
\begin{equation}\label{finalineeq}
    \frac{E[{\mathcal M}^2]}{E[{\mathcal M}]^2} \leq \ 1 + \bigg(\Phi^c\big( \frac{\varrho}{\sigma} + \frac{\sigma}{2} \big)\bigg)^{-2} \frac{E[ (M^*)^2]}{E[M^*]^2}.
\end{equation}
In other words, it suffices to upper bound the CV of $\max_{i = 1,..., D} \sup_{t \in [0, \mathcal{J}]} e^{-r t} G^i_t,$ where the inner maximum is taken over the continuous interval $[0, \mathcal{J}].$   We now complete the proof of Lemma\ \ref{expA2}.  We consider two cases (where the first case is essentially trivial).  Case I: $D = 1.$ In this case, we have the following bounds:
\[
E\bigg[ \max_{i = 1,..., D} \sup_{t \in [0, \mathcal{J}]} e^{-r t} G^i_t \bigg] \ \geq\ 1\ \ \ \textrm{trivially by considering $t = 0$},
\]
and 
\begin{align*}
 & E\bigg[ \bigg(\max_{i = 1,..., D} \sup_{t \in [0, \mathcal{J}]} e^{-r t} G^i_t\bigg)^2 \bigg] 
\\ & \ \ =\ \ E\bigg[ \bigg( \sup_{t \in [0, \mathcal{J}]} e^{-r t} G^i_t\bigg)^2 \bigg] 
 \\& \ \ \leq \ \ E\bigg[\sup_{t \in [0, \infty)}  \exp\big(-2(\varrho + \sigma^2/2) t + 2\sigma W^i_t \big)\bigg]
 \\& \ \ = \ \ \frac{2\sigma^2}{\varrho - 1.5 \sigma^2}.
\end{align*}
Thus when $D = 1$, $E[(M^*)^2]/E[M^*]^2 \leq \frac{2\sigma^2}{\varrho - 1.5 \sigma^2}.$ 

\indent Next consider the primary case, namely $D \geq 2.$ Let $y^* \stackrel{\Delta}{=} \sup \{y: P(M_i \geq y ) = \frac{1}{D}\}$ be the $(1 - \frac{1}{D})$-quantile of the distribution of each $M_i$ (identical for each $i$). We next prove that $M^*$ is concentrated around $y^*$ in an appropriate multiplicative sense. More precisely, we prove that
\begin{align}
    P(M^* \geq  y^*  e^{\Delta}) & \ \leq\ 10 e^{- \frac{\varrho + \sigma^2/2}{\sigma^2} \Delta }\ \textrm{for all}\ \Delta \geq 0; \label{ineq2}
    \\P(M^* \geq  y^*) & \ \geq\  1 - \frac{1}{e}. \label{ineq3}
\end{align}  
By definition, it is easy to show the following relation between the c.d.f. of $M^*$ and the c.d.f. of $M_1:$
\[
P(M^* \geq y) = 1 - \bigg(1 - P(M_1 > y)\bigg)^D.
\]
(\ref{ineq3}) follows almost immediately, as $P(M^* \geq y^*) = 1 - (1 - \frac{1}{D})^D$, which can be shown by a straightforward calculus exercise to be at least $1 - \frac{1}{e}$ for all $D \geq 1$.  This completes the proof of (\ref{ineq3}).  
\indent We now prove (\ref{ineq2}).  We first relate $P(M_i \geq  y^*  e^{\Delta})$ for $\Delta > 0$ to $P(M_i \geq  y^*).$  Since $M_i$ is the maximum of an exponentiated B.m. with drift, we have the closed-form formula for all $i$ and $y > 0$ (see e.g. \cite{boukai1990explicit})
\[
P(M_i > y) = e^{- 2 \log y \frac{\varrho + \sigma^2/2}{\sigma^2}}\Phi^c\bigg(\frac{\log y - (\varrho + \sigma^2/2) \mathcal{J}}{\sigma \sqrt{\mathcal{J}}}\bigg) + \Phi^c\bigg(\frac{\log y + (\varrho + \sigma^2/2) \mathcal{J}}{\sigma \sqrt{\mathcal{J}}}\bigg).
\]
We now bound $P(M_i > y^* e^{\Delta})$ in terms of $P(M_i > y^*),$ proving 
\begin{equation}\label{interim11}
P(M_i > y^*e^\Delta ) \leq P(M_i > y^*)\times 10 e^{- 2 \Delta  \frac{\varrho + \sigma^2/2}{\sigma^2}  }.
\end{equation}

  For any $\Delta \geq 0$,
\begin{align*}
  &   P(M_i > y^* e^{\Delta})
  \\& \ \ \ \ =\  e^{- 2 (\log y^* + \Delta) \frac{\varrho + \sigma^2/2}{\sigma^2}}\Phi^c\bigg(\frac{\log y^* + \Delta - (\varrho + \sigma^2/2) \mathcal{J}}{\sigma \sqrt{\mathcal{J}}}\bigg) 
  \\& \ \ \ \  \ \ \ + \ \ \ \Phi^c\bigg(\frac{\log y^* + \Delta + (\varrho + \sigma^2/2) \mathcal{J}}{\sigma \sqrt{\mathcal{J}}}\bigg)
  \\& \ \ \ \ =\  e^{- 2 \Delta  \frac{\varrho + \sigma^2/2}{\sigma^2}  }\times e^{- 2 \log y^*  \frac{\varrho + \sigma^2/2}{\sigma^2}}\Phi^c\bigg(\frac{\log y^* + \Delta - (\varrho + \sigma^2/2) \mathcal{J}}{\sigma \sqrt{\mathcal{J}}}\bigg)
  \\& \ \ \ \  \ \ \ + \ \ \ \Phi^c\bigg(\frac{\log y^* + \Delta + (\varrho + \sigma^2/2) \mathcal{J}}{\sigma \sqrt{\mathcal{J}}}\bigg)
  \\& \ \ \ \ \leq \  e^{- 2 \Delta  \frac{\varrho + \sigma^2/2}{\sigma^2}  }\times e^{- 2 \log y^*  \frac{\varrho + \sigma^2/2}{\sigma^2}}\Phi^c\bigg(\frac{\log y^* - (\varrho + \sigma^2/2) \mathcal{J}}{\sigma \sqrt{\mathcal{J}}}\bigg)
  \\& \ \ \ \  \ \ \ + \ \ \ \Phi^c\bigg(\frac{\log y^* + \Delta + (\varrho + \sigma^2/2) \mathcal{J}}{\sigma \sqrt{\mathcal{J}}}\bigg).
\end{align*}
Let $u \stackrel{\Delta}{=} \frac{\log y^* + (\varrho + \sigma^2/2) \mathcal{J}}{\sigma \sqrt{\mathcal{J}}}.$  Then the above becomes 
\begin{align*}
    & P(M_i > y^* e^{\Delta}) \ \ \ \ \leq \  e^{- 2 \Delta  \frac{\varrho + \sigma^2/2}{\sigma^2}  }\times e^{- 2 \log y^*  \frac{\varrho + \sigma^2/2}{\sigma^2}}\Phi^c\bigg(\frac{\log y^* - (\varrho + \sigma^2/2) \mathcal{J}}{\sigma \sqrt{\mathcal{J}}}\bigg)
  \\& \ \ \ \  \ \ \ + \ \ \ \Phi^c\bigg(u + \frac{\Delta}{\sigma\sqrt{\mathcal{J}}}\bigg)
\end{align*}
Noting that the above equalities (evaluated at $\Delta = 0$) imply 
$$ P(M_i > y^*) = e^{- 2 \log y^*  \frac{\varrho + \sigma^2/2}{\sigma^2}}\Phi^c\bigg(\frac{\log y^* - (\varrho + \sigma^2/2) \mathcal{J}}{\sigma \sqrt{\mathcal{J}}}\bigg)
+ \Phi^c(u),$$
we conclude that to derive the desired bound \ref{interim11}, 
it suffices to prove that
\begin{equation}\label{eqtail}
    \Phi^c(u + \frac{\Delta}{\sigma\sqrt{\mathcal{J}}})\ \ \leq\ \ 10 e^{- 2 \Delta  \frac{\varrho + \sigma^2/2}{\sigma^2}  }\times \Phi^c(u).
\end{equation}
We proceed by a case analysis.  For $u > 2$, first note that (by the definition of $u$ and fact that $y^* > 1$) one has $u \geq \frac{\varrho + \frac{\sigma^2}{2}}{\sigma} \sqrt{{\mathcal J}}$, and thus $\frac{\Delta}{\sigma \sqrt{{\mathcal J}}} \geq \frac{(\varrho + \frac{\sigma^2}{2}) \Delta}{\sigma^2 u}$.  It then follows from Lemma \ref{lmnormaltail} that

\begin{align*}
    &\ \frac{\Phi^c(u + \frac{\Delta}{\sigma\sqrt{\mathcal{J}}})}{\Phi^c(u)} \ = R(u, \frac{\Delta}{\sigma\sqrt{\mathcal{J}}}) 
   \\& \ \leq R\bigg(u, \frac{(\varrho + \sigma^2/2)\Delta}{\sigma^2 u}\bigg).
\end{align*}
We now use Lemma \ref{lmnormaltail} to prove suitable bounds on $R$.
\begin{align*}
    R\bigg(u, \frac{(\varrho + \sigma^2/2)\Delta}{\sigma^2 u}\bigg) &\ = \ \frac{\Phi^c(u + \frac{(\varrho + \sigma^2/2)\Delta}{\sigma^2 u})}{\Phi^c(u)}
    \\&\ \leq \ \ \frac{\frac{1}{u + \frac{(\varrho + \sigma^2/2)\Delta}{\sigma^2 u}}\exp(-\frac{(u + \frac{(\varrho + \sigma^2/2)\Delta}{\sigma^2 u})^2}{2})}{\bigg(\frac{1}{u} - \frac{1}{u^3}\bigg) e^{-\frac{u^2}{2}}}
     \\&\ \leq \ \ \frac{u^2}{u^2-1}\exp\bigg(-\Delta\frac{\varrho + \sigma^2/2}{\sigma^2}\bigg)
     \\&\ \leq \ \ 2 e^{-\Delta\frac{\varrho + \sigma^2/2}{\sigma^2}}. \ \ \ \ \ \ \ \ \ \ \textrm{(because $u > 2$)}
\end{align*}
For $u \leq 2,$ by Lemma \ref{lmnormaltail} we directly bound
\begin{align*}
    R\bigg(u, \frac{(\varrho + \sigma^2/2)\Delta}{\sigma^2 u}\bigg) &\ \leq\  R\bigg(0, \frac{(\varrho + \sigma^2/2)\Delta}{2 \sigma^2 }\bigg)\ \ \ \ \ \ \  \ \ \ \ \textrm{(because $0 \leq u \leq 2$)}
    \\& \ \ = \frac{\Phi^c(\frac{(\varrho + \sigma^2/2)\Delta}{2 \sigma^2 })}{\Phi^c(0)}
    \\& \ \ =  \ 2\Phi^c(\frac{(\varrho + \sigma^2/2)\Delta}{2 \sigma^2 }) 
   \\& \ \ \leq 10 e^{-\frac{(\varrho + \sigma^2/2)\Delta}{ \sigma^2 }}, 
\end{align*}
the final inequality following from Lemma\ \ref{lmnormaltail} and some straightfoward algebra and a case analysis the details of which we omit.  Combining the above completes the proof of Equation \ref{eqtail}, which in turn completes the proof of (\ref{interim11}).  To complete the proof of (\ref{ineq2}), we note that a simple union bound implies $P(M^* \geq  y^*  e^{\Delta}) \leq D P(M_1 \geq y^* e^{\Delta})$.  Combining with (\ref{interim11}) and the definition of $y^*$ then completes the proof of (\ref{ineq2}).
\indent With (\ref{ineq2}) and (\ref{ineq3}) in hand, we are able to bound $E[M^*]$ and $E[(M^*)^2].$ On the one hand, 
\begin{align*}
 E[M^*] & \ \geq \ E[y^* I(M^* \geq y^*)]
\\& \ \ \geq (1 - \frac{1}{e}) y^*.
\end{align*}
On the other hand, noting that (\ref{ineq2}) implies that 
$$P\big( \frac{M^*}{y^*} \geq z \big) \leq 10 z^{- \frac{\varrho + \frac{\sigma^2}{2}}{\sigma^2}}\ \textrm{for all}\ z \geq 1,$$ 
we have
\begin{align*}
    E[(M^*)^2] & \ = \ \int_{ 0}^\infty P\bigg((M^*)^2 > t\bigg) d t 
    \\& \ \ =  \ \int_{ 0}^{(y^*)^2} P\bigg((M^*)^2 > t\bigg) d t  + \int_{ (y^*)^2}^{\infty} P\bigg((M^*)^2 > t\bigg) d t
    \\& \ \ =  \ (y^*)^2  + \int_{1}^{\infty} P\bigg((M^*)^2 > u (y^*)^2  \bigg) (y^*)^2 du
    \\& \ \ \leq  \ (y^*)^2  + (y^*)^2 \int_{1}^{\infty} P\bigg( M^* > \sqrt{u} y^* \bigg) du
    \\& \ \ \leq  \ (y^*)^2  + (y^*)^2 \int_{1}^{\infty} 10 u^{- \frac{\varrho + \frac{\sigma^2}{2}}{2 \sigma^2}} du
    \\& \ \ = \ (y^*)^2  + 10 (y^*)^2 \frac{2 \sigma^2}{\varrho - 1.5 \sigma^2}
    \\& \ \ \leq (y^*)^2 \big(1 +  20 \frac{\sigma^2}{\varrho - 1.5 \sigma^2} \big).
\end{align*}
Combining the above, we conclude that in all cases 
$$E[(M^*)^2]/E[M^*]^2 \leq (1 - \frac{1}{e})^{-2} (1 + \frac{20 \sigma^2}{\varrho - 1.5 \sigma^2}).$$
Combining with Equation (\ref{finalineeq}), we conclude that 
$$\frac{E[{\mathcal M}^2]}{E[{\mathcal M}]^2} \ \leq \ 1 + \bigg(\Phi^c\big( \frac{\varrho}{\sigma} + \frac{\sigma}{2} \big)\bigg)^{-2} (1 - \frac{1}{e})^{-2} (1 + \frac{20 \sigma^2}{\varrho - 1.5 \sigma^2}).$$
As it follows from Lemma\ \ref{lmnormaltail} and a straightforward case analysis (the details of which we omit) that $\Phi^c(x) \geq .02 e^{-x^2}$ for all $x \geq 0$, and combining with some straightforward algebra then completes the proof.
\qed
\endproof

\begin{lemma}[Bermudan max call satisfies Assumption \ref{assumptionB}]\label{expB}
For any Bermudan max call instance, the runtime to generate one (conditioned) trajectory and to compute the associated stopping rewards along the trajectory is bounded by $10 T D$ (for all $t$ and $\gamma_{[t]}$), under our computational model. Namely, $C + G \leq 10 T D$.
\end{lemma}
\proof{Proof}
In this example, $Y_{[T]}$ is distributed as $D$ independent discretized geometric B.m.  Thus conditional on any $t$ and trajectory $\gamma_{[t]}$, the remainder of the trajectory will be (conditionally) distributed as $D$ independent geometric B.m. (sampled discretely at times $t+1,\ldots,T$), initialized at the values consistent with $\gamma_{t}$ (i.e. the conditioned sample path at time $t$ and considering the Markovian structure of geometric B.m.).  Thus these trajectories can be simulated in a straightfoward iterative manner by generating $(T - t) \times D$ independent standard normals, and iteratively scaling, exponentiating, and multiplying with the $D$ respective processes to construct the remaining trajectories.  As the reward functions are simple functions of the maximum of these processes (at any given time), they too can be evaluated (for any given $Y_{[T]}$) efficiently in our computational model.  It may be easily verified that these straightforward steps can be implemented in time bounded by $10 T D$ in our computational model, and we omit the details.
\qed
\endproof
With Lemmas \ref{expA2} and \ref{expB} proven, we complete the proof of Corollary \ref{cormaxcall}.
\proof{Proof of Corollary \ref{cormaxcall}}
By applying a constant shift $\kappa$ to the stopping reward, namely $\max_{i = 1,..., D} e^{-r t} (G^i_t - \kappa)^+ + \kappa,$ the optimal value increases by $\kappa$ and the optimal stopping time remains unchanged. It is thus equivalent for us to solve a stopping problem associated with this new payoff. Let $\textsc{opt}^s = \textsc{opt} + \kappa$ denote the optimal value under the new payoff. The squared CV of the path-wise maximum of this new payoff is $\frac{E[M^2]}{E[M]^2}$ where $M = \max_{t \in \{t_0,...,t_T\}}\max_{i = 1,..., D} e^{-r t} (G^i_t - \kappa)^+ + \kappa.$ By Lemma \ref{expA2}, it is bounded by a constant 
\[
\gamma_0 = 10^4 e^{2 (\frac{\varrho + \frac{\sigma^2}{2}}{\sigma})^2} (1 + \frac{20 \sigma^2}{\varrho - 1.5 \sigma^2}),
\]
satisfying Assumption \ref{assumptionA2}. By Lemma \ref{expB}, unit sampling and evaluation cost is bounded by $C + G \leq 10 T D,$ satisfying Assumption \ref{assumptionB}. By Theorem \ref{thm:mainmax}, Algorithm $\mathbf{A}$ returns a random number satisfying 
\[
P\bigg(\frac{|X - \textsc{opt}^s|}{\textsc{opt}^s} \leq \epsilon \bigg) \ \geq \ 1 - \delta.
\]
Since $\textsc{opt}^s = \textsc{opt} + \kappa,$ it follows that $X - \kappa$ satisfies the precision requirement of Corollary \ref{cormaxcall}. Combining with Theorem\ \ref{thm:mainmax} and some straightforward algebra then completes the proof.
\qed
\endproof

\section{Proof that assuming non-negativity is w.l.o.g.}\label{sec:ecnonneg}
In this section we provide a short proof of the observation that it is w.l.o.g. to assume that the reward process is non-negative.  We note that a similar observation was made recently in 
\citet{baas2023sampling}, which applies our methodology to the computation of the Gittins index.
\begin{observation}
Every OS problem in which $Z_t$ is integrable for all $t$ can be easily transformed into an essentially equivalent problem in which $Z_t \geq 0$ for all $t$.
\end{observation}
\proof{Proof}
Note that one can always rewrite $Z_t$ as $Z_t = E[\min_s Z_s | {\mathcal F}_t] + (Z_t - E[\min_s Z_s | {\mathcal F}_t])$.  Applying the optional stopping theorem, we deduce that $\sup_{\tau \in {\mathcal T}} E[Z_{\tau}] = E[\min_t Z_t] + \sup_{\tau \in {\mathcal T}} E\big[ (Z_{\tau} -  E[\min_s Z_s | {\mathcal F}_{\tau}]),$ where the process $\lbrace Z_t - E[\min_s Z_s | {\mathcal F}_t], t = 1,\ldots,T \rbrace$ is non-negative. 
\qed
\endproof

\end{document}